\documentclass[letterpaper, 10pt, journal, twocolumn]{ieeetran}

\usepackage{epsfig}

\usepackage{amsmath,graphicx,amsfonts,amssymb,epsfig,subfig,mathrsfs,mathtools}
\usepackage{ntheorem}
\usepackage{cuted}
\usepackage{arydshln}
\usepackage{blkarray}

\usepackage{color}
\usepackage{multirow}
\usepackage{multicol}
\usepackage{lipsum}
\usepackage{rotating}
\usepackage{graphicx}
\usepackage{psfrag}
\usepackage{algorithm}

\usepackage[utf8]{inputenc}
\usepackage{xspace}

\usepackage{algorithmicx}
\usepackage{algpseudocode}
\usepackage{cite}
\usepackage{tikz}
\usepackage{cancel}
\usepackage{textcomp}
 \definecolor{lin}{RGB}{240,0,0}
 \definecolor{paleblue}{RGB}{0,9,255}
\usepackage[colorlinks=true,linkcolor=lin,citecolor=paleblue,pdfa]{hyperref}
\usepackage{comment}
\includecomment{comment}

\newcommand{\map}[3]{#1: #2 \rightarrow #3}
\newcommand{\setdef}[2]{\{#1 \; | \; #2\}}
\newcommand{\st}{\ensuremath{\operatorname{s.t.}}}
\newcommand{\real}{\ensuremath{\mathbb{R}}}
\newcommand{\prob}{\ensuremath{\mathbb{P}}}
\newcommand{\realnonnegative}{\ensuremath{\mathbb{R}}_{\ge 0}}

\newcommand{\until}[1]{\{1,\dots, #1\}}
\newcommand{\subscr}[2]{#1_{\textup{#2}}}
\newcommand{\supscr}[2]{#1^{\textup{#2}}}
\newcommand{\vect}[1]{\boldsymbol{#1}}
\newcommand{\vectorones}[1]{\vect{1}_{#1}}
\newcommand{\vectorzeros}[1]{\vect{0}_{#1}}
\newcommand{\Norm}[1]{\|#1\|}
\newcommand{\trans}[1]{{#1}^\top}
\newcommand{\obj}{\ensuremath{\operatorname{obj}}}
\newcommand{\sgn}{\ensuremath{\operatorname{sgn}}}
\newcommand{\conv}{\ensuremath{\operatorname{conv}}}
\newcommand{\dom}{\ensuremath{\operatorname{dom}}}
\newcommand{\proj}{\operatorname{proj}}
\newcommand{\untilinterval}[2]{\{{#1},\dots, {#2}\}}

\newcommand{\reviseone}[1]{{{#1}}}
\newcommand{\revisetwo}[1]{{{#1}}}

\newcommand{\ODAA}{\textsc{OnDA} Algorithm\xspace}
\newcommand{\ODAAfull}{Online Data Assimilation Algorithm\xspace}
\newcommand{\CGA}{\textsc{Certificate Generation Algorithm}\xspace}
\newcommand{\ICA}{\textsc{Incremental Covering Algorithm}\xspace}
\newcommand{\CGAaco}{\textsc{C-Gen} Algorithm\xspace}
\newcommand{\ICAaco}{\textsc{I-Cover} Algorithm\xspace}

\DeclareMathOperator*{\argmax}{argmax}
\DeclareMathOperator*{\argmin}{argmin}

\makeatletter
\newtheoremstyle{breaknote}
  {\item{\theorem@headerfont
          ##1\ ##2\theorem@separator}\hskip\labelsep\relax}
  {\item{\theorem@headerfont
          ##1\ ##2\ (##3)\theorem@separator}\hskip\labelsep\relax}
\makeatother
\theoremstyle{breaknote}
\newtheorem{assumption}{Assumption}[section]
\newtheorem{theorem}{Theorem}[section]

\newtheorem{lemma}{Lemma}[section]

\newtheorem{remark}{Remark}[section]

\allowdisplaybreaks

\title{\reviseone{Data assimilation and online optimization} with performance guarantees}

\author{Dan Li$^{1}$ and Sonia Mart{\'\i}nez$^{1}$
\thanks{*This research
was developed with funding from the DARPA (Lagrange) contract
N660011824027. The views, opinions and/or findings expressed are those
of the author and should not be interpreted as representing the
official views or policies of the Department of Defense or the
U.S. Government.}
\thanks{$^{1}$D. Li and S. Mart{\'\i}nez are with the Department of Mechanical and Aerospace Engineering, University of California San Diego, La Jolla, CA 92092, USA
        {\tt\small lidan@ucsd.edu; soniamd@ucsd.edu}}
}

\begin{document}

\maketitle

\begin{abstract}
  This paper considers a class of real-time stochastic optimization
  problems dependent on an unknown probability distribution.  In the
  considered scenario, data is streaming frequently while trying to
  reach a decision. Thus, we aim to devise a procedure that
  incorporates samples (data) of the distribution
  sequentially and adjusts decisions accordingly. We approach this
  problem in a distributionally robust optimization framework and
  propose a novel \ODAAfull (\ODAA) for this purpose. This algorithm
  guarantees out-of-sample performance of decisions with high
  probability, and gradually improves the quality of the decisions by
  incorporating the streaming data. We show that the \ODAA \,
  converges under a sufficiently slow data streaming rate, and provide
  a criteria for its termination after certain number of data have
  been collected. Simulations illustrate the results.
\end{abstract}

\section{Introduction}
Online data assimilation can benefit many applications that require
real-time decision making under uncertainty, such as optimal target
tracking, sequential planning problems, and robust quality control. In
these problems, uncertainty is often represented by a multivariate
random variable that has an unknown distribution.  To quantify
uncertainty and make reliable decisions, one often needs to gather a
large number of samples in advance. Such requirement, however, is hard
to achieve under scenarios where acquiring samples is expensive, or
when real-time decisions must be made. Alternatively, recent methods
such as distributionally robust optimization (DRO) have attracted
recent attention due to their capability to provide out-of-sample
performance guarantees with a finite number of samples. However, when
the data is collected over time, it remains unclear what the best
procedure is to assimilate the data in the ongoing optimization
process.  Motivated by this, this work studies how to incorporate
finitely streaming data into a DRO problem, while guaranteeing
out-of-sample performance via the generation of time-varying
certificates.

\textit{Literature Review.} Optimization under uncertainty is a vast
research area, and as such, available methods include stochastic
optimization~\cite{AS-DD-AR:14} and robust
optimization~\cite{AB-LEG-AN:09}. Recently, data-driven
distributionally robust optimization has regained popularity thanks to
its out-of-sample performance guarantees, see e.g.~\cite{PME-DK:17,
  RG-AJK:16} and~\cite{AC-JC:17-allerton,AC-JC:17-tac}, for a
distributed algorithm counterpart, and references therein. In this
setup, one defines a set of distributions or \textit{ambiguity set},
which contains the true distribution of the data-generating system
with high probability. Then, the out-of-sample performance of the
data-driven decision is obtained as the worst-case optimization over
the ambiguity set.
An attractive way of designing these sets is to consider a ball in the
space of probability distributions centered at a reference or
most-likely distribution constructed from the available data. In the
space of distributions, the popular distance metric is the Prokhorov
metric~\cite{EE-GI:06}, $\phi$-divergence~\cite{RJ-YG:16} and the
Wasserstein distance~\cite{PME-DK:17}. In particular, the
work~\cite{PME-DK:17} \reviseone{presents a tractable reformulation
  of DRO via Wasserstein balls, and is extended
  in~\cite{AC-JC:17-allerton} to a distributed setting.}
However, the available \reviseone{problem reformulation}
in~\cite{PME-DK:17} and \reviseone{the distributed algorithm}
in~\cite{AC-JC:17-allerton} do not consider the update of the decision
over time \reviseone{and streaming data}, which is the focus of
this work.
\reviseone{To design a tractable algorithm incorporating streaming data,}
our work connects to various
convex optimization methods~\cite{SB-LV:04,DPB:15} such as the Frank-Wolfe
(FW) Algorithm (e.g., conditional gradient algorithm), the Subgradient
Algorithm, and their variants, see e.g.~\cite{PW:70,CH:74,SLJ-MJ:15}
and references therein. Our emphasis on the convergence of the
data-driven decision obtained through a sequence of optimization
problems contrasts with typical algorithms developed for single
(non-updated) problems.

\textit{Statement of Contributions.} In this paper, we propose a new
\ODAAfull (\ODAA) to solve decision-making problems subject to
uncertainty. The distribution of the uncertainty is unknown and the
algorithm adjusts decisions based on realizations of the stochastic
variable sequentially revealed over time.  The new algorithm addresses
four challenges: 1) the evaluation of the out-of-sample performance of
every possible online decision; 2) the adaptation to online,
increasingly-larger data sets to reach a decision with
performance guarantees with increasingly higher probabilities; 3) the
availability of an online decision vector with performance guarantees
at any time; 4) the capability of handling sufficiently large
streaming data sets.

To address 1), we start from a DRO problem setting. This leads to a
worst-case optimization over an ambiguity set or neighborhood of the
empirical distribution constructed from a data set.  To solve this
intractable problem, we reformulate it into an equivalent convex
optimization over a simplex. This enables us to explore the simplex
vertex set to find a certificate (a value bounding the cost) of a
given decision with certain confidence.  When the data is streaming,
we consider a sequence of DRO problems and their equivalent convex
reformulations employing increasingly larger data sets. Thus, as the
data streams, the associated problems are defined over simplices of
increasingly larger dimension.  The similarities of the feasible sets
allow us to assimilate the data via specialized Frank-Wolfe Algorithm
variants, thus solving 2) via a \CGA (\CGAaco) described in
Section~\ref{sec:CG}. Further, to seek for decisions that approach to
the minimizers of the optimization problem, the \ODAA adapts its
iterations online via a Subgradient Algorithm as described in
Section~\ref{sec:LowJ}.  We show in Section~\ref{sec:DataAssimilate}
that the resulting \ODAA is finitely convergent in the sense that the
confidence of the out-of-sample performance guarantee for the
generated data-driven decision converges to $1$ as the number of data
samples increases to a sufficiently large but finite value. Under this
scheme, a data-driven decision with certain performance guarantee is
also available any time as soon as the algorithm finishes generating
the first certificate for the initial decision, which resolves the
challenge 3). To expedite the algorithm and deal with challenge 4), we
develop in Section~\ref{sec:Heu} an \ICA (\ICAaco) to obtain
low-dimensional ambiguity sets. These new sets are based on a weighted
version of the empirical distribution and thus close to the full
empirical distribution \reviseone{of the
  data.}
We finally illustrate the performance of the proposed \ODAA in
Section~\ref{sec:Sim}, with and without the \ICAaco. A preliminary study
of this work has appeared in~\cite{DL-SM:18-cdc}.

\section{Preliminaries}\label{sec:pre}
This section introduces  basic notations and convexity
  definitions\footnote{Let $\real^m$, $\realnonnegative^m$ and
  $\real^{m \times d}$ denote respectively the $m$-dimensional
  Euclidean space, the $m$-dimensional nonnegative orthant, and the
  space of $m \times d$ matrices. We let $\vect{x} \in \real^m$ denote
  a column vector of dimension $m$, while $\trans{\vect{x}}$
  represents its transpose. We say $\vect{x} \geq 0$, if all its
  entries are nonnegative. We use the shorthand notations
  $\vectorzeros{m}$ for the column vector $\trans{(0,\cdots,0)} \in
  \real^m$, $\vectorones{m}$ for the column vector
  $\trans{(1,\cdots,1)} \in \real^m$, and $\vect{I}_m \in \real^{m
    \times m}$ for the identity matrix. We use either subscripts or
  parentheses superscripts to index vectors, i.e., $\vect{x}_k \in
  \real^m$ or $\vect{x}^{(k)} \in \real^m$, for $k \in
  \{1,2,\ldots,n\}$.  We use $(\vect{x},\vect{y}) \in \real^{m+d}$ to
  indicate the concatenated column vector from $\vect{x} \in \real^m$
  and $\vect{y} \in \real^d$. The $1$-norm of the vector $\vect{x} \in
  \real^m$ is denoted by $\Norm{\vect{x}}$.  We define the
  $m$-dimensional Euclidean ball centered at $\vect{x} \in \real^m$
  with radius $\omega$ as the set
  $B_{\omega}(\vect{x}):=\setdef{\vect{y} \in \real^m}{
    \Norm{\vect{y}-\vect{x}} \leq \omega}$.  Given a set of points $I$
  in $\real^m$, we let $\conv(I)$ indicate its convex hull. The
  gradient of a real-valued function $\map{f}{\real^m}{\real}$ is
  written as $\nabla_{\vect{x}}f(\vect{x})$. The $\supscr{i}{th}$
  component of the gradient vector is denoted by
  $\nabla_if(\vect{x})$. We use $\dom f$ to denote the domain of the
  function $f$, i.e., $\dom f:=\setdef{\vect{x} \in \real^m}{-\infty <
    f(\vect{x}) <+ \infty}$. We call the function $f$ \textit{proper}
  if $\dom f \ne \O$.  We say a function $\map{F}{\mathcal{X}\times
    \mathcal{Y}}{\real}$ is \textit{convex-concave} on
  $\mathcal{X}\times \mathcal{Y}$ if, for any point
  $(\tilde{\vect{x}},\tilde{\vect{y}}) \in \mathcal{X}\times
  \mathcal{Y}$, $\vect{x} \mapsto F(\vect{x},\tilde{\vect{y}})$ is
  convex and $\vect{y} \mapsto F(\tilde{\vect{x}},\vect{y})$ is
  concave. We use the notation $\map{\sgn}{\real}{\real}$, $x \mapsto
  \{-1,0,1 \}$ denote the sign function.  Finally, the projection
  operator
  $\map{\proj_{\mathcal{Y}}(\mathcal{X})}{\mathcal{X}}{\mathcal{Y}}$
  projects the set $\mathcal{X}$ onto $\mathcal{Y}$ under the
  Euclidean norm.}, including some from Probability Theory \reviseone{to}
describe the distributionally robust optimization framework
following~\cite{PME-DK:17}.

Let $(\Omega,\mathcal{F},\prob)$ be a probability space, with $\Omega$
the sample space, $\mathcal{F}$ a $\sigma$-algebra on $\Omega$, and
$\prob$ the associated probability distribution. Let
$\map{\xi}{\Omega}{\real^m}$ be an induced multivariate random
variable.
We denote by $\mathcal{Z} \subseteq \real^m$ the support of the random
variable $\xi$ and denote by $\mathcal{M}(\mathcal{Z})$ the space of
all probability distributions supported on $\mathcal{Z}$ with finite
first moment. In particular, $\prob \in \mathcal{M}(\mathcal{Z})$. To
measure the distance between distributions in
$\mathcal{M}(\mathcal{Z})$,
in this paper we use the dual characterization of the Wasserstein
metric~\cite{KLV-RGS:58} $\map{d_W}{\mathcal{M}(\mathcal{Z}) \times
  \mathcal{M}(\mathcal{Z})}{\realnonnegative}$, defined by
\begin{equation*}
  d_W(\mathbb{Q}_1,\mathbb{Q}_2):= \sup_{f \in
    \mathcal{L}}\int_{\mathcal{Z}} {f(\xi)\mathbb{Q}_1(d{\xi})}
  -\int_{\mathcal{Z}} {f(\xi)\mathbb{Q}_2(d{\xi})} ,
\end{equation*}
where $\mathcal{L}$ is the space of all Lipschitz functions defined on
$\mathcal{Z}$ with Lipschitz constant 1. A closed Wasserstein ball of
radius $\omega$ centered at a distribution $\prob \in
\mathcal{M}(\mathcal{Z})$ is denoted by
$\mathbb{B}_{\omega}(\prob):=\setdef{\mathbb{Q} \in
  \mathcal{M}(\mathcal{Z})}{d_W(\prob,\mathbb{Q}) \le
  \omega}$. \reviseone{We define the Dirac measure at $x_0 \in \Omega$
  as $\map{\delta_{\{x_0\}}}{\Omega}{\{0,1 \}}$. For any set $A \in
  \mathcal{F}$, we let $\delta_{\{x_0\}}(A)=1$, if $x_0 \in A$,
  otherwise $0$.}

\section{Problem Description}\label{sec:ProbStat}
Consider a decision-making problem under of the form
\begin{equation}
  \inf\limits_{\vect{x} \in \real^d}{\mathbb{E}_{\prob} [f(\vect{x},\xi)] },
  \label{eq:P}\tag{P}
\end{equation}
where $\vect{x} \in \real^d$ is the decision variable, the uncertainty random
variable $\map{\xi}{\Omega}{\real^m}$ is induced by the probability
space $(\Omega,\mathcal{F},\prob)$, and the expectation of $f$ is
taken w.r.t.~the unknown distribution $\prob \in
\mathcal{M}(\mathcal{Z})$.  We aim to develop an \ODAAfull (\ODAA)
that efficiently adapts \reviseone{iterations on} decisions $\vect{x}$
of~\eqref{eq:P} with streaming data. \reviseone{The streaming data are
  sequentially available iid realizations of the random variable $\xi$
  under $\prob$, denoted by ${\xi}_{n}$, $n=1,2,\ldots$. This defines
  a sequence of streaming data sets, ${\Xi}_{n} \subseteq
  {\Xi}_{n+1}$, for each $n$. W.l.o.g.  assume that each ${\Xi}_{n+1}$
  consists of just one more new data point, i.e., ${\Xi}_{n+1}=
  {\Xi}_{n} \cup \{{\xi}_{n+1}\}$ and ${\Xi}_{1}=\{ {\xi}_{1} \}$. In
  the following, we refer to the time slot between the updates
  ${\Xi}_n$ and ${\Xi}_{n+1}$ as the
  \textit{$\supscr{n}{th}$-time period} and to its rate of change as
  the \textit{data-streaming rate}.}

\reviseone{In practice, we cannot evaluate the objective function
  of~\eqref{eq:P} because $\prob$ is unknown.} We call a
  decision $\vect{x} \in \real^d$ \textit{a proper data-driven
    decision} of~\eqref{eq:P}, if
  its \textit{out-of-sample performance}, defined by
  ${\mathbb{E}_{\prob} [f(\vect{x},\xi)] }$, satisfies the
  \textit{performance guarantee}
\begin{equation}
  \mathbf{P}^n({\mathbb{E}_{\prob} [f(\vect{x},\xi)] } \leq
  {J}_{n}(\vect{x}))\geq 1- \beta_{n} ,
\label{eq:perfgua}
\end{equation}
where the expected cost upper bound or \textit{certificate} ${J}_{n}(\vect{x})$
is a
function that indicates the goodness of $\vect{x}$ under the
data set ${\Xi}_{n}$.  If $\vect{x}$ is adopted during the
$\supscr{n}{th}$ time period,
then
${\mathbb{E}_{\prob} [f(\vect{x},\xi)] } \leq
{J}_{n}(\vect{x})$ is an event that depends on the $n$ samples
in $\Xi_n$, and $\mathbf{P}^n$ denotes the probability with respect to
these.
The \textit{confidence} $1- \beta_{n} \in (0,1) \subset \real$ governs
the choice of $\vect{x}$ and the resulting certificate
${J}_{n}(\vect{x})$. \reviseone{In words, the
inequality~\eqref{eq:perfgua} establishes that, given finite data $n$,
the performance of the decision under the unknown distribution will
not surpass the upper-bound certificate $J_n(\vect{x})$ with high
probability. In the following section, we will determine the values
${J}_{n}$ via the solution of a parameterized maximization problem over
$\vect{x}$. Therefore,} finding an approximate certificate will be
much easier than finding the exact one.  Based on this, we call
$\vect{x}$ \textit{$\epsilon_1$-proper}, if it
satisfies~\eqref{eq:perfgua} with
\reviseone{${J}_n^{\epsilon_1}(\vect{x})$} such that
\reviseone{${J}_{n}(\vect{x}) \leq {J}_{n}^{\epsilon_1}(\vect{x}) +
  \epsilon_1
  $.}
\reviseone{Thus, the} approximates
\reviseone{${J}_{n}^{\epsilon_1}(\vect{x})$}
\reviseone{also} provide upper bounds to the optimal value
of~\eqref{eq:P} with high confidence $1-\beta_{n}$.
  \begin{figure}[tb]
     \centering
  \psfragscanon
      \includegraphics[width=0.49\textwidth]{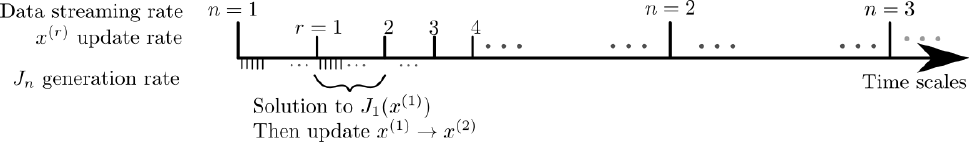}
                 \caption{ \revisetwo{ \footnotesize{Time scales of \ODAAfull}}}
     \label{fig:Cover}
   \end{figure}

To sum up, \reviseone{for any} $n$, given a confidence level
$1-\beta_{n}$, our goal is to approach to an $\epsilon_1$-proper
data-driven decision with a low certificate. \reviseone{Later in
  Section~\ref{sec:LowJ}, we will show, under assumptions on $f$,
  that \revisetwo{both Problem~\eqref{eq:P} and} these certificates ${J}_{n}$ are convex.
  To find a
  decision with a low certificate, we will} call any proper
data-driven decision \textit{$\epsilon_2$-optimal}, labeled as
$\vect{x}_n^{\epsilon_2}$, if ${J}_{n}(\vect{x}_n^{\epsilon_2})
-{J}_{n}(\vect{x}) \leq \epsilon_2$ for all $\vect{x} \in \real^d$.
Then, for any $\epsilon_2$-optimal and $\epsilon_1$-proper data-driven
decision $\vect{x}_n^{\epsilon_2}$ with certificate
${J}_{n}^{\epsilon_1}(\vect{x}_n^{\epsilon_2})$ and $\epsilon_1\ll
\epsilon_2$, we have the guarantee
\begin{equation}
  \mathbf{P}^n({\mathbb{E}_{\prob} [f(\vect{x}_n^{\epsilon_2},\xi)] } \leq
  {J}_{n}^{\epsilon_1}(\vect{x}_n^{\epsilon_2}) + \epsilon_1)\geq 1- \beta_{n}.
\label{eq:epsiperfgua}
\end{equation}
\revisetwo{Then, any decision $\vect{x}_n^{\epsilon_2}$ ensures a high-confidence, potentially-low objective value of~\eqref{eq:P}, upper bounded by ${J}_{n}^{\epsilon_1}(\vect{x}_n^{\epsilon_2}) + \epsilon_1$.}

\revisetwo{
  \subsubsection{Motivating Example in Portfolio Optimization}
  Consider an agent who does short-term trading, i.e., she is to
  select a minute-based portfolio weight $\vect{x}:=(x,1-x)$, $0\leq x
  \leq 1$, for two risky assets that give random returns, with return
  rates $\xi:=(\xi_1,\xi_2)$ following some unknown distribution
  $\prob$. The agent aims to select $\vect{x}$ such that the expected
  profit is maximized, or equivalently, she seeks to solve
\begin{equation*}
  \min\limits_{ 0\leq x \leq 1}   \mathbb{E}_{\prob}[{- \xi_1 x - \xi_2 (1-x) }].
  \tag{P0} \label{eq:P0}
\end{equation*}
Assume that $\prob$ is unknown and independent from the selection of
the portfolio and that, at every minute, the agent has access to
return rates, which are iid samples of $\prob$.
Due to the independence of $\prob$ and $\vect{x}$,
Problem~\eqref{eq:P0} is equivalent to
$\min\limits_{ 0\leq x \leq 1}   \mathbb{E}_{\prob}[{- \xi_1 x - \xi_2
  (1-x) - \trans{\xi}\xi  }].$
To adapt this problem to our unconstrained setting, we consider
an approximation
\begin{equation*}
\min\limits_{  x }   \mathbb{E}_{\prob}[{- \xi_1 x - \xi_2 (1-x) - \trans{\xi}\xi  }] -  \rho \left( \log(x) +\log(1-x)\right),
\end{equation*}
where $\rho>0$ is some penalty for the constraint terms. This problem
is in form~\eqref{eq:P}, fitting in our problem setting with
\begin{equation*}
  f(\vect{x},\xi):= {- \xi_1 x - \xi_2 (1-x) } -  \rho \left( \log(x) +\log(1-x)\right)  - \trans{\xi}\xi.
\end{equation*}
Our proposed algorithm allows the agent to make online decisions
$\vect{x}_n^{\epsilon_2}$ that minimize the objective with high
confidence. }
\subsection{Main Algorithm Goal and its High-level  Procedure}
We describe now the \reviseone{goal} of the \ODAA \,\reviseone{that
  handles a streaming sequence of
  data sets with $n\in \until{N}$.}
Let tolerance parameters $\epsilon_1$ and $\epsilon_2$ be given and
let us choose strictly \reviseone{increasing} confidence levels
$\{1 -\beta_{n}\}_{n=1}^{N}$
  such that $\sum_{n=1}^{\infty}
\beta_{n} <\infty$
 whenever $N \rightarrow \infty$. The
algorithm aims to find a sequence of $\epsilon_2$-optimal and
$\epsilon_1$-proper decisions $\{\vect{x}_n^{\epsilon_2}\}_{n=1}^{N}$
associated with the sequence of the certificates
$\{{J}_n^{\epsilon_1}(\vect{x}_n^{\epsilon_2})
\}_{n=1}^{N}$ so that~\eqref{eq:epsiperfgua} holds for all $n \in
\until{N}$. Additionally, as the data set streams to infinite
cardinality, i.e., $N \rightarrow \infty$, there exists a large enough
but finite $n_0$ such that the algorithm returns a final
$\vect{x}_{n_0}^{\epsilon_2}$ after processing the data set
${\Xi}_{n_0}$. The final decision $\vect{x}_{n_0}^{\epsilon_2}$
guarantees performance almost surely, that is,
$\mathbf{P}^{n_0}({\mathbb{E}_{\prob}
  [f(\vect{x}_{n_0}^{\epsilon_2},\xi)] } \leq
{J}_{n_0}^{\epsilon_1}(\vect{x}_{n_0}^{\epsilon_2}) + \epsilon_1)=1$,
with a certificate
${J}_{n_0}^{\epsilon_1}(\vect{x}_{n_0}^{\epsilon_2})$ close to the
optimal objective value of Problem~\eqref{eq:P}.

\reviseone{To achieve this goal, the algorithm will output a sequence
  of decisions, $\{\vect{x}^{(r)}\}_{r=1}^{\infty}$, on a time scale
  that is faster than the data-streaming rate. We refer to this as the
  \textit{decision-update rate} and we use a parenthetical superscript
  $(r)$ to denote its iterations.  New data arrival will reset
  the \ODAA's subroutines to update the sub-sequences of decisions in
  each $\supscr{n}{th}$ time period efficiently. We denote these
  sub-sequences as $\{\vect{x}^{(r)}\}_{r=r_n}^{r_{n+1}}$, where
  $\vect{x}^{(r_n)}$ is the initial decision adapted from the the time
  period $n-1$. These updates will require the computation of
  certificates (cost upper bounds) and the progressive reduction of these
 bounds.

  The computation of certificates is carried out by the \CGA (or
  \CGAaco, for short). Given a current decision value,
  $\vect{x}^{(r)}$, and data set $\Xi_n$, the \CGAaco finds an
  $\epsilon_1$ certificate $J_n^{\epsilon_1}(\vect{x}^{(r)})$ and a
  worst-case distribution associated with the data set. It operates on
  a faster time scale than the decision-update rate, the so-called
  \textit{certificate-generation rate}. Upon the receipt of new data,
  this algorithm will reset as described in Section~\ref{sec:CG}.

  The second process of the \ODAA relies on iterating decisions to
  reduce the values of the functions
  $J_{n}^{\epsilon_1}(\vect{x})$. This employs the Subgradient
  Algorithm and is described in Section~\ref{sec:LowJ}. A more
  thorough description of how new data triggers a reset in the
  algorithm is described in the following sections. A summary of the
  \ODAA can be found in Section~\ref{sec:DataAssimilate} as well as a
  descriptive table.}

\section{Certificate Design}\label{sec:CD}
\reviseone{In this section, we present a tractable formulation of
  certificates ${J}_{n}(\vect{x})$ and its approximation
  ${J}_{n}^{\epsilon_1}(\vect{x})$ for a fixed $\vect{x} = \vect{x}^{(r)}$, as
  described in~\eqref{eq:perfgua} and~\eqref{eq:epsiperfgua},
  respectively. To achieve this, we first
  follow~\cite{AC-JC:17-allerton,AC-JC:17-tac,PME-DK:17} on DRO to
  find certificates ${J}_{n}$. This defines a parameterized
  maximization problem for ${J}_{n}$, called
  Problem~\eqref{eq:convJn}. Then we reformulate~\eqref{eq:convJn} as
  Problem~\eqref{eq:JoverSimplex}, a convex optimization problem over
  a simplex, for efficient solutions of approximated certificates
  ${J}_{n}^{\epsilon_1}$ in the next section.}

To \reviseone{design ${J}_{n}$, a reasonable attempt is to}
use the data ${\Xi}_{n}$ to estimate an empirical
distribution, $\hat{\prob}^n$, and let
${\mathbb{E}_{\mathbb{\hat{P}^{\it n}}} [f(\vect{x},\xi)] }$ be the
candidate certificate for the performance
guarantee~\eqref{eq:perfgua}. More precisely, assume that the data set
${\Xi}_{n}$ are uniformly sampled from $\prob$. The discrete empirical
probability measure associated with ${\Xi}_{n}$ is the following
$  \hat{\prob}^n := \frac{1}{n}\sum_{k=1}^{n} \delta_{\{{\xi}_{k}\}},$
where $\delta_{\{ {\xi}_{k}\}}$ is the Dirac measure at
${\xi}_{k}$. The candidate certificate is
\begin{equation*}
  \supscr{J}{sae}_{n}(\vect{x}):={\mathbb{E}_{\mathbb{\hat{P}^{\it{n}}}}
    [f(\vect{x},\xi)] }
  =\frac{1}{n}\sum_{k=1}^{n}f(\vect{x},{\xi}_{k}).
\label{eq:candJ}
\end{equation*}
The above approximation $\hat{\prob}^n$ of $\prob$, also known as the
\textit{sample-average estimate}, makes $\supscr{J}{sae}_{n}$ easy to
compute. However, such value only results in an approximation of the
\reviseone{unknown} out-of-sample performance
\reviseone{$\mathbb{E}_{\mathbb{P}}[f(\vect{x},\xi)]$}.
Following~\cite{AC-JC:17-allerton,PME-DK:17}, we are to determine an
\textit{ambiguity set} $\mathcal{P}_n$ containing all the possible
probability distributions supported on $\mathcal{Z} \subseteq \real^m$
that can generate ${\Xi}_{n}$ with high confidence. Then,
\reviseone{one can} consider the worst-case expectation of
\reviseone{$f(\vect{x},\xi)$}
\reviseone{with respect to} all
distributions contained in $\mathcal{P}_n$. The solution to such
problem offers an upper bound for the out-of-sample performance with
high probability in the form of~\eqref{eq:perfgua}, and we
refer to this upper bound as the certificate of the decision
$\vect{x}$.

In order to quantify the \reviseone{ambiguity set and} certificate for an $\epsilon_1$-proper
data-driven decision, we denote by
$\subscr{\mathcal{M}}{lt}(\mathcal{Z}) \subset
\mathcal{M}(\mathcal{Z})$ the set of light-tailed probability measures
in $\mathcal{M}(\mathcal{Z})$, and introduce the following assumption
for $\prob$
\begin{assumption}[Light tailed unknown distributions]
  It holds that $\prob \in
  \subscr{\mathcal{M}}{lt}(\mathcal{Z})$, i.e., there exists an
  exponent $a>1$ such that $b:= \mathbb{E}_{\prob}
  [\exp(\Norm{\xi}^a)] < \infty$. \label{assump:2}
\end{assumption}
\begin{remark}[Class of distributions satisfying Assumption~\ref{assump:2}]
 {\reviseone{\rm Any distribution with an exponentially decaying tail satisfies this
  assumption, such as Gaussian, subGaussian, exponential, and geometric
  distributions. Any distribution with a compact support $\mathcal{Z}$
  will trivially satisfy the assumption. In engineering problems,
  the values of random variables are usually truncated to a compact set
  and hence the Assumption~\ref{assump:2} is automatically satisfied.}}
 \end{remark}
Assumption~\ref{assump:2} validates the following modern measure
concentration result, which provides an intuition for considering
the Wasserstein ball $\mathbb{B}_{\epsilon}(\hat{\prob}^n)$ of center
$\hat{\prob}^n$ and radius $\epsilon$ as the ambiguity set
$\mathcal{P}_n$.

\begin{theorem}[Measure concentration~{\cite[Theorem~2]{NF-AG:15}}]
  If $\prob \in \subscr{\mathcal{M}}{lt}(\mathcal{Z})$, then
\begin{equation}
  \mathbf{P}^n\{ d_W(\prob,\hat{\prob}^n) \geq \epsilon \}
  \leq \left\{ {\begin{array}{*{20}{l}}
        c_1 e^{-c_2n \epsilon^{\max \{2,m\}}}, \quad & \textrm{if} \;
        \epsilon \leq 1,    \\
        c_1 e^{-c_2n \epsilon^a} , \quad &  \textrm{if} \; \epsilon > 1,
\end{array}} \right.
\label{eq:radB}
\end{equation}
for all $n \geq 1$, $m\ne 2$, and $\epsilon >0$, where $c_1$, $c_2$
are positive constants that only depend on $m$, $a$ and $b$. \hfill $\square$
\label{thm:MeasConc}
\end{theorem}

Equipped with this result,
we are able to provide the certificate that ensures the performance
guarantee in~\eqref{eq:perfgua}, for any decision $\vect{x} \in \real^d$.

\begin{lemma}[Certificate in Performance
  Guarantee~\eqref{eq:perfgua}]
  Given ${\Xi}_{n}:=\{{\xi}_{k} \}_{k=1}^{n}$,
\reviseone{$\beta_n \in (0,1)$ and $\vect{x} \in \real^d$},
	let
\begin{equation}
  \epsilon(\beta_{n}) := \left\{ {\begin{array}{*{20}{l}}
        \left( \frac{\log(c_1 \beta_{n}^{-1})}{c_2 n}
        \right)^{1/\max\{2,m\}},
        \; & \textrm{if} \; n \geq \frac{\log(c_1 \beta_{n}^{-1})}{c_2}, \\
        \left( \frac{\log(c_1 \beta_{n}^{-1})}{c_2 n} \right)^{1/a},
        \; & \textrm{if} \; n < \frac{\log(c_1 \beta_{n}^{-1})}{c_2},
\end{array}} \right.
\label{eq:epsirad}
\end{equation}
and
$\mathcal{P}_n:=\mathbb{B}_{\epsilon(\beta_{n})}(\hat{\prob}^n)$. Then
the following certificate satisfies the performance guarantee
in~\eqref{eq:perfgua} for all $\vect{x} \in \real^d$
\begin{equation}
  {J}_{n}(\vect{x})
  := \sup\limits_{\mathbb{Q}
    \in \mathcal{P}_n} {\mathbb{E}_{\mathbb{Q}} [f(\vect{x},\xi)] }.
\label{eq:cert}
\end{equation}
\label{lemma:certgener}
\end{lemma}
To get ${J}_{n}$ in~\eqref{eq:cert}, one needs to solve an
infinite-dimensional optimization problem. Luckily,
Problem~\eqref{eq:cert} can be reformulated into a
finite-dimensional convex  problem as follows.

\begin{theorem}[Convex reduction of~\eqref{eq:cert}{~\cite[Application of Theorem
  4.4]{PME-DK:17}}]
Under
Assumption~\ref{assump:2}, \reviseone{on $\prob$ being light-tailed},
for all $\beta_{n} \in (0,1)$ the value of the certificate
in~\eqref{eq:cert} for \reviseone{a given} decision $\vect{x}$ under
\reviseone{a} data set ${\Xi}_{n}$ is equal to the optimal value of the
following optimization problem
\begin{equation}
\begin{aligned}
  {J}_{n}(\vect{x}):=
  \sup\limits_{ \vect{y}
  } \; &\;
  \frac{1}{n}\sum_{k=1}^{n}
  f(\vect{x},{\xi}_{k}-\vect{y}_k) ,\\
  \st \quad &\frac{1}{n} \sum_{k=1}^{n} \Norm{\vect{y}_k} \le
  \epsilon(\beta_{n}),
\end{aligned}
\label{eq:convJn}\tag{\reviseone{P1$_{n}$}}
\end{equation}
where \reviseone{each component of the concatenated variable
  $\vect{y}:=(\vect{y}_1,\ldots,\vect{y}_n)$ is in $\real^m$, and the
  parameter} $\epsilon(\beta_{n})$ is the radius of
$\mathbb{B}_{\epsilon(\beta_{n})}$ calculated
from~\eqref{eq:epsirad}. Moreover, given any feasible point
$\vect{y}^{(l)}:=(\vect{y}^{(l)}_1,\ldots,\vect{y}^{(l)}_{n})$
of~\eqref{eq:convJn}, indexed by $l$, define a finite atomic
probability measure at $\vect{x}$ in the Wasserstein ball
$\mathbb{B}_{\epsilon(\beta_{n})}$ of the form
\begin{equation}
\begin{aligned}
\mathbb{Q}_{n}^{(l)}(\vect{x}):=\frac{1}{n}\sum_{k=1}^{n} \delta_{\{{\xi}_{k}-\vect{y}^{(l)}_k \}}.
\label{eq:11}
\end{aligned}
\end{equation}
Now, denote by $\mathbb{Q}_{n}^{\star}(\vect{x})$ the
distribution in~\eqref{eq:11} constructed by an optimizer
$\vect{y}^{\star}:=(\vect{y}^{\star}_1,\ldots,\vect{y}^{\star}_n)$
of~\eqref{eq:convJn} and evaluated over $\vect{x}$. Then,
$\mathbb{Q}_{n}^{\star}$ is a worst-case distribution that can
generate the data set ${\Xi}_{n}$ with (high) probability no less than
$1-\beta_{n}$. \hfill $\square$
\label{thm:cvxredu}
\end{theorem}
  \begin{remark} { \reviseone{\rm Theorem~\ref{thm:cvxredu} provides a way of
      computing certificates of~\eqref{eq:perfgua} as the solution to
      a parameterized optimization problem for a decision
      $\vect{x}$.  In addition, it constructs a worst-case
      distribution that achieves the worst-case bound.}}
 \end{remark}
 \reviseone{To enable efficient online solutions of approximated certificates ${J}_{n}^{\epsilon_1}$, let us define parameterized functions $\map{h_k}{\real^m}{\real}$}
\begin{center}
  $h_k(\vect{y}):=f(\vect{x},{\xi}_{k}-\vect{y}), \quad \reviseone{k \in \until{n}}$,
\end{center}
\reviseone{and} consider the following convex optimization problem
over a simplex
\begin{equation}
  \begin{aligned}
    {J}_{n}(\vect{x}):= \max\limits_{ \vect{u}, \vect{v}
    } \; &\; \frac{1}{n}\sum_{k=1}^{n}
    h_k( \vect{u}_k-\vect{v}_k) ,\\
    \st \quad & \; (\vect{u}, \vect{v}) \in n\epsilon(\beta_{n}) \Delta_{2mn},
\end{aligned}
\label{eq:JoverSimplex}\tag{\reviseone{P2$_{n}$}}
\end{equation}
where the concatenated variable $(\vect{u}, \vect{v})$ is composed of
 $\vect{u}:=(\vect{u}_1,\ldots,\vect{u}_n)$ and
$\vect{v}:=(\vect{v}_1,\ldots,\vect{v}_n)$ \reviseone{with $\vect{u}_k$, $\vect{v}_k\in \real^m$ for all $k \in \until{n}$}; and
the scalar $n\epsilon(\beta_{n})$ regulates the size of the feasible
set via scaling of the unit simplex $\Delta_{2mn}:=\setdef{(\vect{u},
  \vect{v}) \in \real^{2mn}}{\trans{\vectorones{2mn}}(\vect{u},
  \vect{v})=1 ,\; \vect{u} \geq 0,\; \vect{v} \geq 0}$.  We denote by
${\Lambda}_{2mn}$ the set of all the extreme points for the simplex
$n\epsilon(\beta_{n})\Delta_{2mn}$.

 The following lemma shows that Problem~\eqref{eq:convJn} and
 Problem~\eqref{eq:JoverSimplex} are equivalent. \reviseone{Thus, we
   can approximately solve~\eqref{eq:JoverSimplex} to find}
 ${J}_n^{\epsilon_1}(\vect{x})$ and
 $\mathbb{Q}_{n}^{\epsilon_1}(\vect{x})$.

 \begin{lemma}[Equivalence of the problem formulation]
   Solving~\eqref{eq:convJn} is equivalent to solving~\eqref{eq:JoverSimplex}
   in the sense that
 \begin{itemize}
 \item[1] For any feasible solution $(\tilde{\vect{u}},\tilde{\vect{v}})$
   of~\eqref{eq:JoverSimplex}, let $\tilde{\vect{y}}:=\tilde{\vect{u}}-\tilde{\vect{v}}$. Then $\tilde{\vect{y}}$
   is feasible for~\eqref{eq:convJn}.
 \item[2] For any feasible solution
   $\tilde{\vect{y}}$
   of~\eqref{eq:convJn}, there exists a feasible point
   $(\tilde{\vect{u}},\tilde{\vect{v}})$
   of~\eqref{eq:JoverSimplex}.
 \item[3] Assume that the point $(\tilde{\vect{u}}^{\star},\tilde{\vect{v}}^{\star})$ is an optimizer of~\eqref{eq:JoverSimplex}. Then by letting $\tilde{\vect{y}}^{\star}:=\tilde{\vect{u}}^{\star}-\tilde{\vect{v}}^{\star}$, the point $\tilde{\vect{y}}^{\star}$
   is also an optimizer of~\eqref{eq:convJn}, with the same optimal
   value.
 \end{itemize}\label{lemma:equiP}
 \end{lemma}

 \reviseone{For this section, we provide proofs of Lemma~\ref{lemma:certgener} and Lemma~\ref{lemma:equiP}
 in the Appendix, and from now on, all the proofs of the lemmas and
 theorems can be found in the Appendix.}

 \section{\CGA} \label{sec:CG} Given a tolerance
 $\epsilon_1$, sequentially available data sets
 $\{{\Xi}_{n}\}_{n=1}^{N}$ and decisions
 $\{\vect{x}^{(r)}\}_{r=1}^{\infty}$, we present in this section the
 \CGA \,(\CGAaco) \reviseone{to obtain approximated
   certificates}
 $\{{J}_n^{\epsilon_1}(\vect{x}^{(r)})\}_{n,r}$
  and
 \reviseone{associated}
 $\epsilon_1$-worst-case distributions
 $\{\mathbb{Q}_{n}^{\epsilon_1}(\vect{x}^{(r)})\}_{n,r}$.
 To achieve this, we \reviseone{first} design \reviseone{for each fixed $\vect{x} = \vect{x}^{(r)}$ the \CGAaco} to solve~\eqref{eq:JoverSimplex} \reviseone{to}
 ${J}_n^{\epsilon_1}(\vect{x})$ efficiently. \reviseone{This is developed}
 via Frank-Wolfe Algorithm
 variants, e.g., the Simplicial Algorithm~\cite{CH:74}
 and the AFWA as \reviseone{described in the Appendix. Then} we analyze
 the convergence of \reviseone{the \CGAaco} under
 $\{{\Xi}_{n}\}_{n=1}^{N}$.

 \subsection{The \CGAaco}

\reviseone{For each fixed $\vect{x}=\vect{x}^{(r)}\in \real^d$} the
algorithm \reviseone{is run at a fast time scale (the
  certificate-update rate), over iterations $l = 0, 1, 2, \dots$. The
  algorithm is then employed inside the \ODAA, so its execution rate is the
  fastest within this algorithm.}  At each iteration
$l$, the \CGAaco generates
$(\vect{u}^{(l)},\vect{v}^{(l)})$, the candidate optimizer
of~\eqref{eq:JoverSimplex}. Let the objective value
of~\eqref{eq:JoverSimplex} at $(\vect{u}^{(l)},\vect{v}^{(l)})$ be
${J}_{n}^{(l)}(\vect{x})$, and, equivalently, write the candidate
optimizer in form of $\vect{y}^{(l)}:=\vect{u}^{(l)}-\vect{v}^{(l)}$
(exploiting the equivalence in
Lemma~\ref{lemma:equiP}).
Each candidate $\vect{y}^{(l)}$ is associated with  a set of
search points denoted by
$I_n^{(l)}:=\{\tilde{\vect{y}}^{[i]}:=\tilde{\vect{u}}^{[i]}-\tilde{\vect{v}}^{[i]},
\; i \in \until{T} \}$, \reviseone{where we use bracket superscript
  $[i]$ to index its elements.}  As we will see later,
the set $I_{n}^{(l)}$ plays \reviseone{a key} role in  generating the certificate
when assimilating data, and is called the \emph{candidate vertex set}.

\reviseone{Given a data set $\Xi_n$, and until new data arrives,} the
\CGAaco
solves the following problems \reviseone{alternatively}
\begin{equation}
\begin{aligned}
 \max\limits_{ \vect{u}, \vect{v}} \;
  & \frac{1}{n}\sum_{k=1}^{n}
  \left\langle  \nabla h_k(\vect{y}^{(l-1)}_k), \cdots \right.\\ & \left.\hspace*{2.0cm} \vect{u}_k-\vect{v}_k -\vect{y}^{(l-1)}_k\right\rangle ,\\
  \st \quad & (\vect{u}, \vect{v})  \in n\epsilon(\beta_{n}) \Delta_{2mn},
\label{eq:LP}
\end{aligned}\tag{LP$_n^{(l)}$}
\end{equation}
\begin{equation}
\begin{aligned}
 \max\limits_{\gamma \in \real^{T} } \; & \frac{1}{n}\sum_{k=1}^{n}
h_k(\sum_{i=0}^{T}\gamma_i  \tilde{\vect{y}}_k^{[i]}) ,\\
\st \quad & \gamma \in \Delta_{T}.
\label{eq:CP}
\end{aligned}\tag{CP$_n^{(l)}$}
\end{equation}
\reviseone{(Note how the solution to one problem parameterizes the
  other.) In this way,} $\vect{y}^{(l-1)}$ (\reviseone{the solution to
  the $\text{CP}_n^{(l-1)}$}) parameterizes the linear
problem~\eqref{eq:LP}. The solution to~\eqref{eq:LP} is then used to
refine the \reviseone{set point}
$I_{n}^{(l)}=\{\tilde{\vect{y}}^{[i]}\}_i$, which spans the
\reviseone{constraint} set $\Delta_T \equiv \conv(I_{n}^{(l)})$ in
problem~\eqref{eq:CP}.
A solution to~\eqref{eq:CP} then determines the
new
$\vect{y}^{(l)}$ of the \reviseone{next} LP problem. \reviseone{This
  process corresponds to lines 3: to 9: in the following \CGAaco
  table.}

{\begin{algorithm}[tbp]
\floatname{algorithm}{\CGAaco}{}
    \caption{CG$(\vect{x}, \{\Xi_n\}_{n=1}^N, \vect{y}^{(0)}, I_{n}^{(0)})$}
\label{alg:cert}
\begin{algorithmic}[1]
\Require{Goes to Step~\ref{step:CG_data} \reviseone{upon data
    arrival, i.e.} $\Xi_{n} \leftarrow \Xi_{n+1}$.}
\State $l \leftarrow 0$; \Comment{Procedure for $\Xi_n$} \label{step:CG_data}
\State Update $\vect{y}^{(l)}$, $I_{n}^{(l)}$, $T$ and $\gamma^{\epsilon_1}$; \Comment{Adapted from $\Xi_{n-1}$}
\Repeat
\State $l \leftarrow l+1$;
\State $(\Omega^{(l)},\eta^{(l)}) \leftarrow$
      ${\textrm{LP}}(\vect{x}, \Xi_n,
      \vect{y}^{(l-1)})$;
\State $I_{n}^{(l)} \leftarrow I_{n}^{(l-1)} \cup \Omega^{(l)}$, $T \leftarrow | {I_{n}^{(l)}}|$;
\State
 $(\gamma^{\epsilon_1}, {J}_{n}^{(l)}(\vect{x})) \leftarrow$ AFWA~\eqref{eq:CP};
\State $\vect{y}^{(l)} \leftarrow \sum_{i=0}^{T}\gamma_i^{\epsilon_1}\tilde{\vect{y}}^{[i]}$, $\tilde{\vect{y}}^{[i]} \in I_{n}^{(l)}$ for each $i$;
\Until{$\eta^{(l)} \le
  \epsilon_1$} \\
\Return
${J}_n^{\epsilon_1}(\vect{x})={J}_{n}^{(l)}(\vect{x})$,
${\vect{y}}^{\epsilon_1}=\vect{y}^{(l)}$,
$\mathbb{Q}_n^{\epsilon_1}(\vect{x})=\frac{1}{n}\sum_{k=1}^{n}
\delta_{\{{\xi}_{k}-{\vect{y}}^{\epsilon_1}_k \}}$.
\end{algorithmic}
\end{algorithm}}
{\begin{algorithm}[tbp]
\floatname{algorithm}{\reviseone{Point Search Algorithm}}{}
    \caption{${\textrm{LP}}(\vect{x}, \Xi_n,
      \vect{y}^{(l-1)})$}
\label{alg:LP}
\begin{algorithmic}[1]
\State Set $\Omega^{(l)}:=\O$;
\State Let $H:=\setdef{(j,k)}{j \in \until{m},\;k\in \until{n}}$;
\State Let $S:=\argmax_{(j,k) \in H} \{ \pm \nabla_j h_k(\vect{y}^{(l-1)}_k) \}$;
\While  {$S \ne \O$, }
\State Pick $(\hslash,\ell) \in S$ and let $\tilde{\vect{y}}=\vectorzeros{mn}$;
 \State Update scalar $\tilde{\vect{y}}_{\hslash\ell}\leftarrow n\epsilon(\beta_{n}){\sgn}(\nabla_{\hslash} \:
  h_{\ell}(\vect{y}^{(l-1)}_{\ell}))$;
\State Update $\Omega^{(l)} \leftarrow \Omega^{(l)} \cup \{\tilde{\vect{y}} \}$;
\State Update $S\leftarrow S \setminus \{(\hslash,\ell)\}$;
\EndWhile
\State Pick any $\tilde{\vect{y}}\in \Omega^{(l)}$ and,
\State set $\eta^{(l)}=\frac{1}{n}\sum_{k=1}^{n}
\left\langle  \nabla h_k(\vect{y}^{(l-1)}_k),\tilde{\vect{y}}_k -\vect{y}^{(l-1)}_k\right\rangle$; \\
\Return the set $\Omega^{(l)}$ and the optimality gap $\eta^{(l)}$.
\end{algorithmic}
\end{algorithm}}
More precisely, at each iteration $l =1,2,\ldots$,
the \CGAaco
first solves~\eqref{eq:LP} using the Point
Search Algorithm,
 which returns the optimal
objective value $\eta^{(l)}$ and the set of maximizers
$\Omega^{(l)}$
such that $\eta^{(l)} \geq
{J}_{n}(\vect{x})-{J}_{n}^{(l)}(\vect{x})$ and
$\Omega^{(l)} \subset \Lambda_{2mn}$. \reviseone{The
value $\eta^{(l)}$ is then used to determine the $\epsilon_1$-suboptimality
condition to the optimal objective of
Problem~\eqref{eq:JoverSimplex} (see below).}
 \reviseone{Meanwhile, the set
  $\Omega^{(l)}$ is used to update candidate vertex set to
  $I_{n}^{(l)}:=I_{n}^{(l-1)} \cup \Omega^{(l)}$, which is used in problem~\eqref{eq:CP}.} In particular,
\reviseone{the Point Search Algorithm}
computes all optimizers by iteratively choosing
a sparse vector with only a positive entry. That is, an extreme point
of the feasible set of~\eqref{eq:LP}, such that the nonzero component
of $(\tilde{\vect{u}}^{(l)}, \tilde{\vect{v}}^{(l)})$ has the largest
absolute gradient component in the linear cost function
of~\eqref{eq:LP}.
Using the obtained $I_{n}^{(l)}$, the algorithm solves the
Problem~\eqref{eq:CP} over the simplex $\Delta_{T}:=\setdef{\gamma \in
  \real^{T}}{\trans{\vectorones{T}}\gamma=1,\; \gamma \geq 0 }$, where \reviseone{$T$ is the cardinality of $I_n^{(l)}$ and}
each component $\gamma_i$ of $\gamma \in \Delta_{T}$ represents the
convex combination coefficient of a candidate vertex
$\tilde{\vect{y}}^{[i]}$\reviseone{$\in I_n^{(l)}$}.
After solving~\eqref{eq:CP} to $\epsilon_1$-optimality via the AFWA
(\reviseone{see Appendix}), an $\epsilon_1$-optimal weighting
$\gamma^{\epsilon_1} \in \Delta_{T}$
with the objective value
${J}_{n}^{(l)}(\vect{x})$
is obtained. A new candidate optimizer $\vect{y}^{(l)}$ is then
calculated by
$\vect{y}^{(l)}=\sum_{i=0}^{T}\gamma_i^{\epsilon_1}\tilde{\vect{y}}^{[i]}$. The
algorithm repeats the process and increments $l$ if the optimality gap
$\eta^{(l)}$ is greater than $\epsilon_1$, otherwise it returns the
certificate
${J}_n^{\epsilon_1}(\vect{x}):={J}_{n}^{(l)}(\vect{x})$,
an $\epsilon_1$-optimal solution
${\vect{y}}^{\epsilon_1}:=\vect{y}^{(l)}$ and an
$\epsilon_1$-worst-case distribution
$\mathbb{Q}_n^{\epsilon_1}(\vect{x}):=\frac{1}{n}\sum_{k=1}^{n}
\delta_{\{{\xi}_{k}-{\vect{y}}^{\epsilon_1}_k \}}$.

\reviseone{When new data arrives, the algorithm will reset by adapting
  the Problem~(P2$_{n+1}$) from Problem~\eqref{eq:JoverSimplex} (line
  2: in the table of the \CGAaco (update from $\Xi_n$ to
  $\Xi_{n+1}$). Note that} adapting the
\CGAaco
to online data sets $\{{\Xi}_{n}\}_{n=1}^{N}$ is inherently difficult
due to the changes in the Problems~\eqref{eq:JoverSimplex}. As the
size of ${\Xi}_{n}$ grows by $1$, the dimension of the
Problem~\eqref{eq:JoverSimplex} increases by $2m$.  To obtain
${J}_n^{\epsilon_1}(\vect{x})$ and
$\mathbb{Q}_{n}^{\epsilon_1}(\vect{x})$ sufficiently fast, we
exploit the relationship among Problems~\eqref{eq:JoverSimplex}, for
different $n$, by adapting the candidate vertex sets
$I_{n}^{(l)}$. Specifically, we initialize the set $I_{n+1}^{(0)}$ for
the new Problem~(P2$_{n+1}$)
by $I_{n}^{(l)}$, constructed from the
previous~\eqref{eq:JoverSimplex}. Suppose that the \CGAaco receives a
new data set ${\Xi}_{n+1} \supset {\Xi}_{n}$ at some intermediate
iteration $l$ with candidate vertex set $I_{n}^{(l)}$. At this stage,
the subset $\conv(I_{n}^{(l)})$ has been explored by the previous
optimization problem, and the gradient information of the objective
function based on the data set ${\Xi}_{n}$ has been partially
integrated. Then, by projecting the set $I_{n}^{(l)}$ onto the set of
extreme points of the new Problem~\eqref{eq:JoverSimplex}, i.e.,
$I_{n+1}^{(0)}:=\proj_{\Lambda_{2m(n+1)}}(\setdef{(\tilde{\vect{y}}^{[i]}
  , \vectorzeros{m})}{\tilde{\vect{y}}^{[i]} \in I_{n}^{(l)}})$, the
subset $\conv(I_{n+1}^{(0)})$ of the feasible set
of~\eqref{eq:JoverSimplex} is already explored.
Such integration contributes to the reduction of the
number of iterations in the \CGAaco for
Problems~\eqref{eq:JoverSimplex}.
This insight gives us a sense of the worst-case efficiency to update a
certificate under the streaming data.
\subsection{Convergence Analysis of the~\CGAaco}
\reviseone{We make} the following assumptions on the local strong
concavity of the function~$f$ and the computation of \reviseone{its
  gradient}
\begin{assumption}[Local strong concavity]
 \reviseone{For any  ${\vect{x}} \in \real^d$
and  ${\xi} \in \real^m$, }
the function {$\map{h}{\real^m}{\real}$}, $\vect{y} \mapsto
f({\vect{x}},{\xi}-\vect{y})$ is differentiable, \revisetwo{concave
  with} a curvature constant $C_h$, and \revisetwo{with} a positive
geometric strong concavity constant $\mu_h$ on
$\Delta_{2mN}$\footnote{ \label{footnote:concavity} \revisetwo{For a
    concave function $\map{h}{\real^m}{\real}$ on $\Delta$, we
    define
 $\displaystyle{C_h:= \sup -\frac{2}{\gamma^2} \left(
     h(\vect{y}^{\star}) - h(\vect{y}) - \left\langle \nabla
       h(\vect{y}), \vect{y}^{\star}-\vect{y} \right\rangle \right),} $ s.t.
      $    \vect{y}^{\star}=\vect{y}+\gamma(\vect{s} - \vect{r}),$ $
      \gamma \in [0,1]$, $ \vect{y}, \vect{s}, \vect{r} \in {\Delta}.$
and   $\displaystyle{  \mu_h:= \inf\limits_{\vect{y} \in \Delta}
  \inf\limits_{\vect{y}^{\star} \in \Delta}   -
  \frac{2}{\Gamma(\vect{y},\vect{y}^{\star})^2} \times } $ $ \left( h(\vect{y}^{\star})-h(\vect{y})- \left\langle \nabla h(\vect{y}), \vect{y}^{\star}-\vect{y} \right\rangle \right), $
    s.t.~$\left\langle \nabla h(\vect{y}), \vect{y}^{\star}-\vect{y}
    \right\rangle > 0,$
  where $\Gamma(\vect{y},\vect{y}^{\star})$ is a step-size measure in
  AFWA. See, e.g.,~\cite{SLJ-MJ:15} for details.
  We say $h$ is locally strongly concave, if $\mu_h>0$.}}.
\label{assump:localconcave}
\end{assumption}

\begin{assumption}[
  \revisetwo{Accessible} gradients]
  For any decision \reviseone{${\vect{x}} \in
  \real^d$}, we denote by $\nabla
  h(\vect{y})$ the gradient of the function
  {$\map{h}{\real^m}{\real}$},
  $\vect{y} \mapsto f({\vect{x}},\vect{y})$ and assume it \revisetwo{is accessible.}
  \label{assump:chgrad}
\end{assumption}
  Under Assumptions~\ref{assump:localconcave} and~\ref{assump:chgrad},
  we show the convergence properties of the \CGAaco.
  \begin{figure}[!tbp]
     \centering
  \psfragscanon
      \includegraphics[width=0.4\textwidth]{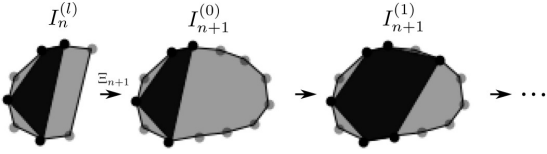}
                 \caption{ \revisetwo{ \footnotesize{\CGAaco Procedure on a projected plane. At each particular time period ($n$ or $n+1$) and iteration $l$, the  dots, shaded region and solid region represent the projection of vertices of $\Delta$, $\Delta$ and $\conv(I)$, respectively. The solid region $\conv(I)$ implicitly expands
                 for solutions to various~\eqref{eq:JoverSimplex}. }}}
     \label{fig:Cover}
   \end{figure}
\begin{theorem}[Convergence of the \CGAaco] \label{thm:ConvergeJhat}
  Let a tolerance $\epsilon_1$ and a decision $\vect{x}$
   be given. Let us choose $\vect{y}^{(0)}=\vectorzeros{m}$ and
  $I_{1}^{(0)}=\O$ as the initial candidate optimizer and candidate
  vertex set for the \CGAaco, respectively. Consider the
  online data sets $\{{\Xi}_{n}\}_{n=1}^{N}$ and the set of
  parameterized functions $\{h_n\}_{n=1}^{N}$. \reviseone{Under}
  Assumption~\ref{assump:localconcave} and
  Assumption~\ref{assump:chgrad}, \reviseone{we have that} for
  \reviseone{all} data set ${\Xi}_{n}$, there exists a parameter
  $\kappa \in (0,1) \subset \real$ such that the worst-case
  computational bound $\phi(n)$ of the \CGAaco, depending
  on $n$, is
 \begin{equation*}
   \phi(n) \leq (2mn) {\log}_{\kappa}(\frac{\epsilon_1}{{J}_n(\vect{x}) -{J}_n^{(0)}(\vect{x})}).
 \end{equation*}
 Moreover, consider that data sets $\{{\Xi}_{n}\}_{n=1}^{N}$ are
 streaming \reviseone{and consider function $\supscr{J}{sae}_{N}(\vect{x})$ defined as in Section~\ref{sec:CD}}. Then there exists a parameter $\bar{\kappa} \in (0,1)
 \subset \real$ and a computational bound
\begin{equation*}
  \bar{\phi}(n):= (2mn)
  {\log}_{\bar{\kappa}}(\frac{\epsilon_1}{{J}_{N}(\vect{x}) -
\supscr{J}{sae}_{N}(\vect{x})})
\end{equation*}
such
that, \reviseone{if the data-streaming rate is slower or equal than $(\bar{\phi}(1))^{-1}$,}
  then the \CGAaco is guaranteed to obtain the
certificates $\{{J}_n^{\epsilon_1}(\vect{x})\}_{n=1}^{N}$ and
$\{\mathbb{Q}_{n}^{\epsilon_1}(\vect{x})\}_{n=1}^{N}$.
\end{theorem}
\reviseone{Theorem~\ref{thm:ConvergeJhat} relates the worst-case
  computational bound of the \CGAaco, executed on the
  certificate-update rate (the fastest of the time scales considered),
  to the data-streaming rate. Note
  that, as $\epsilon_1$ decreases, the bound $\bar{\phi}(1)$ increases
  and therefore the smaller the data-streaming rate has to be so that
  the certificates can be generated by the algorithm.  In practice,
  the \CGAaco tends to find the smallest implicit feasible set that
  contains an optimal solution
  of~\eqref{eq:JoverSimplex}.
  This means that the computation of the \CGAaco generally performs
  better than its worst-case bound as in
  Theorem~\ref{thm:ConvergeJhat} and so it can handle data-streaming
  rates faster than $(\bar{\phi}(1))^{-1}$.  In the sequel, we assume
  that the \CGAaco converges with a rate that is faster than the
  worst-case bound in Theorem~\ref{thm:ConvergeJhat}.  }
\begin{remark}[Effects of Assumption~\ref{assump:localconcave} and~\ref{assump:chgrad}]
  \revisetwo{{ \rm The essential ingredients for convergence of
      the~\CGA are 1) the concavity of $h$, which ensures
      that~\eqref{eq:JoverSimplex} is a convex problem, and 2)
      accessible gradients of $h$, which allows for computations to a
      solution of~\eqref{eq:JoverSimplex}. To obtain a fast, linear
      convergence rate as in Theorem~\ref{thm:ConvergeJhat}, we assume
      that $h$ is strongly concave on the simplex $\Delta_{2mN}$
      located at each data point $\xi \in {\Xi}_{n}$. Intuitively, as
      ${\Xi}_{n}$ comes from $\prob$,
      Assumption~\ref{assump:localconcave} eventually requires $h$ to
      be strongly concave on a subset of the support $\mathcal{Z}$ of
      $\prob$ where the high-probability outcomes are concentrated
      onto. Otherwise, if $h$ is concave but not locally strongly
      concave, or if the gradients of $h$ are inaccessible (e.g., when
      only non-biased gradient estimate of $h$ are available),
  the convergence of AFWA, as described in  Theorem~\ref{thm:AFWconv}, reduces to a sublinear rate. This, in turn, reduces the computational bound $\bar{\phi}(n)$ in Theorem~\ref{thm:ConvergeJhat} to a bound of order $\mathcal{O}(1/\epsilon_1)$.} }
\end{remark}
\begin{remark}[Example in Portfolio Optimization]
\revisetwo{{\rm  The portfolio problem in Section~\ref{sec:ProbStat} results in a strongly concave $h$ which implies the local strong concavity as required by Assumption~\ref{assump:localconcave}. Let $\vect{y}:=(y_1,y_2)$ and, for any given data point $\xi_k:=(\xi_{k,1},\xi_{k,2})$, $\nabla h$ is accessible and computed by
\begin{equation*}
  \nabla h(\vect{y}):= \begin{pmatrix}
    x + 2( \xi_{k,1}+\xi_{k,2} -y_1 -y_2 ) \\
    1-  x + 2( \xi_{k,1}+\xi_{k,2} -y_1 -y_2 )
\end{pmatrix}.
\end{equation*}
}}
\end{remark}
\section{\reviseone{Sub-optimal decisions with guarantees}}
\label{sec:LowJ}
In this section, we aim to construct a sub-sequence of
$\epsilon_2$-optimal data-driven decisions
$\{\vect{x}_n^{\epsilon_2}\}_{n=1}^{N}$, associated with the
$\epsilon_2$-lowest certificates
$\{{J}_n^{\epsilon_1}(\vect{x}_n^{\epsilon_2})\}_{n=1}^{N}$
over time.
We achieve this by means of the Subgradient Algorithm to derive an
$\epsilon_1$-proper decision sequence
$\{\vect{x}^{(r)}\}_{r=r_n}^{r_{n+1}}$;
and the concatenation of $\{\vect{x}^{(r)}\}_{r=r_n}^{r_{n+1}}$
\reviseone{for different $n$}
to obtain $\{\vect{x}_n^{\epsilon_2}\}_{n=1}^{N}$.

To construct an $\epsilon_1$-proper decision sub-sequence
$\{\vect{x}^{(r)}\}_{r=r_n}^{r_{n+1}}$
let us consider the following problem
\begin{equation*}
  {J}_n^{\star}:=
  \inf\limits_{\vect{x} \in \real^d} {J}_{n}(\vect{x}),
\label{eq:certstar}
\end{equation*}
where the function ${J}_{n}(\vect{x})$ is defined as in
either~\eqref{eq:cert} or~\eqref{eq:convJn}, and we assume the
approximation of ${J}_{n}(\vect{x})$, ${J}_n^{\epsilon_1}(\vect{x})$,
can be evaluated as in Section~\ref{sec:CG}.

To solve this Problem to
${J}_n^{\epsilon_1}(\vect{x}_n^{\epsilon_2})$, we have the following
assumption on the convexity of $f$
\begin{assumption}[Convexity in $\vect{x}$]
  The function $\map{f_{\xi}}{\real^d}{\real}$
  ${\vect{x}}\mapsto{f(\vect{x},\xi)}$
  is
  convex for all $\xi \in \real^m$. \label{assump:cvx}
\end{assumption}

Assumption~\ref{assump:cvx} results in convexity of
${J}_{n}(\vect{x})$ as follows.
\begin{lemma}[Convexity of ${J}_{n}(\vect{x})$]
  If Assumption~\ref{assump:cvx} \reviseone{(convexity in $\vect{x}$)}
    holds, then for each $n \in
  \until{N}$ the certificate ${J}_{n}(\vect{x})$ defined
  by~\eqref{eq:cert} is convex in $\vect{x}$.
\label{lemma:cvxJn}
\end{lemma}

Lemma~\ref{lemma:cvxJn} allows us to apply the Subgradient
Algorithm~\cite{SMR:99,TL-MP-AS:03,DPB-AN-AO:03a} to obtain
$\vect{x}_{n}^{\epsilon_2}$ via $\{\vect{x}^{(r)}\}_{r=r_n}^{r_{n+1}}$
and the following lemma.
\begin{lemma}[Easy estimate of the $\epsilon$-subgradients
  of ${J}_{n}(\vect{x})$]
  Let the tolerance $\epsilon_1$ and time period $n$ be given. For any
  decision
  $\vect{x}^{(r)}$,
  we denote an $\epsilon_1$-optimal solution and
  $\epsilon_1$-worst-case distribution of~\eqref{eq:convJn} by
  $\vect{y}^{\epsilon_1}$ and
  $\mathbb{Q}_n^{\epsilon_1}(\vect{x}^{(r)})$,
  respectively.
  Let us consider the function $\map{g_{n}^{r}}{\real^d}{\real^d}$,
  defined as
\begin{equation*}
  g_{n}^{r}(\vect{x}):=
  \frac{d}{d \vect{x}}{\mathbb{E}_{\mathbb{Q}_n^{\epsilon_1}(\vect{x}^{(r)})} [f(\vect{x},\xi)] }.
\end{equation*}
\reviseone{Denote} an $\epsilon$-subdifferential of
${J}_n(\vect{x})$ at $\vect{x}$, by $\partial_{\epsilon}
{J}_{n}(\vect{x})$.
Then, for all $\epsilon \geq \epsilon_1$ we have the following
\begin{equation*}
	g_{n}^{r}(\vect{x}^{(r)}) \in \partial_{\epsilon} {J}_{n}(\vect{x}^{(r)}),
\end{equation*}
or equivalently, for every $\vect{z} \in {\dom}\: {J}_{n}$ and
$\epsilon \geq \epsilon_1$, we have
\begin{equation*}
  {J}_{n}(\vect{z})\geq {J}_{n}(\vect{x}^{(r)}) + \trans{g_{n}^{r}(\vect{x}^{(r)})}(\vect{z}-\vect{x}^{(r)}) - \epsilon.
\end{equation*}

Moreover, for any $\tilde{\vect{x}} \in \real^d$, there exists $\eta>0$
such that for all $\epsilon \geq \eta$ the following relation
holds
\begin{equation*}
  g_{n}^{r}(\tilde{\vect{x}}) \in \partial_{\epsilon} {J}_{n}(\tilde{\vect{x}}).
\end{equation*}
\label{lemma:access_subg}
\end{lemma}
\reviseone{Note how Lemma~\ref{lemma:access_subg} employs the discrete
  distribution $\mathbb{Q}_n^{\epsilon_1}$ generated from
  the \CGAaco in the computation of an
  $\epsilon$-subgradient function of ${J}_{n}$.}
Thus, the Subgradient Algorithm can be employed to reach an
$\epsilon_1$-proper data-driven decision with a lower certificate.

To do this, we make use of the scaled
\reviseone{$\epsilon$-subgradient} direction \reviseone{for the update
  of decisions $\{\vect{x}^{(r)}\}_{r=r_n}^{r_{n+1}}$, as follows}
\begin{equation}
\begin{aligned}
  \vect{x}^{(r+1)}=\vect{x}^{(r)}-\alpha^{(r)}
  \frac{g_{n}^{r}(\vect{x}^{(r)})}{\max\{
    \Norm{g_{n}^{r}(\vect{x}^{(r)})} \; ,
    \; 1 \} },
\label{eq:subg}
\end{aligned}
\end{equation}
where \reviseone{the nonnegative step size rule $\{ \alpha^{(r)}\}_r$
  is determined in advance. Later in the next subsection we will see
  how the choice of a step size rule affects the convergence of
  the Subgradient
  Algorithm to an $\vect{x}_{n}^{\epsilon_2}$.}

The Subgradient Algorithm requires access of $\{
g_{n}^{r}\}_{r=r_n}^{r_{n+1}}$, which are obtained from \CGAaco.  To
reduce the number of computations, we estimate the candidate
subgradient functions $\{ g_{n}^{r}\}_{r=r_n}^{r_{n+1}}$ as follows.
Let $\subscr{\epsilon}{SA} \geq \epsilon_1$ be a specified
tolerance. At some iteration $r \geq r_n$, assume that an
$\epsilon_1$-optimizer $\vect{y}^{\epsilon_1}$ and
$\epsilon_1$-worst-case distribution
$\mathbb{Q}_{n}^{\epsilon_1}(\vect{x}^{(r)})$ are obtained from
the \CGAaco. Using $\mathbb{Q}_{n}^{\epsilon_1}(\vect{x}^{(r)})$, we
calculate the function $g_{n}^{r}$ at $\vect{x}^{(r)}$ and perform the
subgradient iteration~\eqref{eq:subg}. At iteration $r+1$ with
$\vect{x}^{(r+1)}$, we firstly check for the suboptimality of Problem
(P1$_{n}^{(r+1)}$) using the initial candidate optimizer
$\vect{y}^{(0)}:=\vect{y}^{\epsilon_1}$ in \reviseone{the Point Search
  Algorithm}.
If the optimality gap $\eta^{(1)}$ is less than
$\subscr{\epsilon}{SA}$, we estimate the candidate subgradient
function $g_{n}^{r+1}$ using $g_{n}^{r}$ and proceed with the subgradient
iteration. Otherwise, we obtain $g_{n}^{r+1}$ from the \CGAaco, which
is again an ${\epsilon}_1$-subgradient function at $\vect{x}^{(r+1)}$.
Thus, we construct a sequence of
$\subscr{\epsilon}{SA}$-subgradient functions $\{
g_{n}^{r}\}_{r=r_n}^{r_{n+1}}$ that achieve an
$\vect{x}_{n}^{\epsilon_2}$ efficiently.
\begin{remark}[Effect of tolerance
    $\subscr{\epsilon}{SA}$]{\reviseone{\rm The tolerance $\subscr{\epsilon}{SA}$
      quantifies whether the function $g_{n}^{r}$ generated by the
      current worst-case distribution can also provide a good estimate
      of the ${\epsilon}$-subgradient at the next iteration point. If
      the function $g_{n}^{r}$ is an ${\epsilon}$-subgradient for
      ${\epsilon}$ small enough, there is no need of employing the \CGAaco
      to obtain a new subgradient function, which will be again an
      ${\epsilon_1}$-subgradient. In practice, we suggest to choose
      $\subscr{\epsilon}{SA} \gg \epsilon_1$ as it reduces the number
      of computations from the \CGAaco.}}
\end{remark}
\subsection{Convergence Analysis for the $\epsilon_2$-optimal decisions}
The following lemma follows from the convergence of the Subgradient
Algorithm applied to our problem scenario.

\begin{lemma}[Convergence of $\subscr{\epsilon}{SA}$-Subgradient
  Algorithm]
  \reviseone{For} each time period $n$ with an initial data-driven decision
  $\vect{x}^{(r_n)}$, assume that the subgradients defined in
  Lemma~\ref{lemma:access_subg}
  are uniformly bounded, i.e., there exists a constant $L>0$ such that
  $\Norm{ g_{n}^{r} } \leq L$ for all $r\geq r_{n}$.

  Given a predefined $\epsilon_2 >0$, let the certificate tolerance
  $\epsilon_1$ and the subgradient tolerance $\subscr{\epsilon}{SA}$
  be such that $0 < \epsilon_1 \leq \subscr{\epsilon}{SA} <
  \epsilon_2/ \mu$ \reviseone{with} $\mu := \max\{ L , \, 1 \}$. \reviseone{Let $\vect{x}_{n}^{\star} \in \argmin_{\vect{x}\in \real^d} {J}_{n}(\vect{x})$.} Then,
  there exists a large enough number $\bar{r}$, \reviseone{depending
    on $\subscr{\epsilon}{SA}$ and the step size rule $\{ \alpha^{(r)}\}_r$}, such that the
  above designed Subgradient Algorithm in~\eqref{eq:subg}
  has the following performance bounds
  \begin{equation*}
    \min\limits_{k\in \untilinterval{r_{n}}{r\reviseone{+r_{n}}}} \{
    {J}_{n}(\vect{x}^{(k)})\} -
    {J}_{n}(\vect{x}_{n}^{\star}) \leq \epsilon_2, \quad \forall
    \; r \geq  \bar{r},
  \end{equation*}
  and terminates at the iteration $r_{n+1}:=\bar{r}\reviseone{+r_{n}}$ with an
  $\epsilon_2$-optimal decision
  \reviseone{by choosing}
  $\vect{x}_{n}^{\epsilon_2} \in \argmin \limits_{k\in
    \untilinterval{r_{n}}{\bar{r}\reviseone{+r_{n}}} } \{
  {J}_{n}(\vect{x}^{(k)})\}$.
  \reviseone{In particular, there exists a large enough parameter $M$
    such that we can select 1) a constant step-size rule given by
  \begin{equation*}
    \alpha^{(i)}:=\frac{M}{\sqrt{\bar{r}+1}}, \; \forall i \in \untilinterval{r_{n}}{\bar{r}+r_{n}},
  \end{equation*}
  where $\bar{r}:=M^{2}\left(
    \frac{\epsilon_2}{\mu}-\subscr{\epsilon}{SA}\right)^{-2}$; or 2)
  a divergent, but square-summable, step-size rule given by
\begin{equation*}
  \alpha^{(i)}:=\frac{M}{i-r_{n}+1}, \; \forall i \in \untilinterval{r_{n}}{\bar{r}+r_{n}},
\end{equation*}
where $\bar{r}= \min \setdef{r\in \mathbb{N}}{M(3-\frac{1}{r+1})\leq
  2(\frac{\epsilon_2}{\mu} - \subscr{\epsilon}{SA})\ln(r+1)}$.  }
\label{lemma:convergeX}
\end{lemma}
\reviseone{In other words, Lemma~\ref{lemma:convergeX} specifies that
  there is a finite, large enough iteration step at which the
  $\subscr{\epsilon}{SA}$-Subgradient Algorithm terminates using
  the estimated $\subscr{\epsilon}{SA}$-subgradient functions.
}
To quantify the effect of the subgradient estimation on
the convergence rate under $\Xi_n$, we have the following
theorem.
\begin{theorem}[Worst-case computational bound for an
  $\vect{x}_{n}^{\epsilon_2}$] \label{thm:epsSA} \reviseone{For} each
  time period $n$ with an initial $\vect{x}^{(r_n)}$, let us consider
  the algorithm setting as in
  Lemma~\ref{lemma:convergeX}. Then, there exist parameters $\kappa
  \in (0,1)$, and $t > \epsilon_1$ such that the
  computational steps $\varphi(n,\bar{r})$ to reach
  $\vect{x}_{n}^{\epsilon_2}$ are bounded by
\begin{equation*}
  \varphi(n,\bar{r}) \leq \phi(n)+ \bar{r} \left( {\log}_{\kappa}(\frac{\epsilon_1}{t}) +1 \right)  ,
\end{equation*}
where $\bar{r}$ are the subgradient steps of
Lemma~\ref{lemma:convergeX}. The value $\phi(n)$
is the worst-case computational bound as in
Theorem~\ref{thm:ConvergeJhat} and one should use $\bar{\phi}(1)$ in
the bound in place of $\phi(n)$ if considering a data-streaming
scenario.
\end{theorem}
\reviseone{Theorem~\ref{thm:epsSA} integrates together the obtained
  bounds for the Subgradient Algorithm as well as the \CGAaco. As the
  result, a worst-case computational bound for the \ODAA, executed on
  the decision-update rate (the second  time scale in
  $(r)$),
  is related to the data-streaming rate. Whenever the data-streaming
  rate is greater than the worst-case bound, the \ODAA provides an
  $\vect{x}_{n}^{\epsilon_2}$ decision together with its estimated
  certificate ${J}_n^{\epsilon_1}(\vect{x}_n^{\epsilon_2})$ for each
  $n$.
     }

     From the Subgradient Algorithm, we provide, for each $\Xi_n$, a
     sequence $\{\vect{x}^{(r)}\}_{r=r_n}^{r_{n+1}}$ that approaches
     an $\vect{x}_{n}^{\epsilon_2}$.
     If \reviseone{the new data set} ${\Xi}_{n+1}$  \reviseone{is received} before reaching
     $\vect{x}_{n}^{\epsilon_2}$, we initialize the next sub-sequence
     \reviseone{obtained by applying the Subgradient Algorithm,}
     using the best decision at current iteration $r$, i.e.,
     $\vect{x}^{(r_{n+1})}:=\supscr{\vect{x}}{best}_{n} \in
     \argmin_{k\in
       \untilinterval{r_{n}}{r}}\{{J}_{n}^{\epsilon_1}(\vect{x}^{(k)})\}$. Then
     by connecting these sequences over $n$, our goal is achieved.

\section{Data Assimilation \reviseone{via \ODAA}} \label{sec:DataAssimilate}

This section summarizes and analyzes our \ODAAfull (\ODAA) for online
data \reviseone{sets} $\{{\Xi}_{n}\}_{n=1}^{N}$.  Specifically, we present the
algorithm procedure, its transient behavior and the convergence
result.

The \ODAA starts from some random initial decision $\vect{x}^{(r)} \in
\real^d$ and a data set ${\Xi}_{n}$, with $r=1$ and $n=1$.  Then,
\reviseone{it first} generates the certificate
${J}_{n}^{\epsilon_1}(\vect{x}^{(r)})$ via
\CGAaco,
then \reviseone{it} executes the Subgradient Algorithm to \reviseone{obtain the} decisions $\{
\vect{x}^{(r+1)}, \vect{x}^{(r+2)}, \ldots \}$ with lower and lower
certificates $\{ {J}_{n}^{\epsilon_1}(\vect{x}^{(r+1)}),
{J}_{n}^{\epsilon_1}(\vect{x}^{(r+2)}), \ldots\}$.  This algorithm has
the anytime property, meaning that the performance guarantee is
provided anytime, as soon as the first $\epsilon_1$-proper data-driven
decision with certificate ${J}_{n}^{\epsilon_1}(\vect{x}^{(r)})$ is
found. If no new data set ${\Xi}_{n+1}$ comes in, the algorithm
terminates as soon as the Subgradient Algorithm terminates at
iteration
$r_{n+1}$.
Otherwise, the algorithm \reviseone{resets the \CGAaco and the
  Subgradient Algorithm to update the} decision using more
data. \reviseone{This} achieves lower certificates with higher
confidence until we obtain the lowest possible certificate and
guarantee the performance almost surely. The details of the whole
algorithm procedure are summarized in the \reviseone{table of}
\ODAA.

{\begin{algorithm}[tbp]
\floatname{algorithm}{\ODAA}{}
 \caption{}
\label{alg:full}
\begin{algorithmic}[1]
\Require{Goes to Step~\ref{step:stream} \reviseone{upon data arrival, i.e.} $\Xi_{n} \leftarrow \Xi_{n+1}$.}
\State Set $\epsilon_1$, $\epsilon_2$, $\subscr{\epsilon}{SA}$, $\Xi_1$, $\vect{x}^{(0)} \in \real^d$, $\vect{y}^{(0)}=\vectorzeros{m}$ and $I_{1}^{(0)}=\O$;
\State $n \leftarrow 1$, $r \leftarrow 1$;
\State $r_n \leftarrow r$; \label{step:stream}
\State $({J}_n^{\epsilon_1}(\vect{x}^{(r)})$, ${\vect{y}}^{\epsilon_1}$, $\mathbb{Q}_n^{\epsilon_1}(\vect{x}^{(r)})) \leftarrow$ \CGAaco; \label{step:CG}
\Repeat
\State $\vect{x}^{(r+1)} \leftarrow  (\vect{x}^{(r)}, g_{n}^{r})$ as in~\eqref{eq:subg}, $r \leftarrow r+1$;
\State $\eta \leftarrow$ Point Search Algorithm;
\If {$\eta > \subscr{\epsilon}{SA}$, }
\State Goes to Step~\ref{step:CG};
\Else
\State Update $g_{n}^{r} \leftarrow g_{n}^{r-1}$;
\If{${J}_{n}^{\epsilon_1}(\vect{x}^{(r)}) <{J}_{n}^{\epsilon_1}(\supscr{\vect{x}}{best}_{n})$, }
\State Update and post ($\supscr{\vect{x}}{best}_{n}$, ${J}_{n}^{\epsilon_1}(\supscr{\vect{x}}{best}_{n})$);
\EndIf
\EndIf
\Until{$\Norm{ \vect{x}^{(r)} -\vect{x}^{(r-1)}}_2 < {\epsilon_2}$;}
\State $r_{n+1} \leftarrow r$;
\State Post $\vect{x}_n^{\epsilon_2}:=\supscr{\vect{x}}{best}_{n}$, ${J}_n^{\epsilon_1}(\vect{x}_n^{\epsilon_2}):={J}_{n}^{\epsilon_1}(\supscr{\vect{x}}{best}_{n})$;
\State Wait for $\Xi_{n+1}$, or Termination if $n=n_0$.
\end{algorithmic}
\end{algorithm}}

The transient behavior of the \ODAA is affected by the data-streaming
rate and the rate of convergence of the intermediate algorithms
(\reviseone{decision-update rate and certificate-update
  rate}).
To further describe these effects \reviseone{in each time period $n$,
  we say that the data-streaming rate is \textit{slow with respect to
    the decision-update rate},} if we can find an
$\vect{x}_n^{\epsilon_2}$ \reviseone{via the \ODAA, where the
  worst-case scenario} is described in
Theorem~\ref{thm:epsSA}. Further, we call \reviseone{it \textit{slow
    with respect to the certificate-update rate},} if we can find at
least one certificate during \reviseone{this time period, where the
  worst-case scenario} is described in
Theorem~\ref{thm:ConvergeJhat}. When the data-streaming rate is
\reviseone{slow w.r.t.~the decision-update rate for all time periods,}
the \ODAA guarantees to find
$\{\vect{x}_n^{\epsilon_2}\}_{n=1}^{N}$. When the
\reviseone{data-streaming rate is slow w.r.t.~the decision-update
  rate} for at least one time period \reviseone{$n_0$}, it guarantees
to find an \reviseone{$\vect{x}_{n_0}^{\epsilon_2}$}. When the
data-streaming rate is \reviseone{slow w.r.t.~the certificate-update
  rate} for at least one time period \reviseone{$n_0$}, the \ODAA
guarantees to find a \reviseone{${J}_{n_0}^{\epsilon_1}$} for an
$\vect{x}^{(r)}$ \reviseone{and $r \ge r_{n_0}$}. When the
\reviseone{data-streaming rate is \textit{not slow w.r.t.~the
    certificate-update rate} for any} time period, the \ODAA will hold
on the newly streamed data set, to make the data-streaming rate
\reviseone{\textit{slow w.r.t.~the decision-update rate}} and
achieve a better data-driven decision efficiently.

Next, we state the convergence result of the \ODAA when the data
streams are \reviseone{slow w.r.t.~the decision-update rate}
for all  time periods.

\begin{theorem}[Finite convergence of the \ODAA]
  \reviseone{Consider tolerances} $\epsilon_1$, $\epsilon_2>0$ and
  \reviseone{streaming data sets $\{{\Xi}_{n}\}_{n=1}^{N}$ with $N <
    \infty$ for a decision making problem~\eqref{eq:P}. Assume that the
    data streams are slow w.r.t.~the to decision-update rate
    for all $n$, i.e., assume the
    length of each time period $n$ is no shorter than
    $\varphi(1,\bar{r})$, where $\varphi$ and $\bar{r}$ are
    described as in Theorem~\ref{thm:epsSA} and
    Lemma~\ref{lemma:convergeX}, respectively.} Then, the \ODAA
  guarantees to find a sequence of $\epsilon_2$-optimal
  $\epsilon_1$-proper data-driven decisions
  $\{\vect{x}_n^{\epsilon_2}\}_{n=1}^{N}$ associated with the sequence
  of the certificates $\{{J}_n^{\epsilon_1}(\vect{x}_n^{\epsilon_2})
  \}_{n=1}^{N}$ so that the performance
  guarantee~\eqref{eq:epsiperfgua} holds for all
  $n$. \reviseone{Furthermore, the values of these certificates are
    guaranteed to be low in high probability. That is, for each $n$
	\begin{equation}
          \mathbf{P}^{n} ( {J}_{n}^{\epsilon_1}(\vect{x}_n^{\epsilon_2}) \leq  J^{\star} + \epsilon_1 + \epsilon_2 + 2\hat{L} \epsilon(\beta_{n}) ) \geq 1- \beta_{n}
	\label{eq:qualityofJ}
\end{equation}
holds, where $J^{\star}:= \inf\limits_{\vect{x} \in
  \real^d}{\mathbb{E}_{\prob} [f(\vect{x},\xi)] }$ is the optimal
objective value for the original problem~\eqref{eq:P}, the parameter
$\hat{L}$ depends on steepness of the function $f$, and the parameter
$\epsilon(\beta_{n})$ is determined as in Lemma~\ref{eq:cert}.  }

In addition, given any tolerance $\epsilon_3$\reviseone{, data
  stream that is slow w.r.t.~the decision-update rate for all $n
  \in \until{N}$ }
with $N \rightarrow \infty$, \reviseone{and $\sum_{n=1}^{\infty}
\beta_{n} <\infty$,} there exists a large enough number
$n_0(\epsilon_3) >0$, such that the algorithm terminates in finite
time with a guaranteed $\epsilon_2$-optimal and $\epsilon_1$-proper
data-driven decision $\vect{x}_{n_0}^{\epsilon_2}$ and a certificate
${J}_{n_0}^{\epsilon_1}(\vect{x}_{n_0}^{\epsilon_2})$ such that the
performance guarantee holds almost surely. \reviseone{That is,}
  \begin{equation}
    \begin{aligned}
      \mathbf{P}^{n_0}({\mathbb{E}_{\prob}
        [f(\vect{x}_{n_0}^{\epsilon_2},\xi)] } \leq
      {J}_{n_0}^{\epsilon_1}(\vect{x}_{n_0}^{\epsilon_2}) +
      \epsilon_1)=1,
    \end{aligned} \label{eq:asperf}
  \end{equation}
  and meanwhile the quality of the designed certificate
  ${J}_{n_0}^{\epsilon_1}(\vect{x}_{n_0}^{\epsilon_2})$ is
  guaranteed. \reviseone{In other words,} for all the rest of the data sets
  $\{{\Xi}_{n}\}_{n=n_0}^{\infty}$, any element in the desired
  certificate sequence
  $\{{J}_n^{\epsilon_1}(\vect{x}_n^{\epsilon_2})
  \}_{n=n_0}^{\infty}$ satisfies
  \begin{equation}
    \begin{aligned}
      \sup_{n \geq n_0}
      {J}_{n}^{\epsilon_1}(\vect{x}_{n}^{\epsilon_2})
      \leq J^{\star} + \epsilon_1 +\epsilon_2 +\epsilon_3.
    \end{aligned}
		\label{eq:perfbound}
  \end{equation}
\label{thm:assimilation}
\end{theorem}
\reviseone{Theorem~\ref{thm:assimilation} quantifies the goodness of
  the certificates that are achievable via the \ODAA, under the
  condition that the data-streaming rate is slow w.r.t.~the
  decision-update rate.
  Intuitively, the smaller the tolerances $\epsilon_i$ are, the lower
  the certificates become. Further, as $N\rightarrow \infty$, the
  smaller the parameter $\epsilon(\beta_{N})$ is and the higher the
  confidence $1 - \beta_N \rightarrow 1$. When infinitely many data
  are streamed in, the theorem implies that we can get arbitrarily
  close to the optimal decision with probability one.}
  \begin{remark}[Selection of tolerances $\epsilon_2$, $\epsilon_1$ and
    $\subscr{\epsilon}{SA}$] { \reviseone{\rm In practice, the tolerance
      $\epsilon_2$ determines the performance bound of ${J}_{n}$,
      which governs the whole algorithm. With a given $\epsilon_2$,
      tolerance $\epsilon_1$ and $\subscr{\epsilon}{SA}$ can be chosen
      following the rule in Lemma~\ref{lemma:convergeX}. Intuitively,
      $\epsilon_1$ can be chosen to be two orders of magnitude smaller
      than $\epsilon_2$, while $\subscr{\epsilon}{SA}$ can be an order
      of magnitude smaller than $\epsilon_2$. These tolerances can
      also be chosen in a data-driven fashion,
      to achieve asymptotic convergence, or a better
      transient behavior of the algorithm.}}
 \end{remark}
 \begin{remark}[Asymptotic behavior of the \ODAA] { \reviseone{\rm
     Theorem~\ref{thm:assimilation} claims the
     convergence of the \ODAA to a decision with a desired
     certificate using large but a finite data set. The smaller
     $\epsilon_3$ is, the larger data set is needed to achieve the
     desired certificate. Because tolerances $\epsilon_1$,$\epsilon_2$
     and $\epsilon_3$ can be chosen arbitrarily small, the certificate
     can indeed approach to $J^{\star}$. However,
     $\epsilon_1$,$\epsilon_2$ may affect the transient behavior of
     the algorithm and the data-streaming rate. In practice, to reach
     $J^{\star}$, these tolerances can be chosen in a data-driven
     fashion, for example diminishing sequences.}}
 \end{remark}
 \section{Data incremental covering
 }\label{sec:Heu}

 In this section, we aim to handle large streaming data sets for
 efficient \ODAAfull (\ODAA). To achieve this, we firstly propose an
 \ICA (\ICAaco). This algorithm leverages the pattern of the data
 points \reviseone{to obtain} a new ambiguity set, denoted by
 $\mathcal{\tilde{P}}_n$. Then, we adapt $\mathcal{\tilde{P}}_n$ for a
 variant of the \ODAA. The resulting algorithm enables us to construct
 subproblems which have a lower dimension than those generated without
 it, and we verify its capability of handling large data sets in
 simulation.

 \subsection{The \ICAaco} Let $\zeta$ and $\omega$ denote the center
 and radius of the Euclidean ball $B_{\omega}(\zeta)$,
 respectively. For each data set $\Xi_n$ and a given $\omega$, let
 $\mathcal{C}_n \subset \Xi_n$ denote the set of points such that
 $\Xi_n \subset \cup_{\zeta \in \mathcal{C}_n} B_{\omega}(\zeta)$. Let
 $p:=|\mathcal{C}_n|$ denote the number of these Euclidean balls. To
 account for the number of data points that are covered by a specific
 ball, we associate each ball $B_{\omega}(\zeta_k)$ a weighting
 parameter $\theta_k$.
 We denote by $\mathcal{Q}_n:=\{\theta_k\}_{k=1}^{p}$ the set of these
 parameters. Then, as data sets $\{{\Xi}_{n}\}_{n=1}^{N}$ are
 sequentially accessible, we are to incrementally cover data sets by
 adapting $\mathcal{C}_n$ and $\mathcal{Q}_n$.

 Formally, the \ICAaco works as follows. Let $\mathcal{C}_0=\O$ and
 $\mathcal{Q}_0=\O$. For the $\supscr{n}{th}$ time period with set
 $\Xi_{n}$, we initialize sets as $\mathcal{C}_n:=\mathcal{C}_{n-1}$
 and $\mathcal{Q}_n:=\mathcal{Q}_{n-1}$. To generate a random cover
 for $\Xi_{n}$, we randomly and sequentially evaluate each newly
 streamed data point. Let $\varsigma \in \Xi_{n} \setminus \Xi_{n-1}$
 denote the data point under consideration. If $\varsigma \notin
 B_{\omega}(\zeta_{k})$ for all $\zeta_{k} \in \mathcal{C}_n$, we
 update $\mathcal{C}_n \leftarrow \mathcal{C}_n \cup
 \{\zeta_{p+1}:=\varsigma\}$, $\mathcal{Q}_{n} \leftarrow
 \mathcal{Q}_{n} \cup \{ \theta_{p+1}:=1 \}$ and $p \leftarrow |\mathcal{C}_n|$.
 If $\varsigma$ is covered by some (at least one) Euclidean balls,
 i.e., $\varsigma \in B_{\omega}(\zeta_{k})$ for some $k$ with
 $\zeta_{k} \in \mathcal{C}_n$, we only update $\mathcal{Q}_{n}$. Let
 $\ell_{\varsigma}$ denote the number of the balls that cover
 $\varsigma$ and let $I_{\varsigma} \subset \until{p}$ denote the
 index set of these balls.
 Then we update elements of $\mathcal{Q}_{n}$ via $\theta_k \leftarrow
 \theta_k+\ell_{\varsigma}^{-1}$ for all $k \in I_{\varsigma}$.  After
 all the new data points have been evaluated in this way, we achieve a cover
 of $\Xi_n$. Then, as the data set streams over time, the algorithm
 incrementally updates the cover and weights. By construction, we see
 that $|\mathcal{C}_n| \leq n$.

 Next, we use $\mathcal{C}_n$ and $\mathcal{Q}_n$ to construct a new
 ambiguity set that results in potentially low dimensional subproblems
 in the \ODAA.

\subsection{Integration of \ICAaco into \ODAA}
Following the \ICAaco, we consider a distribution $\tilde{\prob}^n$
associated with ${\Xi}_{n}$, as follows
\begin{equation}
  \tilde{\prob}^n := \frac{1}{n}\sum_{k=1}^{p} \theta_k \delta_{\{\zeta_{k}\}},
\label{eq:weightedempmeas}
\end{equation}
where $\delta_{\{ \zeta_{k}\}}$ is a Dirac measure at the center of
the covering ball $B_{\omega}(\zeta_{k})$ and $\theta_k$ is the
associated weight of $B_{\omega}(\zeta_{k})$. We claim the
distribution $\tilde{\prob}^n$ is close to the empirical distribution
$\hat{\prob}^n$ under the Wasserstein metric, using the following
lemma.

\begin{lemma}[Distribution $\tilde{\prob}^n$ is a good estimate of
  $\hat{\prob}^n$]
  Let the radius $\omega$ of the Euclidean ball be chosen. Then the
  distribution $\tilde{\prob}^n$ constructed by the \ICAaco on
  ${\Xi}_{n}$ is close to $\hat{\prob}^n$ under the Wasserstein
  metric, i.e., $d_W(\hat{\prob}^n,\tilde{\prob}^n) \leq \omega$.
\label{lemma:closeof2P}
\end{lemma}

Then equipped with Lemma~\ref{lemma:closeof2P} and
Theorem~\ref{thm:MeasConc} on the measure of concentration result, we
can provide the certificate that ensures the performance guarantee
in~\eqref{eq:perfgua}.

\begin{lemma}[Tractable certificate generation for $\vect{x}$ with
  Performance Guarantee~\eqref{eq:perfgua} using $\tilde{\prob}^n$]
  Given ${\Xi}_{n}:=\{{\xi}_{k} \}_{k=1}^{n}$, $\beta_{n} \in (0,1)$,
  \reviseone{$\vect{x} \in
    \real^d$,}
  and the radius $\omega$ of the covering balls. Define the new
  ambiguity set
  $\mathcal{\tilde{P}}_n:=\mathbb{B}_{\tilde{\epsilon}(\beta_{n})}(\tilde{\prob}^n)$
  where the center of the Wasserstein ball $\tilde{\prob}^n$ is
  defined in~\eqref{eq:weightedempmeas} and the radius
  $\tilde{\epsilon}(\beta_{n}):=\epsilon(\beta_{n})
  +\omega$.
  Then the following certificate satisfies~\eqref{eq:perfgua} for all
 $\vect{x} \in \real^d$
\begin{equation}
  {J}_{n}(\vect{x})
  := \sup\limits_{\mathbb{Q}  \in \mathcal{\tilde{P}}_n} {\mathbb{E}_{\mathbb{Q}} [f(\vect{x},\xi)] }.
\label{eq:tildePcert}
\end{equation}
Further, under the same assumptions required in
Theorem~\ref{thm:cvxredu} we have the new version
of~\eqref{eq:convJn} as follows
\begin{equation}
\begin{aligned}
  {J}_{n}(\vect{x}):=
  \sup\limits_{\vect{y}_1,\ldots,\vect{y}_{p} \in \real^m}&
  \frac{1}{n}\sum_{k=1}^{p}
  \theta_k f(\vect{x},\zeta_{k}-\vect{y}_k) ,\\
  \st \quad &\frac{1}{n} \sum_{k=1}^{p} \theta_k \Norm{\vect{y}_k} \le
  \tilde{\epsilon}(\beta_{n}),
\end{aligned}
\label{eq:tildePconvJn}\tag{\small $\reviseone{\tilde{\rm{P1}}_{n}}$}
\end{equation}
and the associated worst-case distribution
$\tilde{\mathbb{Q}}_{n}^{\star}(\vect{x})$ is a weighted
version of $\mathbb{Q}_{n}^{\star}(\vect{x})$ in
Theorem~\ref{thm:cvxredu}, i.e.,
\begin{equation*}
  \tilde{\mathbb{Q}}_{n}^{\star}(\vect{x}) := \frac{1}{n}\sum_{k=1}^{p} \theta_k \delta_{\{\zeta_{k} - \vect{y}^{\star}_k\}},
\end{equation*}
where $\vect{y}^{\star}:=(\vect{y}^{\star}_1,\ldots,\vect{y}^{\star}_{p})$ is an
optimizer of~\eqref{eq:tildePconvJn}.
\label{lemma:tildePcertgener}
\end{lemma}

\begin{remark}[New version of~\eqref{eq:JoverSimplex}]
{\rm  The equivalent formulation of Problem~\eqref{eq:tildePconvJn}
   is a
  new version of~\eqref{eq:JoverSimplex}, defined as follows
\begin{equation}
\begin{aligned}
 {J}_{n}(\vect{x}):= \max\limits_{{
\substack{\vect{u}_1,\ldots,\vect{u}_p \in \real^m \\ \vect{v}_1,\ldots,\vect{v}_p \in \real^m}
 } }& \frac{1}{n}\sum_{k=1}^{p}
\reviseone{\tilde{h}_k}( \frac{\vect{u}_k - \vect{v}_k}{\theta_k}) ,\\
\st \quad &  (\vect{u}, \vect{v}) \in n \tilde{\epsilon}(\beta_{n}) {\Delta}_{2mp},
\end{aligned}
\label{eq:tildePJoverSimplex}\tag{\small $\reviseone{\tilde{\rm{P2}}_{n,p}}$}
\end{equation}

where for each $k \in \until{p}$, $\zeta_{k} \in
\mathcal{C}_n$ and $\vect{x} \in \real^d$, we define
$\map{\reviseone{\tilde{h}_k}}{\real^m}{\real}$ as
\begin{center} {$\reviseone{\tilde{h}_k}(\vect{y}):=\theta_k
    f(\vect{x},\zeta_{k}-\vect{y})$.} \end{center}
}
\end{remark}

With the constructed ambiguity set $\mathcal{\tilde{P}}_n$ and
certificate function ${J}_{n}$, the the
developed algorithms in Section~\ref{sec:CG} and
Section~\ref{sec:LowJ} are valid to solve
Problem~\eqref{eq:tildePJoverSimplex}. And the main
Theorem~\ref{thm:assimilation} on the finite convergence of the \ODAA
is valid for the certificate function ${J}_{n}$ where the only
difference is that the quality of the certificate for
$\vect{x}_{n}^{\epsilon_2}$ in~\eqref{eq:perfbound} is replaced
by
\begin{equation*}
    \begin{aligned}
      \sup_{n \geq n_0}
      {J}_{n}^{\epsilon_1}(\vect{x}_{n}^{\epsilon_2})
      \leq J^{\star} + \epsilon_1 +\epsilon_2 +\epsilon_3 + 2 (1-\frac{p_{n_0}}{n_0}) \hat{L}\omega,
    \end{aligned}
		\label{eq:tildePperfbound}
  \end{equation*}
  where $n_0$ is the number of the data set in ${\Xi}_{n_0}$ and
  $p_{n_0}$ indicates the number of Euclidean balls that cover
  ${\Xi}_{n_0}$.

\section{Simulation results}\label{sec:Sim}
In this section, we demonstrate the application of the
\reviseone{proposed algorithms on two case studies}, with a
potentially large streaming data set.

\subsection{Study 1: The Effect of the \ICAaco}
In order to visualize the effect of the \ICAaco, here we solve a toy
problem in form of~\eqref{eq:P} using \ODAA, with and without the
\ICAaco respectively.  Let $x \in \real$ be the variable for
Problem~\eqref{eq:P}. Assume there are $N=200$ data points
$\{{\xi}_{k}\}_{k=1}^{N}$ streaming into the algorithm. Assume each
time period is one second, and for each second $k$ we only stream in
one data point ${\xi}_{k} \in \real^{3}$, where ${\xi}_{k}$ is a
realization of the unknown distribution $\prob$. The $\prob$ we use
for simulation is a multivariate weighted Gaussian mixture
distribution with three centers, where each center has mean
$\mu_1=(2,-4,3)$, $\mu_2=(-3,5,0)$, $\mu_3=(0, 0,-6)$, variance
$\sum_1=\rm{diag}(1,3,2)$, $\sum_2=2 \cdot \vect{I}_3$,
$\sum_3=\vect{I}_3$, and weights 0.25, 0.5, 0.25, respectively. Let
the cost function to be $f({x},\xi):= x^2 -\trans{\xi}\xi$, the
confidence be $1-\beta_{n}:=1-0.95e^{1-\sqrt{n}}$ and we use the
parameter $c_1=2$, $c_1=1$ to design the radius $\epsilon(\beta_{n})$
of the Wasserstein ball in~\eqref{eq:epsirad}. The radius of the
Euclidean ball for the \ICAaco is $w=1.5$. We sample the initial
decision $x^{(0)}$ from the uniform distribution $[-10,10]$. The
tolerance for the algorithm is $\epsilon_1=10^{-5}$,
\reviseone{$\epsilon_2 = 10^{-4}$}.

\begin{figure*}[htp]
  \subfloat[\footnotesize{ \revisetwo{Data set cover at
    $\supscr{200}{th}$ time period. }}]{
	\begin{minipage}[c][0.7\width]{
	   0.32\textwidth}
	   \centering
	   \includegraphics[width=1\textwidth]{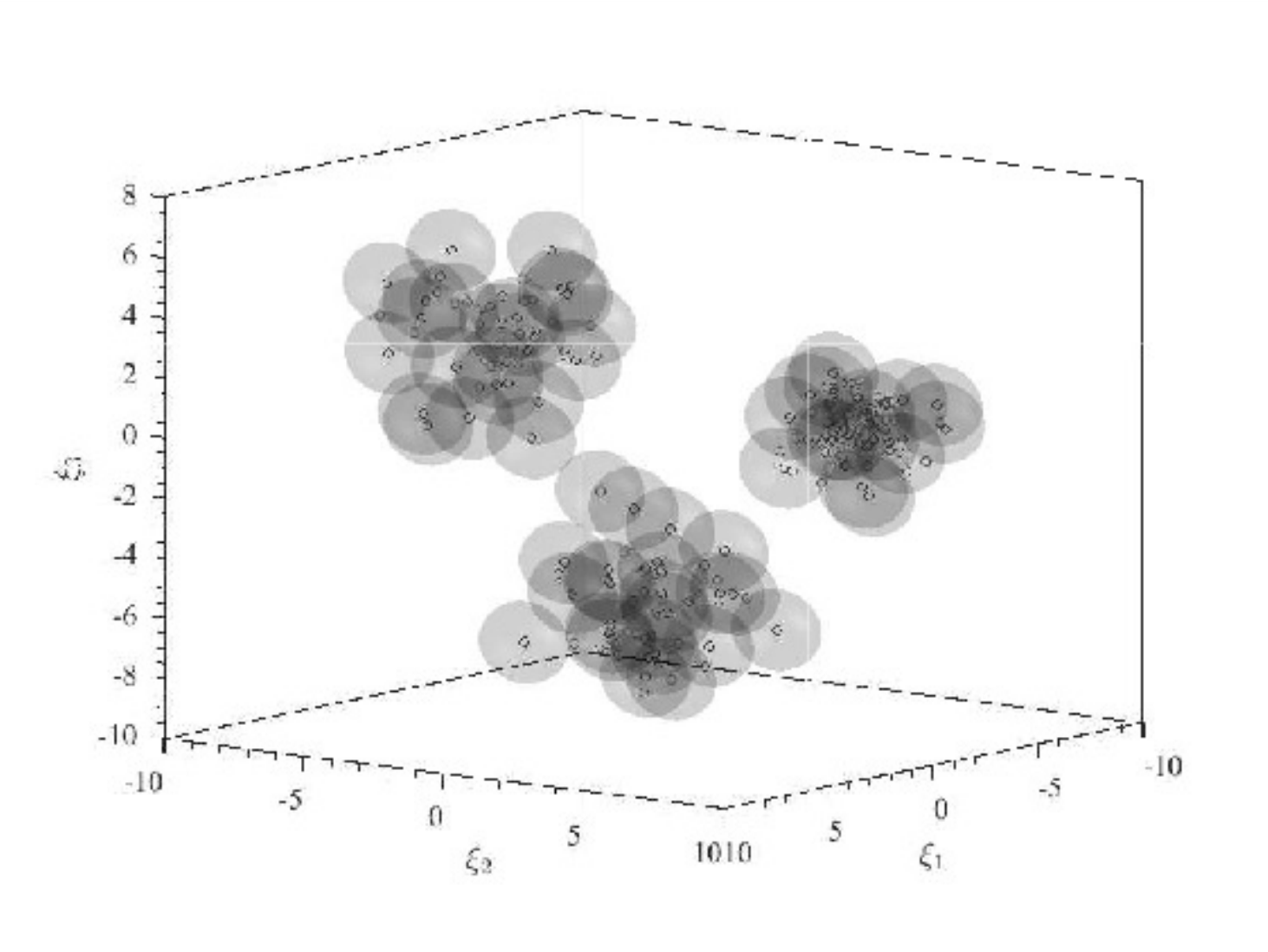}
     \label{fig:Cover}
	\end{minipage}}
 \hfill
  \subfloat[\footnotesize{\revisetwo{Size of $\tilde{\prob}^n$ ($\hat{\prob}^n$) and count of $\tilde{\rm{CP}}_n^{(l)}$ (${\rm{CP}}_n^{(l)}$)
  solved with(out) \ICAaco. }}]{
	\begin{minipage}[c][0.7\width]{
	   0.32\textwidth}
	   \centering
	   \includegraphics[width=0.92\textwidth]{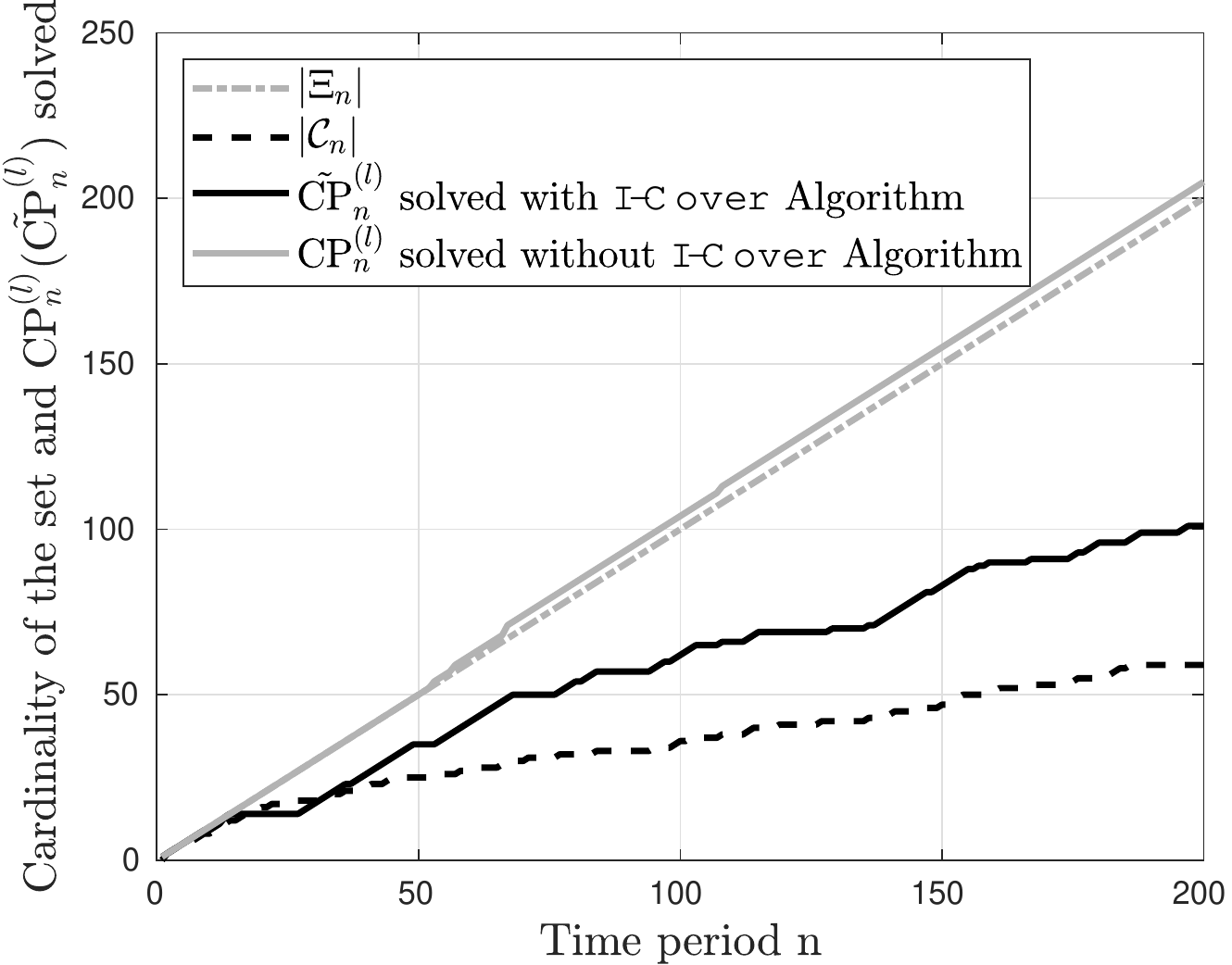}
     \label{fig:CP}
	\end{minipage}}
 \hfill
  \subfloat[\footnotesize{\revisetwo{Relative error of ${J}^{\epsilon_1}_{n}(\vect{x}^{(r)})$ over time, with the
    \ICAaco. } }]{
	\begin{minipage}[c][0.7\width]{
	   0.32\textwidth}
	   \centering
	   \includegraphics[width=0.98\textwidth]{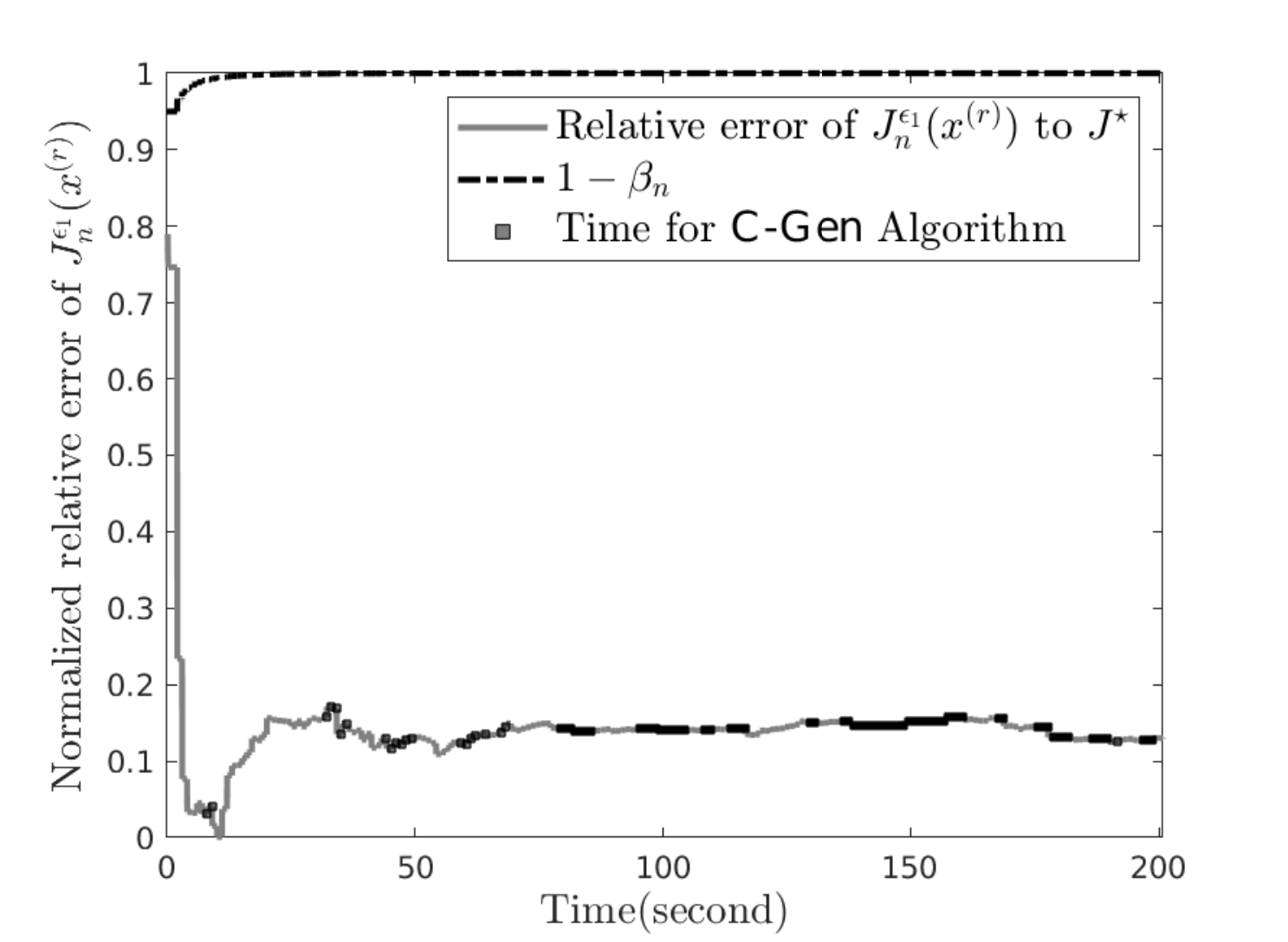}
      \label{fig:CertICA}
	\end{minipage}}  \\
  \subfloat[\footnotesize{\revisetwo{Relative error of ${J}^{\epsilon_1}_{n}(\vect{x}^{(r)})$ over time, without the \ICAaco. } }]{
	\begin{minipage}[c][0.7\width]{
	   0.32\textwidth}
	   \centering
	   \includegraphics[width=0.98\textwidth]{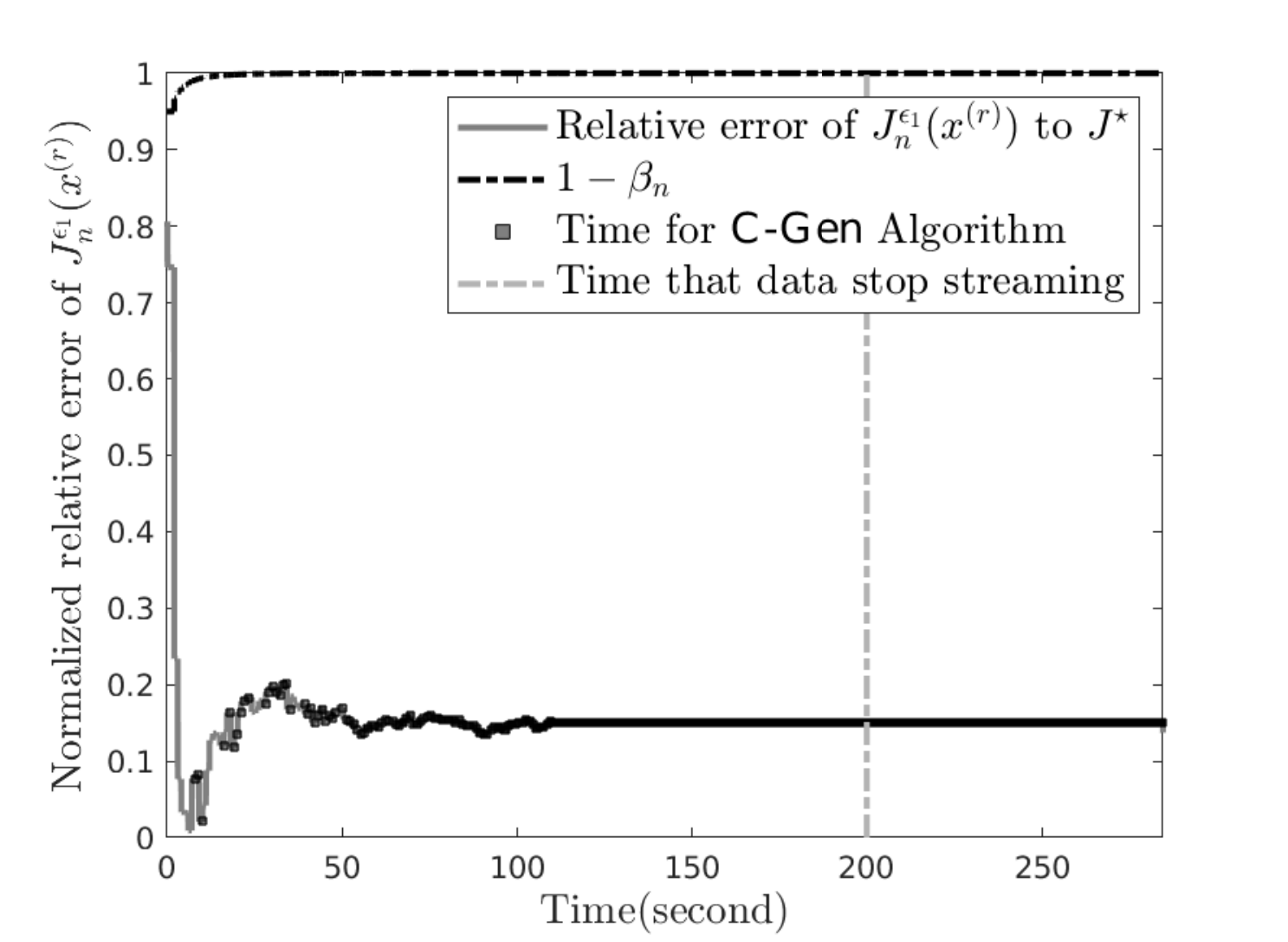}
      \label{fig:CertnoICA}
	\end{minipage}}
 \hfill
  \subfloat[\footnotesize{\revisetwo{Size of $\hat{\prob}^n$ ($\tilde{\prob}^n$) indicated by $|\Xi_n|$ ($|\mathcal{C}_n|$) and count of $\tilde{\rm{CP}}_n^{(l)}$ solved. }}]{
	\begin{minipage}[c][0.7\width]{
	   0.32\textwidth}
	   \centering
	   \includegraphics[width=0.95\textwidth]{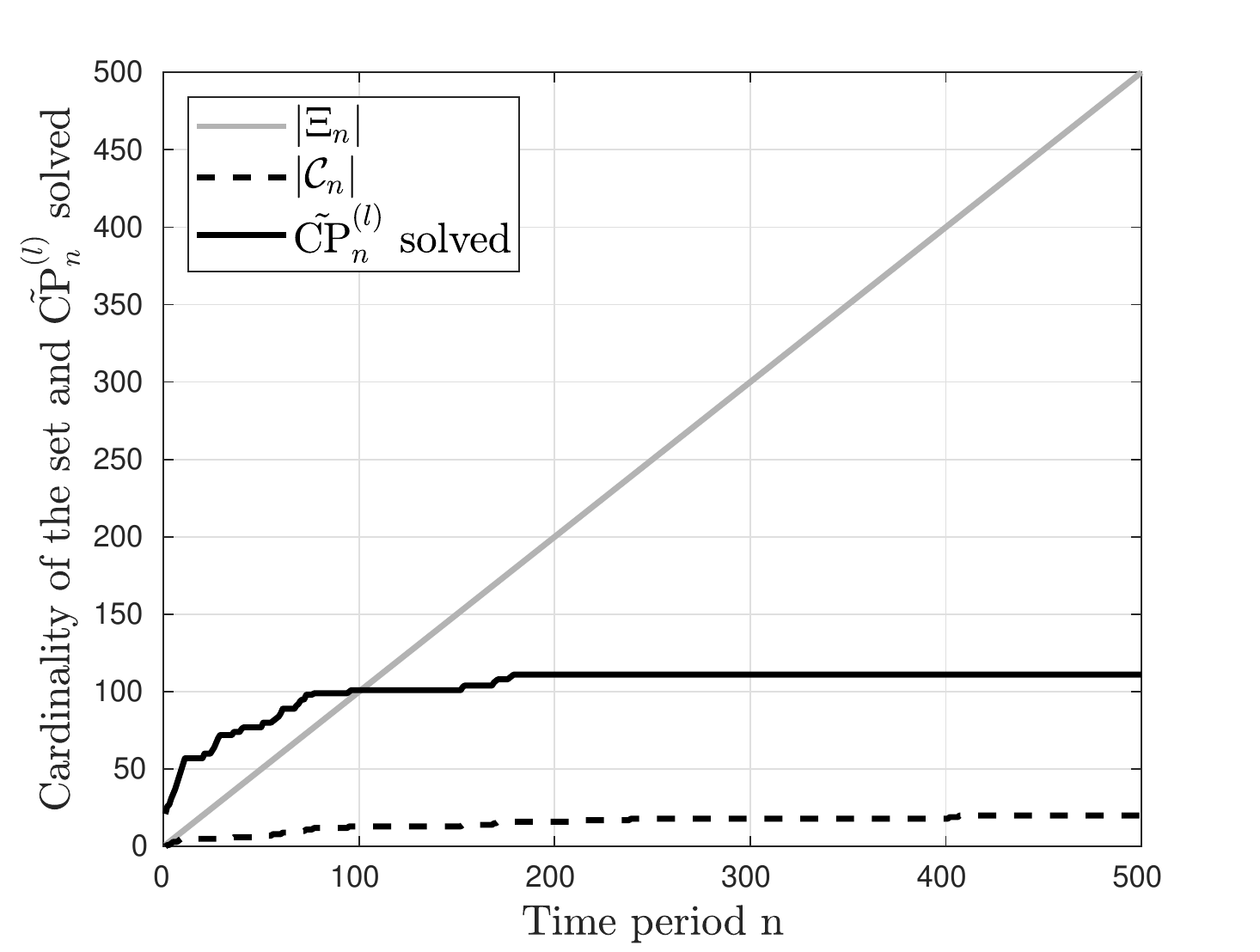}
     \label{fig:CP500}
	\end{minipage}}
 \hfill
  \subfloat[\footnotesize{\revisetwo{Relative error of ${J}^{\epsilon_1}_{n}(\vect{x}^{(r)})$ over time, with the
    \ICAaco. } }]{
	\begin{minipage}[c][0.7\width]{
	   0.32\textwidth}
	   \centering
	   \includegraphics[width=0.98\textwidth]{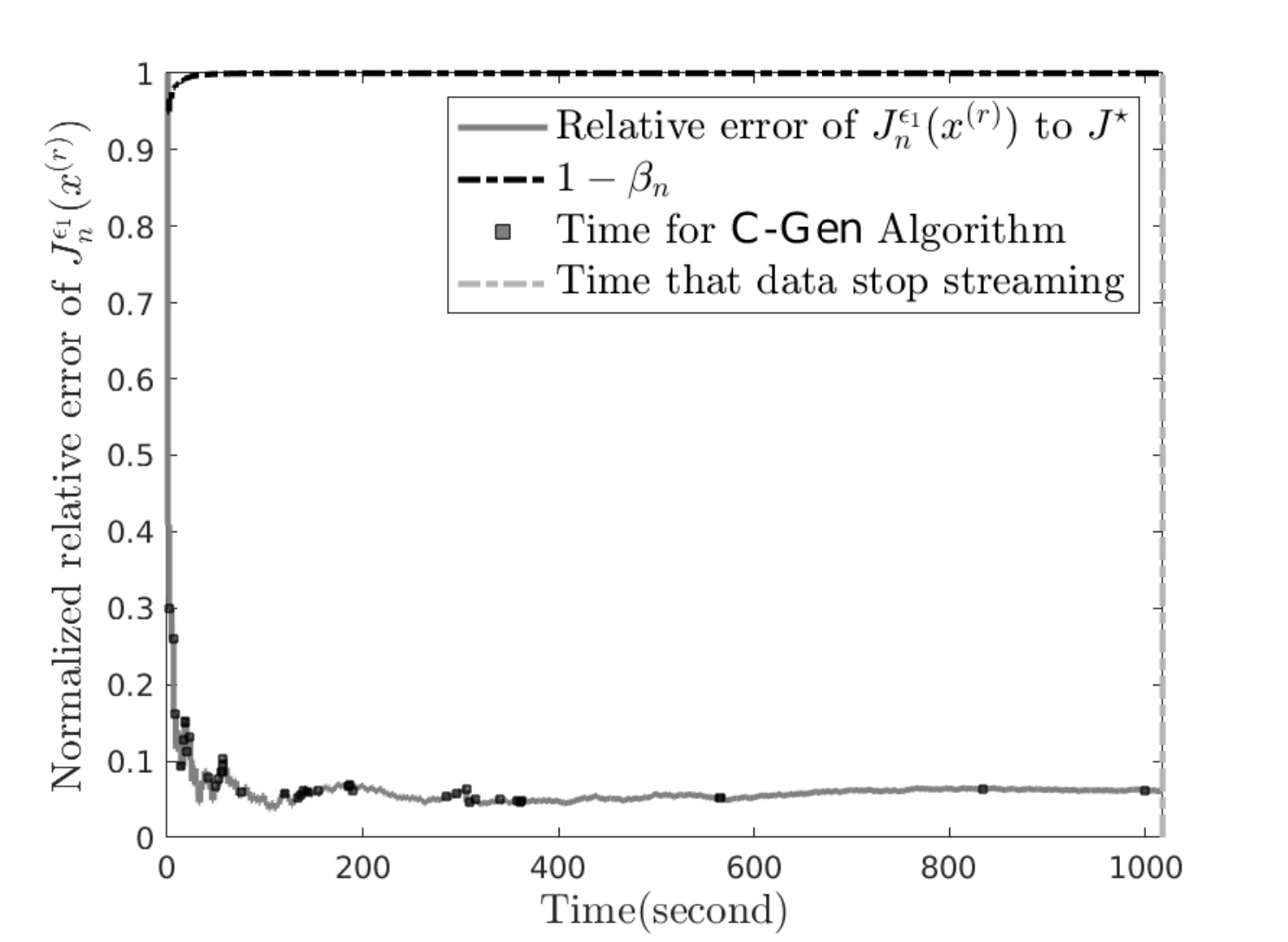}
      \label{fig:CertICA500}
	\end{minipage}}
 \caption{ \footnotesize \revisetwo{ Simulation results of the \ODAAfull, with and without the \ICA}}
\end{figure*}

Figure~\ref{fig:Cover} and Figure~\ref{fig:CP} demonstrate the effect
of the \ICAaco in the \ODAA. Specifically, Figure~\ref{fig:Cover} shows the
incremental data covering at the end of the $\supscr{200}{th}$ time
period in the $(\xi_1,\xi_2,\xi_3)$ coordinates. The large shaded
area are $59$ Euclidean balls $B_{\omega}$ with their centers
$\{\zeta_{k}\}_{k=1}^{59}$ denoted by some of the
small circles, where
all these small circles constitute the streamed data set
${\Xi}_{200}:=\{{\xi}_{k}\}_{k=1}^{200}$. In Figure~\ref{fig:CP}, the
\reviseone{gray} dashed line represents the number of the data points used as
centers of the empirical distribution $\hat{\prob}^n$ over time and
the black dashed line is that for distribution
$\tilde{\prob}^n$. Clearly as the data streams over time, the number
$p:=|\mathcal{C}_n|$ is significantly smaller than $n:=|{\Xi}_n|$,
which results in the size of Problem~\eqref{eq:tildePJoverSimplex}
\reviseone{being} much smaller than that
of~\eqref{eq:JoverSimplex}. Further, the \reviseone{gray solid}
line
counts the total number of subproblems (CP$_{n}^{(l)}$) solved to
generate certificates over time and the black \reviseone{solid} line
represents that for subproblems ($\tilde{\rm{CP}}_{n}^{(l)}$) \reviseone{in solution to~\eqref{eq:tildePJoverSimplex}}. These
subproblems search the explicit solution for the
$\epsilon_1$-worst-case distribution and consume the major computing
resources in the \ODAA. It can be seen that the number of
($\tilde{\rm{CP}}_{n}^{(l)}$) solved over time is on average only half
of the (CP$_{n}^{(l)}$) in each time period. Together, the dimension
and total number of subproblems ($\tilde{\rm{CP}}_{n}^{(l)}$) solved
with the \ICAaco is significantly smaller than \reviseone{that} without~it.

To evaluate the quality of the obtained $\epsilon_1$-proper
data-driven decision with the streaming data, we estimate the
optimizer of~\eqref{eq:P}, $\vect{x}^{\star}$, by minimizing the
average value of the cost function $f$ for a validation data set of
$\subscr{N}{val}=10^4$ data points randomly generated from the
distribution $\prob$ (in the simulation case $\prob$ is known). We
take the resulting objective value as the estimated optimal objective
value for Problem~\eqref{eq:P}, i.e.,
$J^{\star}:={J}^{\star}({\vect{x}}^{\star})$. We calculate
${J}^{\star}({\vect{x}}^{\star})$ using the underline distribution
$\prob$, serving as the true but unknown scale to evaluate the
goodness of the certificate obtained throughout the algorithm.

Figure~\ref{fig:CertICA} and Figure~\ref{fig:CertnoICA} show the
evolution of the certificate sequence
$\{{J}_n^{\epsilon_1}(\vect{x}^{(r)}) \}_{n=1,r=1}^{N,\infty}$
with the \ICAaco and \reviseone{that} without the \ICAaco, respectively. Here, the
optimal decision of~\eqref{eq:P} is trivially ${\vect{x}}^{\star}=0$,
and for both algorithms the subgradient counterpart of the \ODAA
returns the optimal decision after the first data point ${\xi}_{1}$ is
used.  Therefore, after a very short period within the first second,
both figures start reflecting the certificate evolution under the
decision sequence $\{\vect{x}^{(r)} \approx 0\}_{r=r_2}^{\infty}$.
The \reviseone{gray solid} line in both Figure~\ref{fig:CertICA} and
Figure~\ref{fig:CertnoICA} show the relative goodness of the
certificates for the currently used $\epsilon_1$-proper data-driven
decision $\vect{x}^{(r)} \approx 0$ calibrated by the estimated
optimal value $J^{\star}$ over time. The \reviseone{black} segments on the \reviseone{gray solid} line
indicate that \reviseone{the \CGAaco is executing for certificates update,}
while at these time intervals the old
certificate ${J}_{n}^{\epsilon_1}(\vect{x}_{n}^{\epsilon_2})$,
associated with the $\epsilon_2$-optimal and $\epsilon_1$-proper
data-driven decision $\vect{x}_{n}^{\epsilon_2}$, is still valid to
guarantee the performance under the old confidence $1-\beta_{n}$. This
situation commonly happens when a new data set ${\Xi}_{n+1}$ is
streamed in and a new certificate
${J}_{n+1}^{\epsilon_1}(\vect{x}^{(r)})$ is yet to be obtained. It
can be seen that after a few samples streamed, both the obtained
certificate becomes close (\reviseone{within $10\%$}) to the estimated
true optimal value $J^{\star}$. In Figure~\ref{fig:CertnoICA} however,
as the data streams over $50$ seconds, the computing cost for updating
certificates becomes significant for the algorithm without the
\ICAaco. After $\supscr{100}{th}$ data point \reviseone{has been
  assimilated}, the certificate ${J}_n^{\epsilon_1}(\vect{x}^{(r)})$
stops updating for all $n \geq 100$. And, further, after all the data
points streamed \reviseone{(in $200$ seconds)}, the algorithm took about $70$ seconds to terminate
the algorithm with certificate
${J}_{200}^{\epsilon_1}(\vect{x}^{(r)})$. This is a clear disadvantage
compared to the algorithm with the \ICAaco, which terminates as soon
as all the data points were taken in.

\subsection{Study 2: \ODAA with Large Streaming Data Sets}
Here, we are to find an $\epsilon_2$-optimal, $\epsilon_1$-proper
decision $\vect{x} \in \real^{30}$ for Problem~\eqref{eq:P}. We
consider $N=500$ iid sample points $\{{\xi}_{k}\}_{k=1}^{N}$ streaming
randomly in between every $1$ to $3$ seconds with each data point
${\xi}_{k} \in \real^{10}$ a realization of $\prob$. We assume that
the unknown distribution is a multivariate Gaussian mixture
distribution with three centers where the components of the mean of
each center are uniformly chosen between $[-10,10]$, and the variance
matrix is $\vect{I}_m$ for each center. We assume the cost function
$\map{f}{\real^{30} \times \real^{10}}{\real}$ to be
$f(\vect{x},\xi):=\trans{\vect{x}}A\vect{x}+\trans{\vect{x}}B\xi+\trans{\xi}C\xi$
with random values for the positive semi-definite matrix $A \in
\real^{30 \times 30}$, $B \in \real^{30 \times 10}$ and negative
definite matrix $C \in \real^{10 \times 10}$. The radius of the
Euclidean ball for the \ICAaco is $w=5$.

Similarly to Figure~\ref{fig:CP}, Figure~\ref{fig:CP500} demonstrates
the incremental construction of the distribution $\tilde{\prob}^n$ and
the accumulated number of Problem~\eqref{eq:tildePJoverSimplex} solved
over time. Clearly, after certain amount of data \reviseone{have been
  assimilated}, the structure of the data set was inferred by the
\ICAaco and the number of Euclidean balls used to cover the data set
is about $20$. Also, after the $\supscr{100}{th}$ time period (from
$100$ to $200$ seconds in this case), the algorithm can \reviseone{validate} new
certificate
without solving any
Problem~\eqref{eq:tildePJoverSimplex}. This feature dramatically
improves the performance of the \ODAA and makes the algorithm flexible
for online settings.

  Similarly to Figure~\ref{fig:CertICA}, Figure~\ref{fig:CertICA500}
  shows the evolution of the certificate sequence
  $\{{J}_n^{\epsilon_1}(\vect{x}^{(r)}) \}_{n=1,r=1}^{N,\infty}$ for
  the decision sequence
  $\{\vect{x}^{(r)}\}_{r=1}^{\infty}$. \reviseone{In the same way} as
  \reviseone{in the }last case study, the obtained certificate becomes
  close to the estimated true optimal value $J^{\star}$
  \reviseone{(within $10\%$)} after about $25$ seconds with
  \reviseone{the assimilation of} $10$ data sets. Also, as
  \reviseone{more} data sets \reviseone{are assimilated}, the update
  of the certificate
  ${J}_n^{\epsilon_1}(\vect{x}^{\epsilon_2}_{n})$ remains
  \reviseone{fast} and the algorithm terminates \reviseone{within a second} after the last data set was streamed in.
\section{Conclusions}\label{sec:Conclude}
In this paper, we have proposed the \ODAAfull (the \ODAA) to solve the
problem in the form of~\eqref{eq:P}, where the realizations of the
unknown distribution (i.e., the streaming data) are collected over
time in order for the real-time data-driven decision of~\eqref{eq:P}
to have guaranteed out-of-sample performance.
The data-driven
decision with the certificate that guarantees out-of-sample
performance are available any time during the execution of the
algorithm, and the optimal data-driven decision are approached with a
(sub)linear convergence rate. The algorithm terminates after
collecting a sufficient amount of data to make good decision. To
facilitate the decision making, an enhanced version of the proposed
algorithm is further constructed, by using an \ICA (the \ICAaco) to
estimate new ambiguity sets over time. We provided sample problems and
showed the actual performance of the proposed \ODAA with the \ICAaco
over time. Future work will generalize the results for weaker
assumptions of the problem and potentially extend the algorithm to
scenarios that include system dynamics.

\appendix
There are mainly two types of Numerical
methods that serve as
the main ingredients of our
\ODAA. One type is given by
Frank-Wolfe Algorithm (FWA) variants and another is the Subgradient
Algorithm.  In this Section, we describe FWA and the Away-step
Frank-Wolfe Algorithm (AFWA) for the sake of completeness.  We combine
AFWA with another variant, the Simplicial Algorithm, in
Section~\ref{sec:CG}. For the Subgradient Algorithm, please refer
to~\cite{SMR:99,TL-MP-AS:03,YN:13}.
\subsection*{Frank-Wolfe Algorithm over a unit simplex}
To solve convex programs over a unit simplex, we introduce the FWA and
AFWA following~\cite{SLJ-MJ:15,CH:74}. Let us denote the
$m$-dimensional unit simplex by $\Delta_m:=\setdef{\lambda \in
  \real^m}{\trans{\vectorones{m}} \lambda=1, \; \lambda \geq 0}$. Let
${\Lambda}_{m}$ be the set of all extreme points for the simplex
$\Delta_m$. Consider the maximization of a concave function
$f(\vect{x})$ subject to $\vect{x} \in \Delta_m$; we refer to this
problem by $(\star)$ and denote by $\vect{x}^{\star}$ an optimizer of
$(\star)$. We call ${\vect{x}}^{\epsilon}$ an $\epsilon$-optimal
solution of $(\star)$, if ${\vect{x}}^{\epsilon} \in \Delta_m$ and
$f(\vect{x}^{\star}) -f({\vect{x}}^{\epsilon})\leq \epsilon$.  The
classical FWA solves problem $(\star)$ to an ${\vect{x}}^{\epsilon}$
via the iterative process as follows.  Let $\vect{x}^{(0)} \in
\Delta_m$ denote a random initial point for FWA.  For each iteration
$k$ with an $\vect{x}^{(k)} \in \Delta_m$, the concavity of $f$
enables $f(\vect{x}^{\star}) \leq f(\vect{x}^{(k)}) + \trans{\nabla
  f(\vect{x}^{(k)})}(\vect{x}^{\star}-\vect{x}^{(k)})$, which implies
$f(\vect{x}^{\star}) \leq f(\vect{x}^{(k)}) + \max_{\vect{x} \in
  \Delta_m} \trans{\nabla
  f(\vect{x}^{(k)})}(\vect{x}-\vect{x}^{(k)})$. Using this property,
we define a FW search point $\vect{s}^{(k)}$ by an extreme point such
that $\vect{s}^{(k)} \in \argmax_{\vect{x} \in \Delta_m} \trans{\nabla
  f(\vect{x}^{(k)})}(\vect{x}-\vect{x}^{(k)})$. With this search point
we define the FW direction at $\vect{x}^{(k)}$ by
$\subscr{d}{FW}^{(k)}:=\vect{s}^{(k)}-\vect{x}^{(k)}$.  The classical
FWA then iteratively finds a FW direction and solves a line search
problem over this direction until an ${\epsilon}$-optimal solution
${\vect{x}}^{\epsilon}:=\vect{x}^{(k)}$ is found, certified by
$\eta^{(k)}:=\trans{\nabla f(\vect{x}^{(k)})}\subscr{d}{FW}^{(k)}
\leq \epsilon$.

It is known that the classical FWA has linear convergence rate if the
cost function $f$ is $\mu$-strongly concave and the optimum is achieved
in the relative interior of the feasible set $\Delta_m$. If the
optimal solution lies on the boundary of $\Delta_m$, then this
algorithm only has a sublinear convergence rate, due to a zig-zagging
phenomenon~\cite{SLJ-MJ:15}.  AFWA is an extension of the FWA that
guarantees the linear convergence rate of the problem $(\star)$ under
some conditions related to the local strong concavity. The main
difference between AFWA and the classical FWA is that the latter
solves the line-search problem after obtaining a ascent direction by
considering all extreme points, while the AFWA chooses a ascend
direction that prevents zig-zagging. We summarize the convergence
properties of the AFWA here. For complete descriptions of the AFWA, we
refer the reader to~\cite{SLJ-MJ:15,DL-SM:18-extended}.
  The detailed FWA and AFWA are shown in the following Algorithm tables.
  {\begin{algorithm}[h]
  \caption{Classical FWA for $(\star)$: FW$(\vect{x}^{(0)},\Delta_m,\epsilon)$.}
  \label{alg:ClassicFW}
  \begin{algorithmic}[1]
    \Ensure $\epsilon$-optimal $\vect{x}^{\epsilon}$;
    \State Set
     $k \leftarrow 0$, $\eta^{(k)}
    \leftarrow +\infty$;
    \Repeat
    \State Pick $\vect{s}^{(k)} \in \argmax\limits_{\vect{x} \in {\Lambda}_{m}}
    \trans{\nabla f(\vect{x}^{(k)})}(\vect{x}-\vect{x}^{(k)})$;
    \State $\subscr{d}{FW}^{(k)} \leftarrow \vect{s}^{(k)}-\vect{x}^{(k)}$;
    \State $\eta^{(k)} \leftarrow \trans{\nabla f(\vect{x}^{(k)})}\subscr{d}{FW}^{(k)}$;
      \State  Pick $\gamma^{(k)} \in \argmax\limits_{\gamma \in [0,1]} f(\vect{x}^{(k)} +
    \gamma \subscr{d}{FW}^{(k)})$;
    \State $\vect{x}^{(k+1)} \leftarrow \vect{x}^{(k)}+\gamma^{(k)} \subscr{d}{FW}^{(k)}$;
    \State $k \leftarrow k+1$;
    \Until $\eta^{(k)} \leq \epsilon$;
    \State { Return $\vect{x}^{(k)}$.}
  \end{algorithmic}
  \end{algorithm}}

\begin{algorithm}[htb!]
\caption{AFWA for $(\star)$: $(\vect{x}^{\epsilon}, {\obj}^{\epsilon}) \leftarrow$ AFW$(f(\vect{x}), \Delta_m,\epsilon)$.}
\label{alg:AFW}
\begin{algorithmic}[1]
\Ensure
$\epsilon$-optimal $\vect{x}^{\epsilon}$ with objective ${\obj}^{\epsilon}$;
\State Set $k \leftarrow 0$, $\eta^{(k)} \leftarrow +\infty$;
\State Pick $\vect{x}^{(k)} \in {\Lambda}_{m}$, $\subscr{I}{Act}^{(k)}:=\{\vect{x}^{(k)} \}$, $p=|\subscr{I}{Act}^{(k)}|$;
\State Let
$ \alpha^{(k)}_v=\begin{cases}
  1/p, \textrm{ if } \vect{v} \in \subscr{I}{Act}^{(k)}, \\
  0 , \textrm{ if } \vect{v} \in {\Lambda}_{m} - \subscr{I}{Act}^{(k)},
\end{cases}$
\Repeat
\State Pick $\vect{s}^{(k)} \in \argmax\limits_{\vect{x} \in {\Lambda}_{m}} \nabla f(\vect{x}^{(k)})^T(\vect{x}-\vect{x}^{(k)})$;
\State $\subscr{d}{FW}^{(k)} \leftarrow \vect{s}^{(k)}-\vect{x}^{(k)}$;
\State Pick $\vect{v}^{(k)} \in \argmin\limits_{\vect{x} \in \subscr{I}{Act}^{(k)}} \nabla f(\vect{x}^{(k)})^T(\vect{x}-\vect{x}^{(k)})$;
\State $\subscr{d}{A}^{(k)} \leftarrow \vect{x}^{(k)}-\vect{v}^{(k)}$; \Comment{Away-step direction}
\If{$\left\langle \nabla f(\vect{x}^{(k)}), \subscr{d}{FW}^{(k)} \right\rangle  \geq \left\langle \nabla f(\vect{x}^{(k)}), \subscr{d}{A}^{(k)} \right\rangle$, }
\State $d^{(k)} \leftarrow \subscr{d}{FW}^{(k)}$;
\State flag $\leftarrow$ True, $\subscr{\gamma}{max} \leftarrow 1$;
\Else \Comment{AFW direction has larger potential ascent}
\State $d^{(k)} \leftarrow \subscr{d}{A}^{(k)}$;
\State $\subscr{\gamma}{max} \leftarrow \alpha^{(k)}_{\vect{v}^{(k)}}/(1-\alpha^{(k)}_{\vect{v}^{(k)}})$;
\EndIf
\State Pick $\gamma^{(k)} \in \argmax\limits_{\gamma \in [0,\subscr{\gamma}{max}]} f(\vect{x}^{(k)} + \gamma d^{(k)})$;
\If{flag is True, }
 \If {$\gamma^{(k)}=1$, }  \Comment{Hit extreme point}
  \State $\subscr{I}{Act}^{(k+1)} \leftarrow \{\vect{s}^{(k)} \}$;
 \Else
  \State $\subscr{I}{Act}^{(k+1)} \leftarrow \subscr{I}{Act}^{(k)} \cup \{\vect{s}^{(k)} \}$;
 \EndIf
\State $\alpha^{(k+1)}_{\vect{s}^{(k)}} \leftarrow (1-\gamma^{(k)})\alpha^{(k)}_{\vect{s}^{(k)}}+\gamma^{(k)}$;
\State $\alpha^{(k+1)}_{\vect{v}} \leftarrow (1-\gamma^{(k)})\alpha^{(k)}_{\vect{v}}$, $\forall \vect{v} \in \subscr{I}{Act}^{(k)} - \{\vect{s}^{(k)} \}$;
\Else
 \If{$\gamma^{(k)}=\subscr{\gamma}{max}$, } \Comment{Hit $\Delta_m$ boundary}
   \State $\subscr{I}{Act}^{(k+1)}\leftarrow\subscr{I}{Act}^{(k)}-\{\vect{v}^{(k)} \}$;
 \Else
   \State $\subscr{I}{Act}^{(k+1)}\leftarrow\subscr{I}{Act}^{(k)}$;
 \EndIf
\State $\alpha^{(k+1)}_{\vect{v}^{(k)}}=(1+\gamma^{(k)})\alpha^{(k)}_{\vect{v}^{(k)}}-\gamma^{(k)}$;
\State $\alpha^{(k+1)}_{\vect{v}}=(1+\gamma^{(k)})\alpha^{(k)}_{\vect{v}}$, $\forall \vect{v} \in \subscr{I}{Act}^{(k)} - \{\vect{v}^{(k)} \}$;
\EndIf
\State $\vect{x}^{(k+1)}\leftarrow\vect{x}^{(k)}+\gamma^{(k)} d^{(k)}$;
\State $k \leftarrow k+1$;
  \Until $\eta^{(k)} \leq \epsilon$;
\State Return $\vect{x}^{\epsilon} \leftarrow \vect{x}^{(k)}$ and ${\obj}^{\epsilon} \leftarrow f(\vect{x}^{\epsilon})$.
\end{algorithmic}
\end{algorithm}

\begin{theorem}[Linear convergence of AFWA{~\cite[Theorem 8]{SLJ-MJ:15}}]
  Suppose the function $f$ has a curvature constant $C_f$ and a geometric
		strong concavity constant $\mu_f$ on $\Delta_m$, as defined in \revisetwo{footnote~\ref{footnote:concavity}}
  Let us define the decay rate $\kappa:=1-{\mu_f}/({4C_f}) \in (0,1) \subset \real$. Then the suboptimality bound at the
  iteration point $\vect{x}^{(k)}$ of the AFWA decreases
  geometrically  as
  $f(\vect{x}^{\star}) -f(\vect{x}^{(k+1)}) \le
  \kappa(f(\vect{x}^{\star}) -f(\vect{x}^{(k)}) )$.
\label{thm:AFWconv}
\end{theorem}
\section*{Proofs}
\subsubsection*{Lemma~\ref{lemma:certgener}}
\begin{IEEEproof}
  Following~\cite{AC-JC:17-allerton,PME-DK:17} and from
  Theorem~\ref{thm:MeasConc}, we prove that ${J}_{n}(\vect{x})$
  is a valid certificate for~\eqref{eq:perfgua}. Knowing
  that~\eqref{eq:epsirad} is obtained by letting the right-hand side
  of~\eqref{eq:radB} to be equal to a given $\beta_{n}$, for each~$n$
  we substitute~\eqref{eq:epsirad} into the right-hand side
  of~\eqref{eq:radB}, yielding $\mathbf{P}^n\{
  d_W(\prob,\hat{\prob}^n) \geq \epsilon(\beta_{n}) \} \leq \beta_{n}
  $ for each $n$. This means that a
  data set ${\Xi}_{n}$ we can construct an empirical probability measure
  $\hat{\prob}^n$
  such that $d_W(\prob,\hat{\prob}^n) \le \epsilon(\beta_{n})$ with
  probability at least $1-\beta_{n}$. Namely, $\mathbf{P}^n\{ \prob
  \in \mathbb{B}_{\epsilon(\beta_{n})}(\hat{\prob}^n) \} \geq 1-
  \beta_{n} $. Thus, for all $\vect{x} \in \real^d$, we
  have $\mathbf{P}^n\{ \prob \in
  \mathbb{B}_{\epsilon(\beta_{n})}(\hat{\prob}^n) \} = \mathbf{P}^n\{
  \mathbb{E}_{\prob} [f(\vect{x},\xi)] \leq
  \sup\limits_{\mathbb{Q} \in \mathcal{P}_n} {\mathbb{E}_{\mathbb{Q}}
    [f(\vect{x},\xi)]} \}=\mathbf{P}^n\{ \mathbb{E}_{\prob}
  [f(\vect{x},\xi)] \leq {J}_{n}(\vect{x}) \} \geq 1-
  \beta_{n}$.
\end{IEEEproof}
\subsubsection*{Lemma~\ref{lemma:equiP}}
\reviseone{
  In the following, we for matrices $A_1 \in \real^{m \times d}$ and
  $A_2 \in \real^{p \times q}$, we let $A_1 \oplus A_2$ denote their
  direct sum. The shorthand notation $\oplus_{i=1}^m A_i$ represents
  $A_1 \oplus \cdots \oplus A_m$. }
\begin{IEEEproof}
  To prove 1, for any feasible solution
  $(\tilde{\vect{u}},\tilde{\vect{v}})$ of~\eqref{eq:JoverSimplex}, we
  compute $\frac{1}{n} \sum_{k=1}^{n}
  \Norm{\tilde{\vect{y}}_k}=\frac{1}{n} \sum_{k=1}^{n}
  \Norm{\tilde{\vect{u}}_k-\tilde{\vect{v}}_k} \le \frac{1}{n}
  \sum_{k=1}^{n} \Norm{\tilde{\vect{u}}_k} + \frac{1}{n}
  \sum_{k=1}^{n} \Norm{\tilde{\vect{v}}_k}= \frac{1}{n}
  \trans{\vectorones{mn}}\tilde{\vect{u}}+\frac{1}{n}
  \trans{\vectorones{mn}}\tilde{\vect{v}} =\frac{1}{n}
  \trans{\vectorones{2mn}}(\tilde{\vect{u}},\tilde{\vect{v}}) =
  \epsilon(\beta_{n})$. Therefore
  $(\tilde{\vect{y}}_1,\ldots,\tilde{\vect{y}}_n)$ is feasible
  for~\eqref{eq:convJn}.

  For 2, we exploit that any feasible solution $\tilde{\vect{y}}$
  of~\eqref{eq:convJn} is a linear combination of the extreme points
  of the constraint set in~\eqref{eq:convJn}. Let us denote the matrix
  $A_n:=[\oplus_{i=1}^{n}I_m, -\oplus_{i=1}^{n}I_m] \in \real^{mn
    \times
    2mn}$.
  By construction of Problem~\eqref{eq:JoverSimplex}, we see that each
  column vector of the matrix $n\epsilon(\beta_{n})A_n$ is a
  concatenated vector of an extreme point of
  Problem~\eqref{eq:convJn}, and that all the extreme points
  of~\eqref{eq:convJn} are included. Then, any feasible solution
  of~\eqref{eq:convJn} can be written as
  $\tilde{\vect{y}}=n\epsilon(\beta_{n})A_n(\hat{\vect{u}},\hat{\vect{v}})$
  where $(\hat{\vect{u}},\hat{\vect{v}})$ is a vector of the convex
  combination coefficients of the extreme points of the constraint set
  in~\eqref{eq:convJn}. Clearly, we have
  $(\hat{\vect{u}},\hat{\vect{v}}) \in \Delta_{2mn}$, i.e.,
  $n\epsilon(\beta_{n})(\hat{\vect{u}},\hat{\vect{v}})$ is in the
  feasible set of the Problem~\eqref{eq:JoverSimplex}. Then, by
  construction
  $(\tilde{\vect{u}},\tilde{\vect{v}}):=n\epsilon(\beta_{n})(\hat{\vect{u}},\hat{\vect{v}})$
  is feasible for~\eqref{eq:JoverSimplex}.

  For 3, since~\eqref{eq:convJn} and~\eqref{eq:JoverSimplex} are the
  same in the sense of (1) and (2), then if
  $(\tilde{\vect{u}}^{\star},\tilde{\vect{v}}^{\star})$ is an
  optimizer of~\eqref{eq:JoverSimplex}, by letting
  $\tilde{\vect{y}}_k^{\star}:=\tilde{\vect{u}}_k^{\star}-\tilde{\vect{v}}_k^{\star}$
  for each $k \in \until{n}$ we know the objective values of the two
  problems coincide. We claim that the optimum of~\eqref{eq:convJn} is
  achieved via the optimizer $\tilde{\vect{y}}^{\star}$. If not, then
  there exists $\hat{\vect{y}}^{\star} \ne \tilde{\vect{y}}^{\star}$
  such that the optimum is achieved with higher value. Then, from the
  construction in (2) we can find a feasible solution
  $(\hat{\vect{u}},\hat{\vect{v}})$ of~\eqref{eq:JoverSimplex} that
  results in a higher objective value. This contradicts the assumption
  that $(\tilde{\vect{u}}^{\star},\tilde{\vect{v}}^{\star})$ is an
  optimizer
  of~\eqref{eq:JoverSimplex}.
\end{IEEEproof}
\subsubsection*{Theorem~\ref{thm:ConvergeJhat}}
\begin{IEEEproof}
  Given tolerance $\epsilon_1$, decision $\vect{x}$
  and \reviseone{any} data set
  ${\Xi}_{n}$ \reviseone{with $n \in \until{N}$}, let $\map{H_n}{\real^{mn}}{\real}$,
  $H_n:=\frac{1}{n}\sum_{k=1}^{n} h_k$ denote the objective function
  of~\eqref{eq:JoverSimplex} and let $\mathcal{S}_n$ denote the family
  of subsets of $\Lambda_{2mn}$. In the procedure of
  \CGAaco, let us consider a sequence of generated
  candidate vertex sets: $I_{n}^{(l)}\subset I_{n}^{(l+1)}$, $l=0,1,2,
  \ldots$ with $I_{n}^{(l)} \in \mathcal{S}_n$.  We show the
  convergence of \CGAaco for any data set
  ${\Xi}_{n}$, by two steps.

  Step 1) The sequence $\{I_{n}^{(l)}\}_l$ is finite and
  the number of iterations is at most $2mn$. For each $l$ and
  candidate optimizer $\vect{y}^{(l-1)}$, we generate a nonempty set
  of search points $\Omega^{(l)}$ with suboptimality gap $\eta^{(l)}$
  via~\eqref{eq:LP}. If $\eta^{(l)} \leq \epsilon_1$, then we
  solved~\eqref{eq:JoverSimplex} to $\epsilon_1$-optimality and $l$ is
  therefore finite, otherwise we update $I_{n}^{(l)}:=I_{n}^{(l-1)}
  \cup \Omega^{(l)}$. Given that the maximal cardinality of each
  $I_{n}^{(l)} \in \mathcal{S}_n$ is bounded by $2mn$, then
  it is sufficient to show $\Omega^{(l)} \cap
  I_{n}^{(l-1)}=\O$.  Because $\vect{y}^{(l-1)}$ is an
  $\epsilon_1$-optimal of~\eqref{eq:CP} under
  ${\conv}(I_{n}^{(l-1)})$, then for any ${\vect{y}} \in
  {\conv}(I_{n}^{(l-1)})$, it holds that $\frac{1}{n} \sum_{k=1}^{n}
  \langle \nabla h_k(\vect{y}^{(l-1)}_k),
  \vect{y}_k-\vect{y}^{(l-1)}_k \rangle \leq \epsilon_1$. Since any
  element in $\Omega^{(l)}$ is such that $\eta^{(l)} > \epsilon_1$,
  then for any ${\vect{y}} \in {\conv}(I_{n}^{(l-1)})$, we have
  ${\vect{y}}\notin \Omega^{(l)}$, which concludes $\Omega^{(l)} \cap
  I_{n}^{(l-1)}=\O$. Further, the cardinality of $\Omega^{(l)}$ is at
  least one for every iteration $l$, then after at most $2mn$ steps
  the cardinality of $I_{n}^{(l)}$ becomes $2mn$, which implies the
  $\epsilon_1$-optimality of~\eqref{eq:JoverSimplex} by the
  $\epsilon_1$-optimality of~\eqref{eq:CP}.

  Step 2)
  The computational bound of \CGAaco is quantified.  To see this,
  consider the problems~$\{\eqref{eq:LP}\}_l$
  and~$\{\eqref{eq:CP}\}_l$. By Assumption~\ref{assump:chgrad} on the
  cheap access of the gradients, the computation of~\eqref{eq:LP} is
  negligible. Thus, the computational bound is given by the sum of the
  steps to solve the $\{\eqref{eq:CP}\}_l$, where the number of
  iterations $l$ is $2mn$ in the worst
  case.

  For each~\eqref{eq:CP} solved by AFWA, index the AFWA iterations by
  $i=0,1,2,\ldots$, let ${\obj}_{i}^{(l)}$ be the objective value at
  each iteration, and assume the optimal objective value is
  ${\obj}_{\star}^{(l)}$. As in Theorem~\ref{thm:AFWconv},
   let
  $\kappa_{\reviseone{n},l} \in (0,1) \subset \real$ be the decay
  parameter related to local strong concavity of $H_n$ over
  $\conv(I_{n}^{(l)})$. Then using the linear convergence rate of the
  AFWA, each~\eqref{eq:CP} achieves the following computational bound
\begin{equation*}
  {\obj}_{\star}^{(l)} - {\obj}_{i}^{(l)} \le \kappa_{\reviseone{n},l}^{i}( {\obj}_{\star}^{(l)} -{\obj}_{0}^{(l)}),
\end{equation*}
where the initial condition ${\obj}_{0}^{(l)}$ results from an
$\epsilon_1$-optimal optimizer of CP at iteration $l-1$, i.e., we can
equivalently denote ${\obj}_{0}^{(l)}$ by
${J}_n^{(l-1)}(\vect{x})$, for all $l \in \until{2mn}$.

Let us consider sequence~$\{\eqref{eq:CP}\}_{l}$ with feasible sets
$\{\conv(I_{n}^{(l)})\}_{l}$. Then we have
\begin{center} {\centering $\conv(I_{n}^{(0)}) \subset
    \conv(I_{n}^{(1)}) \subset \cdots \subset \conv(I_{n}^{(2mn)})$.
  } \end{center} This results into monotonically decaying parameters
and ($\epsilon_1$-)optimal objective values, as given in the following
\begin{center} {
    $0 <\kappa_{\reviseone{n},1} \leq \kappa_{\reviseone{n},2} \leq \cdots \leq \kappa_{\reviseone{n},2mn}<1$, \\
    ${J}_n^{(0)}(\vect{x}) \leq{J}_n^{(1)}(\vect{x}) \leq \cdots \leq {J}_n^{(2mn)}(\vect{x})$, \\
    ${\obj}^{(0)}_{\star} \leq {\obj}_{\star}^{(1)} \leq \cdots \leq
    {\obj}_{\star}^{(2mn)}$.  }
\end{center}
Using the previous notation, we can identify
${J}_n^{(2mn)}(\vect{x}) \equiv
{J}_n^{\epsilon_1}(\vect{x})$, ${\obj}^{(0)}_{\star} \equiv
{J}_n^{(0)}(\vect{x})$, and
${\obj}_{\star}^{(2mn)}\equiv{J}_n(\vect{x})$. Let us denote
$\kappa:=
\max_{\reviseone{n},l}\{\kappa_{\reviseone{n},l}\}$.  Then, by solving each~$\eqref{eq:CP}$
to $\epsilon_1$-optimality, it leads to the accumulated computational
steps $\phi(n):=\sum_{l} i_l$, where each $i_l$ is the computation
step for $\epsilon_1$-optimal~\eqref{eq:CP} that satisfies the
following inequality
 \begin{equation*}
   \kappa^{i_l}( {J}_n(\vect{x}) -{J}_n^{(0)}(\vect{x})) \leq \epsilon_1, \quad l \in \until{2mn}.
\end{equation*}
Finally, in the worst-case scenario, the computational bound of
the \CGAaco is
 \begin{equation*}
   \phi(n) \leq (2mn) {\log}_{\kappa}(\frac{\epsilon_1}{{J}_n(\vect{x}) -{J}_n^{(0)}(\vect{x})}).
\end{equation*}
Next, we show the convergence of the \CGAaco under online
data sets $\{{\Xi}_{n}\}_{n=1}^{N}$. Similarly to the proof for the
computational bound for a given $n$, we can compute the worst-case
bound under $\{{\Xi}_{n}\}_{n=1}^{N}$, by summing over the steps
required to solve the~$\{\eqref{eq:CP}\}_{n,l}$. This leads to the
stated bound $\bar{\phi}(n)$, where the empirical cost
$\supscr{J}{sae}_{N}(\vect{x})$ serves as the cost of initial
condition $\vect{y}^{(0)}:=\vectorzeros{2mN}$.  In this way, when
 the data-streaming rate is slower \reviseone{or equal} than $(\bar{\phi}(1))^{-1}$,
we
claim that \CGAaco can always find the certificate
for each data set ${\Xi}_{n}$. This is because in each time period $n$, we only have $2mn$ extreme points, and $2m(n-1)$ has been
explored due to the adaptation of the candidate vertex set
$I_{n}^{(0)}$.
\end{IEEEproof}
\subsubsection*{Lemma~\ref{lemma:cvxJn}}
\begin{IEEEproof}
  For any $\vect{x}$, $\vect{y} \in \real^d$ and $t \in [0,1] \subset
  \real$, we have  $\vect{z}=t \vect{x}+(1-t)\vect{y} \in
  \real^d$ and an optimizer of~\eqref{eq:cert},
  $\mathbb{Q}_n^{\star}(\vect{z})$, such that {\begin{equation*}
\begin{split}
  {J}_{n}(\vect{z})\leq &
  {\mathbb{E}_{\mathbb{Q}_n^{\star}(\vect{z})}
    [t f(\vect{x},\xi) + (1-t) f(\vect{y},\xi)] } \\
  =& t{\mathbb{E}_{\mathbb{Q}_n^{\star}(\vect{z})}
    [f(\vect{x},\xi)] +
    (1-t) \mathbb{E}_{\mathbb{Q}_n^{\star}(\vect{z})}[f(\vect{y},\xi)] } \\
  \le & t {J}_{n}(\vect{x})+(1-t) {J}_{n}(\vect{y}).
\end{split}
\end{equation*}}
\end{IEEEproof}
\subsubsection*{Lemma~\ref{lemma:access_subg}}
\begin{IEEEproof}
  Let us consider the function
  ${\mathbb{E}_{\mathbb{Q}_n^{\epsilon_1}(\vect{x}^{(r)})}
    [f(\vect{x},\xi)] }$. Using Assumption~\ref{assump:cvx} on
  convexity of $f$ in $\vect{x}$, we have for any $\vect{z} \in
  {\dom}\: {J}_{n}$ the following relation
\begin{equation*}
{\mathbb{E}_{\mathbb{Q}_n^{\epsilon_1}(\vect{x}^{(r)})} [f(\vect{z},\xi)] }
 \geq {J}_n^{\epsilon_1}(\vect{x}^{(r)}) + \trans{g_{n}^{r}(\vect{x}^{(r)}) }(\vect{z}-\vect{x}^{(r)}).
\end{equation*}
Knowing that ${J}_n^{\epsilon_1}(\vect{x}^{(r)}) \geq  {J}_n(\vect{x}^{(r)}) - \epsilon_1$ and ${J}_n(\vect{z}) \geq  {\mathbb{E}_{\mathbb{Q}_n^{\epsilon_1}(\vect{x}^{(r)})} [f(\vect{z},\xi)] }$,
this concludes the first part of the proof.

To show the second part, similarly, we also have for any $\tilde{\vect{x}}$, $\vect{z} \in {\dom}\: {J}_{n}$ the following relation
\begin{equation*}
{\mathbb{E}_{\mathbb{Q}_n^{\epsilon_1}(\vect{x}^{(r)})} [f(\vect{z},\xi)] }
 \geq {\mathbb{E}_{\mathbb{Q}_n^{\epsilon_1}(\vect{x}^{(r)})} [f(\tilde{\vect{x}},\xi)] }+ \trans{g_{n}^{r}(\tilde{\vect{x}}) }(\vect{z}-\tilde{\vect{x}}).
\end{equation*}
Using Point Search Algorithm,
we achieve an $\eta>0$ such that ${\mathbb{E}_{\mathbb{Q}_n^{\epsilon_1}(\vect{x}^{(r)})} [f(\tilde{\vect{x}},\xi)] } \geq {J}_n(\tilde{\vect{x}}) - \eta$. Finally, by similar statement as in the first part, we claim $g_{n}^{r}(\tilde{\vect{x}}) \in \partial_{\epsilon} {J}_{n}(\tilde{\vect{x}})$.
\end{IEEEproof}
\subsubsection*{Lemma~\ref{lemma:convergeX}}
\begin{IEEEproof}
In the $\supscr{n}{th}$ time period, let us consider subgradient iterates $i$ for all $r_n \leq i \leq r$
\begin{equation*}
\begin{split}
  & \Norm{\vect{x}^{(i+1)} - \vect{x}_{n}^{\star}}^2 =  \Norm{ \vect{x}^{(i)}-\vect{x}_{n}^{\star} -\alpha^{(i)}
\frac{g_{n}^{i}(\vect{x}^{(i)})}{\max\{
    \Norm{g_{n}^{i}(\vect{x}^{(i)})} \; ,
    \; 1 \} }  }^2 \\
  &=  \Norm{ \vect{x}^{(i)} -\vect{x}_{n}^{\star} }^2 +(\alpha^{(i)})^2 \min \{ \Norm{ g_{n}^{i}(\vect{x}^{(i)}) }^2 , 1 \} \\
  & \qquad -2 \alpha^{(i)}
\frac{
    \trans{g_{n}^{i}(\vect{x}^{(i)})}(\vect{x}^{(i)}
    -\vect{x}_{n}^{\star})}{\max\{
    \Norm{g_{n}^{i}(\vect{x}^{(i)})} \; ,
    \; 1 \} } .
\end{split}
\end{equation*}
From Lemma~\ref{lemma:access_subg}, we know that
${J}_{n}(\vect{x}_{n}^{\star})\geq {J}_{n}(\vect{x}^{(i)}) + \trans{g_{n}^{i}(\vect{x}^{(i)})}(\vect{x}_{n}^{\star}-\vect{x}^{(i)}) - \subscr{\epsilon}{SA}$
for all $\vect{x}^{(i)}$. Then, we have
\begin{equation*}
\begin{split}
  &\Norm{ \vect{x}^{(i+1)} -\vect{x}_{n}^{\star} }^2 \leq (\alpha^{(i)})^2 \min \{ \Norm{g_{n}^{i}(\vect{x}^{(i)})  }^2 , 1 \}  \\
  &+\Norm{ \vect{x}^{(i)} -\vect{x}_{n}^{\star} }^2 +
  \frac{2 \alpha^{(i)} ({J}_{n}(\vect{x}_{n}^{\star})
    -{J}_{n}(\vect{x}^{(i)}) +\subscr{\epsilon}{SA} )}{\max\{
    \Norm{ g_{n}^{i}(\vect{x}^{(i)}) } \; ,
    \; 1 \} } .
\end{split}
\end{equation*}
Combining the inequalities over iterations from $r_n$ to $r$ \reviseone{gives}
\begin{equation*}
\begin{split}
& 0 \leq   \Norm{ \vect{x}^{(r_{n})} -\vect{x}_{n}^{\star} }^2 + \sum_{i=r_{n}}^r {(\alpha^{(i)})^2 \min \{ \Norm{ g_{n}^{r}(\vect{x}^{(i)})   }^2 \; , \; 1 \}} \\
& +\sum_{i=r_{n}}^r { \frac{2 \alpha^{(i)} ({J}_{n}(\vect{x}_{n}^{\star}) -{J}_{n}(\vect{x}^{(i)}) +\subscr{\epsilon}{SA} )}{\max\{ \Norm{ g_{n}^{r}(\vect{x}^{(i)})   } \; , \; 1 \} } } \\
& \leq  \Norm{ \vect{x}^{(r_{n})} -\vect{x}_{n}^{\star} }^2  +2 \subscr{\epsilon}{SA} \sum_{i=r_{n}}^r {\alpha^{(i)} }+ \sum_{i=r_{n}}^r {(\alpha^{(i)})^2} \\
& + \sum_{i=r_{n}}^r { \frac{2 \alpha^{(i)} ({J}_{n}(\vect{x}_{n}^{\star}) -{J}_{n}(\vect{x}^{(i)}))}{\max\{ \Norm{g_{n}^{r}(\vect{x}^{(i)})   } \; , \; 1 \} } }.
\end{split}
\end{equation*}
Then, using the fact that
\begin{equation*}
\begin{split}
&\sum_{i=r_{n}}^r { \frac{2 \alpha^{(i)} ({J}_{n}(\vect{x}_{n}^{\star}) -{J}_{n}(\vect{x}^{(i)}))}{\max\{ \Norm{ g_{n}^{r}(\vect{x}^{(i)})   } \; , \; 1 \} }}  \leq \\
& \sum_{i=r_{n}}^r { \frac{-2 \alpha^{(i)} \min\limits_{k\in \untilinterval{r_{n}}{r}} \{ {J}_{n}(\vect{x}^{(k)}) -{J}_{n}(\vect{x}_{n}^{\star}) \}}{\max\{ \Norm{ g_{n}^{r}(\vect{x}^{(i)})   } \; , \; 1 \} } } \leq \\
& -2(\sum_{i=r_{n}}^r { \alpha^{(i)}}) \frac{\min\limits_{k\in \untilinterval{r_{n}}{r}} \{{J}_{n}(\vect{x}^{(k)})\} -{J}_{n}(\vect{x}_{n}^{\star})}{\mu},
\end{split}
\end{equation*}
and the previous iteration,
we have
\begin{equation*}
\begin{split}
&\min\limits_{k\in \untilinterval{r_{n}}{r}} \{{J}_{n}(\vect{x}^{(k)})\} -{J}_{n}(\vect{x}_{n}^{\star}) \leq \\
&\frac{\mu  \Norm{ \vect{x}^{(r_{n})} -\vect{x}_{n}^{\star} }^2 + \mu \sum_{i=r_{n}}^r {(\alpha^{(i)})^2}}{2(\sum_{i=r_{n}}^r { \alpha^{(i)}})}+ \mu \subscr{\epsilon}{SA}.
\end{split}
\end{equation*}
\reviseone{Next, it remains to select a step-size rule $\{\alpha^{(i)}
  \}_{i=r_n}^r$ such that 1) the above right hand side term is upper
  bounded by $\epsilon_2$, and 2) the number of subgradient iterations
  $\bar{r}$ as described in the lemma is bounded. Note that the
  selection procedure is not unique, so we propose two step-size rules
  to obtain an explicit expression of $\bar{r}$.

  For any data set $\Xi_n$, let us select a sufficiently large value
  $M$ to be the diameter of the decision domain of interest, i.e.,
\begin{equation*}
  \Norm{ \vect{x}^{(r_{n})} -\vect{x}_{n}^{\star} } \leq M, \; \forall \, n \in \until{N}.
\end{equation*}
Then the step size rule and $\bar{r}$ has to satisfy the following
\begin{equation}
  \frac{ M^2 + \sum_{i=r_{n}}^{\bar{r}+r_{n}} {(\alpha^{(i)})^2}}{2(\sum_{i=r_{n}}^{\bar{r}+r_{n}} { \alpha^{(i)}})}+ \subscr{\epsilon}{SA} < \frac{\epsilon_2}{\mu}.
\label{eq:subrequire}
\end{equation}

We first consider a constant step-size rule, and select the step size as follows
\begin{equation*}
  \alpha^{(i)}:=\frac{M}{\sqrt{\bar{r}+1}}, \; \forall i \in \untilinterval{r_{n}}{\bar{r}+r_{n}}.
\end{equation*}
Then, to satisfy~\eqref{eq:subrequire}, we determine
\begin{equation*}
  \bar{r}:=M^{2}\left( \frac{\epsilon_2}{\mu}-\subscr{\epsilon}{SA}\right)^{-2}.
\end{equation*}

Alternatively, consider the divergent but square-summable step-size
rule, i.e.~$\sum_{i=r_{n}}^{\infty} { \alpha^{(i)}} = \infty$,
$\sum_{i=r_{n}}^{\infty} {(\alpha^{(i)})^2} < \infty$. For this class
of step-size rules, as $\bar{r}$ increases to $\infty$, we have the
left-hand side term of~\eqref{eq:subrequire} goes to $\mu
\subscr{\epsilon}{SA} < \epsilon_2$, then there exists a large enough
but finite number $\bar{r}$, such that~\eqref{eq:subrequire} holds. To
see this explicitly, we select the step-size rule to be the harmonic
sequence as follows
\begin{equation*}
  \alpha^{(i)}:=\frac{M}{i-r_{n}+1}, \; \forall i \in \untilinterval{r_{n}}{\bar{r}+r_{n}}.
\end{equation*}
Now we upper bound the numerator and lower bound the denominator
of~\eqref{eq:subrequire} using the following fact
\begin{equation*}
  \sum\limits_{i=1}^{\bar{r}+1} \frac{1}{i^{2}} \leq 2- \frac{1}{\bar{r}+1}, \quad \sum\limits_{i=1}^{\bar{r}+1} \frac{1}{i} \geq \ln(\bar{r}+1).
\end{equation*}
Then we determine $\bar{r}$ to be the following
\begin{equation*}
  \bar{r}= \min \setdef{r \in \mathbb{N}}{M(3-\frac{1}{r+1})\leq 2(\frac{\epsilon_2}{\mu} - \subscr{\epsilon}{SA})\ln(r+1)}.
\end{equation*}
This concludes the proof.
}
\end{IEEEproof}
\subsubsection*{Theorem~\ref{thm:epsSA}}
\begin{IEEEproof}
The computational bound to achieve an $\vect{x}_{n}^{\epsilon_2}$
strongly depends on the subgradient iterations
$\bar{r}:=r_{n+1}-r_{n}$ in Lemma~\ref{lemma:convergeX} and the number
of subgradient functions $\{ g_{n}^{r}\}_{r=r_n}^{r_{n+1}}$
constructed via the \CGAaco. To characterize this bound,
we quantify the computational steps for $\{
g_{n}^{r}\}_{r=r_n}^{r_{n+1}}$ next.

For each time period $n$, let us assume the \CGAaco has
explored the feasible set of~\eqref{eq:JoverSimplex} when obtaining
the initial certificate
${J}_n^{\epsilon_1}(\vect{x}^{(r_n)})$.
This procedure consumes a worst-case computational time $\phi(n)$, (or
$\bar{\phi}(1)$ if a data-streaming scenario), as stated in Theorem~\ref{thm:ConvergeJhat}. After
this initial step, every time the Subgradient Algorithm needs to
execute \CGAaco at some $r \leq r_{n+1}$,
\CGAaco will solve a unique~\eqref{eq:CP} and return
$\mathbb{Q}_n^{\epsilon_1}$ for an ${\epsilon}_1$-subgradient function
$g_{n}^{r}$ at $\vect{x}^{(r)}$.
Let ${\textrm{CP}}_r$ denote the unique~\eqref{eq:CP} solved at
$\vect{x}^{(r)}$. Then, to quantify the computational steps for $\{
g_{n}^{r}\}_{r=r_n}^{r_{n+1}}$, we compute the sum of the steps to
solve $\{ {\textrm{CP}}_r\}_r$.

Let us denote the number of steps solving ${\textrm{CP}}_r$ by $i_r$,
for all $r \in \{ r_n,\ldots, r_{n+1} \}$. Then, we aim to quantify
$i_{r+1}$ for $g_{n}^{r+1}$. To achieve this, let us assume a
subgradient function $g_{n}^{r}$ is computed at an iteration $r$.
Then we perform a subgradient iteration~\eqref{eq:subg} and obtain an
$\vect{x}^{(r+1)}$. By using a subgradient estimation strategy, we
obtain the optimality gap $\eta^{(1)}$ via Point Search Algorithm,
denoted by $\bar{\eta}_{r+1}:=\eta^{(1)}$. This gap $\bar{\eta}_{r+1}$
enables us to quantify the distance between the initial objective
value and the optimal objective value of~${\textrm{CP}}_{r+1}$. When
$\bar{\eta}_{r+1} \leq \subscr{\epsilon}{SA}$, the algorithm uses
the estimated subgradient function and $i_{r+1}=0$. Otherwise, the
computational steps can be calculated via convergence of AFWA
for~${\textrm{CP}}_{r+1}$, by
$   \kappa^{i_{r+1}} \bar{\eta}_{r+1} \leq \epsilon_1,
$
where $\kappa$, or using $\bar{\kappa}$ for the data-streaming case,
is determined as in Theorem~\ref{thm:ConvergeJhat}. Let us consider a
threshold value $t_r$
\begin{equation*}
  t_r := \left\{ {\begin{array}{*{20}{l}}
        \epsilon_1, \; & \textrm{if} \; \bar{\eta}_{r} \leq \subscr{\epsilon}{SA}, \\
        \bar{\eta}_{r}, \; & \textrm{o.w.}
\end{array}} \right.
\end{equation*}
Then we can represent each value $i_{r}$ by
$ i_{r}={\log}_{\kappa}(\frac{\epsilon_1}{ t_r  }), \; r \in \{ r_n,\ldots, r_{n+1} \}.
$
Let us denote $t := \max_{r}\{ t_r\}$. Then, the computational steps
for $\{ g_{n}^{r}\}_{r=r_n}^{r_{n+1}}$, $\sum_{r} i_r$, are bounded by
$\sum_{r} i_r \leq \bar{r} {\log}_{\kappa}(\frac{\epsilon_1}{t}).$
Finally, the computational steps to achieve an
$\vect{x}_{n}^{\epsilon_2}$, denoted by $\varphi(n,\bar{r}):=
\phi(n)+\sum_{r} i_r + \bar{r}$, are bounded as
$\varphi(n,\bar{r}) \leq \phi(n)+ \bar{r} \left(
  {\log}_{\kappa}(\frac{\epsilon_1}{t}) +1 \right)$.
Again, one should use $\bar{\phi}(1)$ in the bound in place of
$\phi(n)$ if considering a data-streaming scenario.
\end{IEEEproof}
\subsubsection*{Theorem~\ref{thm:assimilation}}
\begin{IEEEproof}
  The first part of the proof is an application of
  Theorem~\ref{thm:ConvergeJhat} and Theorem~\ref{thm:epsSA}. For
  any data set ${\Xi}_{n}$ and the initial data-driven decision
  $\vect{x}^{(r_{n})}$, by Theorem~\ref{thm:ConvergeJhat} we can
  show $\vect{x}^{(r_{n})}$ to be $\epsilon_1$-proper, via
  finding ${J}_n^{\epsilon_1}(\vect{x}^{(r_{n})})$ such that
  $\mathbf{P}^n({\mathbb{E}_{\prob} [f(\vect{x}^{(r_{n})},\xi)] }
  \leq {J}_{n}^{\epsilon_1}(\vect{x}^{(r_{n})}) +
  \epsilon_1)\geq 1- \beta_{n} $.  Then using
  Theorem~\ref{thm:epsSA}, an $\epsilon_2$-optimal
  $\epsilon_1$-proper data-driven decision
  $\vect{x}^{\epsilon_2}_{n}$ with certificate
  ${J}_n^{\epsilon_1}(\vect{x}_n^{\epsilon_2})$ can be
  achieved. Therefore the performance guarantee~\eqref{eq:epsiperfgua}
  holds for $\vect{x}^{\epsilon_2}_{n}$, i.e.,
  $\mathbf{P}^n({\mathbb{E}_{\prob} [f(\vect{x}_n^{\epsilon_2},\xi)]
  } \leq {J}_{n}^{\epsilon_1}(\vect{x}_n^{\epsilon_2}) +
  \epsilon_1)\geq 1- \beta_{n}$.

  \reviseone{In the following, we show the certificate
    ${J}_{n}^{\epsilon_1}(\vect{x}_n^{\epsilon_2})$ can be upper
    bounded in high probability, for each $n$.}

	  First, let $\vect{x}^{\delta}$ denote the $\delta$-optimal solution
  of~\eqref{eq:P}, i.e., ${\mathbb{E}_{\prob}
    [f(\vect{x}^{\delta},\xi)] } \leq J^{\star} + \delta$. By
  construction of the certificate in the algorithm we have
  ${J}_{n}^{\epsilon_1}(\vect{x}_{n}^{\epsilon_2}) \leq
  {J}_{n}(\vect{x}_{n}^{\epsilon_2}) \leq
  {J}_{n}(\vect{x}_{n}^{\star}) +\epsilon_2 \leq
  {J}_{n}(\vect{x}^{\delta})+\epsilon_2 \leq
  {J}_{n}^{\epsilon_1}(\vect{x}^{\delta}) +\epsilon_1 +\epsilon_2$
  for all $n$, where the first inequality holds because ${J}_{n}$
  is the function that achieves the supreme of Problem~\eqref{eq:cert}
  while ${J}_{n}^{\epsilon_1}(\vect{x}_{n}^{\epsilon_2})$ is
  the objective value for a feasible distribution
  $\mathbb{Q}_n^{\epsilon_1}(\vect{x}^{\epsilon_2})$, the second
  inequality holds because $\vect{x}_{n}^{\epsilon_2}$ is
  $\epsilon_2$-optimal, the third inequality holds because
  $\vect{x}_{n}^{\star}$ is a minimizer of the certificate
  function ${J}_{n}$, the last inequality holds because the
  \CGAaco for certificate generation guarantees the
  existence of ${J}_{n}^{\epsilon_1}(\vect{x}^{\delta})$ such that
  ${J}_{n}(\vect{x}^{\delta}) \leq
  {J}_{n}^{\epsilon_1}(\vect{x}^{\delta}) +\epsilon_1$, with an
  distribution $\mathbb{Q}_n^{\epsilon_1}({\vect{x}}^{\delta})$
  satisfying
  $d_W(\hat{\prob}^n,\mathbb{Q}_n^{\epsilon_1}({\vect{x}}^{\delta}))
  \leq \epsilon(\beta_{n})$.

  Next, we exploit the connection between
  ${J}_{n}^{\epsilon_1}(\vect{x}^{\delta})$ and $J^{\star}$. By
  Assumption~\ref{assump:localconcave} on the concavity of $f$ in
  $\xi$, there exists a constant $\hat{L}>0$ such that $f(x,\xi) \le
  \hat{L}(1+ \Norm{\xi})$ holds for all $x \in \real^d$ and $\xi \in
  \mathcal{Z}$. Then by the dual representation of the Wasserstein
  metric from Kantorovich and Rubinstein~\cite{KLV-RGS:58,PME-DK:17}
  we have
  ${J}_{n}^{\epsilon_1}(\vect{x}^{\delta}):={\mathbb{E}_{\mathbb{Q}_n^{\epsilon_1}({\vect{x}}^{\delta})}
    [f(\vect{x}^{\delta},\xi)]} \leq {\mathbb{E}_{\prob}
    [f(\vect{x}^{\delta},\xi)] } +
  \hat{L}d_W(\prob,\mathbb{Q}_n^{\epsilon_1}({\vect{x}}^{\delta}))$.
  In order to quantify the last term, we apply the triangle
  inequality, which gives us
  $d_W(\prob,\mathbb{Q}_n^{\epsilon_1}({\vect{x}}^{\delta})) \leq
  d_W(\prob,\hat{\prob}^n) +
  d_W(\hat{\prob}^n,\mathbb{Q}_n^{\epsilon_1}({\vect{x}}^{\delta}))$. Then
  by the performance guarantee we have $\mathbf{P}^n\{
  d_W(\prob,\hat{\prob}^n) \leq \epsilon(\beta_{n}) \} \geq 1-
  \beta_{n}$, and by the the way of constructing
  $\mathbb{Q}_n^{\epsilon_1}({\vect{x}}^{\delta})$ we have
  $d_W(\hat{\prob}^n,\mathbb{Q}_n^{\epsilon_1}({\vect{x}}^{\delta}))
  \leq \epsilon(\beta_{n})$. These inequalities result in
  $\mathbf{P}^n\{d_W(\prob,\mathbb{Q}_n^{\epsilon_1}(\vect{x}^{\delta}))
  \leq 2\epsilon(\beta_{n}) \} \geq 1- \beta_{n}$. We use now this
  bound to deal with the last term in the upper bound of
  ${J}_{n}^{\epsilon_1}(\vect{x}^{\delta})$. In particular, we have
  $\mathbf{P}^{n}\{ {J}_{n}^{\epsilon_1}(\vect{x}^{\delta}) \leq
  {\mathbb{E}_{\prob} [f(\vect{x}^{\delta},\xi)] } + 2\hat{L}
  \epsilon(\beta_{n}) \} \geq 1- \beta_{n}$ for all
  $n$. \reviseone{Using the obtained inequality
    ${J}_{n}^{\epsilon_1}(\vect{x}_{n}^{\epsilon_2}) \leq
    {J}_{n}^{\epsilon_1}(\vect{x}^{\delta}) +\epsilon_1 +\epsilon_2$
    and knowing $\delta$ can be arbitrary small, we achieved the goal
    as in~\eqref{eq:qualityofJ}.

    Now, it remains to find an $n_0$, associated with an
    $\epsilon_2$-optimal and $\epsilon_1$-proper data-driven decision
    $\vect{x}_{n_0}^{\epsilon_2}$, such that the almost sure
    guarantee~\eqref{eq:asperf} and bound~\eqref{eq:perfbound} of the
    certificate ${J}_{n_0}^{\epsilon_1}(\vect{x}_{n_0}^{\epsilon_2})$
    can be guaranteed for the termination of the \ODAA as
    $N\rightarrow \infty$. We achieve this by two steps.

    First, we show the almost sure performance guarantee when the data
    set is sufficiently large.} For any time period $n$, the algorithm
 \reviseone{finds}  $\vect{x}^{\epsilon_2}_{n}$ with the performance
  guarantee~\eqref{eq:epsiperfgua}, which can be equivalently written
  as $\mathbf{P}^n({\mathbb{E}_{\prob}
    [f(\vect{x}_n^{\epsilon_2},\xi)]} \geq
  {J}_{n}^{\epsilon_1}(\vect{x}_n^{\epsilon_2}) + \epsilon_1)\leq
  \beta_{n}$. As $\sum_{n=1}^{\infty} \beta_{n} < \infty$, from the
  $\supscr{1}{st}$ Borel-Cantelli Lemma we have that
  $\mathbf{P}^{\infty}\{ {\mathbb{E}_{\prob}
    [f(\vect{x}_n^{\epsilon_2},\xi)] } \geq
  {J}_{n}^{\epsilon_1}(\vect{x}_n^{\epsilon_2}) + \epsilon_1 \;
  \textrm{occurs infinitely many often} \} =0$. That is, almost surely
  we have that ${\mathbb{E}_{\prob} [f(\vect{x}_n^{\epsilon_2},\xi)] }
  \geq {J}_{n}^{\epsilon_1}(\vect{x}_n^{\epsilon_2}) + \epsilon_1$
  occurs at most for finite number of $n$. Thus, there exists a
  sufficiently large $n_1$, such that for all $n \geq n_1$, we have
  ${\mathbb{E}_{\prob} [f(\vect{x}_n^{\epsilon_2},\xi)] } \leq
  {J}_{n}^{\epsilon_1}(\vect{x}_n^{\epsilon_2}) + \epsilon_1$ occurs
  almost surely, i.e., $\mathbf{P}^{n}({\mathbb{E}_{\prob}
    [f(\vect{x}_{n}^{\epsilon_2},\xi)] } \leq
  {J}_{n}^{\epsilon_1}(\vect{x}_{n}^{\epsilon_2}) + \epsilon_1)=1$ for
  all $n \geq n_1$. Later if we pick $n_0 \geq n_1$, then the almost
  sure performance guarantee holds for such
  $\vect{x}_{n_0}^{\epsilon_2}$ and
  ${J}_{n_0}^{\epsilon_1}(\vect{x}_{n_0}^{\epsilon_2})$.

	\reviseone{Second, we show a tight certificate bound can be achieved almost surely. Consider performance bound~\eqref{eq:qualityofJ}. As
  $\epsilon(\beta_{n})$ decreases and goes to $0$ as $n \rightarrow
  \infty$, there exists $n_2$ such that $2\hat{L}\epsilon(\beta_{n})
  \leq \epsilon_3$ holds for all $n \geq n_2$. Therefore, we have
	$\mathbf{P}^{n}\{ {J}_{n}^{\epsilon_1}(\vect{x}_n^{\epsilon_2}) \leq  J^{\star} + \epsilon_1 + \epsilon_2 + \epsilon_3 \} \geq 1- \beta_{n}$ for all $n \geq n_2$, or equivalently, $\mathbf{P}^{n}\{ {J}_{n}^{\epsilon_1}(\vect{x}_n^{\epsilon_2}) \geq  J^{\star} + \epsilon_1 + \epsilon_2 + \epsilon_3 \} \leq \beta_{n}$. As $\sum_{n=1}^{\infty} \beta_{n} < \infty$, then the $\supscr{1}{st}$
  Borel-Cantelli Lemma applies to this situation.
  Thus we claim that there exists a sufficiently large $n_3$ such that
  for all $n \geq \max \{n_2,n_3\}$ we have almost surely, ${J}_{n}^{\epsilon_1}(\vect{x}_n^{\epsilon_2}) \leq  J^{\star} + \epsilon_1 + \epsilon_2 + \epsilon_3$.}

Then, by letting $n_0:= \max \{n_1,n_2,n_3 \}$ we have almost sure
  performance guarantee~\eqref{eq:asperf} and almost surely, the bound~\eqref{eq:perfbound}.
\end{IEEEproof}
\subsubsection*{Lemma~\ref{lemma:closeof2P}}
\begin{IEEEproof}
  The proof is an application of the dual characterization of
  the Wasserstein distance. Let us consider
\begin{equation*}
\begin{split}
  d_W(\hat{\prob}^n,\tilde{\prob}^n)&= \sup_{f \in \mathcal{L}}\{\int_{\mathcal{Z}} {f(\xi)\hat{\prob}^n(d{\xi})} -\int_{\mathcal{Z}} {f(\xi)\tilde{\prob}^n(d{\xi})} \}, \\
  & = \frac{1}{n}\sup_{f \in \mathcal{L}}\{ \sum_{k=1}^{n}
  f({\xi}_{k}) - \sum_{k=1}^{p} \theta_k f(\zeta_{k}) \}.
\end{split}
\end{equation*}
By partitioning the data set ${\Xi}_{n}$ into $\mathcal{C}_n$ and
${\Xi}_{n} \setminus \mathcal{C}_n$ for each summation term, we have
\begin{equation*}
\begin{aligned}
  \sum_{k=1}^{n} f({\xi}_{k})=&\sum_{\varsigma \in \mathcal{C}_n} f(\varsigma) + \sum_{\varsigma \in {\Xi}_{n} \setminus \mathcal{C}_n } f(\varsigma), \\
  \sum_{k=1}^{p} \theta_k f(\zeta_{k})=&\sum_{k=1}^{p} f(\zeta_{k})
  +\sum_{\varsigma \in {\Xi}_{n} \setminus \mathcal{C}_n }
  \ell_{\varsigma}^{-1} \sum_{k \in I_{\varsigma}} f(\zeta_{k}).
\end{aligned}
\end{equation*}
Canceling the first summation term gives us the following
\begin{equation*}
\begin{split}
  d_W(\hat{\prob}^n,\tilde{\prob}^n)&= \frac{1}{n} \sup_{f \in
    \mathcal{L}}\{ \sum_{\varsigma \in {\Xi}_{n} \setminus
    \mathcal{C}_n } \ell_{\varsigma}^{-1} \sum_{k \in I_{\varsigma}}
  f(\varsigma)- f(\zeta_{k})  \},
\\
  & \leq \frac{1}{n} \sup_{f \in \mathcal{L}}\{ \sum_{\varsigma \in
    {\Xi}_{n} \setminus \mathcal{C}_n } \ell_{\varsigma}^{-1} \sum_{k
    \in I_{\varsigma}} |f(\varsigma)- f(\zeta_{k})|  \} ,\\
  & \leq \frac{1}{n} \sum_{\varsigma \in {\Xi}_{n} \setminus \mathcal{C}_n } \ell_{\varsigma}^{-1} \sum_{k \in I_{\varsigma}} \Norm{\varsigma - \zeta_{k}}  ,\\
  & \leq \frac{1}{n} \sum_{\varsigma \in {\Xi}_{n} \setminus \mathcal{C}_n } \ell_{\varsigma}^{-1} \sum_{k \in I_{\varsigma}} \omega = \frac{n-p}{n} \omega \leq \omega , \\
\end{split}
\end{equation*}
where the first inequality is derived taking component-wise absolute
values; the second inequality is due to $f$ being in the space of
Lipschitz functions defined on $\mathcal{Z}$ with Lipschitz constant
1; and the third inequality is due to~$\varsigma \in B_{\omega}(\zeta_{k})$.
\end{IEEEproof}
\subsubsection*{Lemma~\ref{lemma:tildePcertgener}}
\begin{IEEEproof}
  From Lemma~\ref{lemma:certgener} we have $\mathbf{P}^n\{
  d_W(\prob,\hat{\prob}^n) \leq \epsilon(\beta_{n}) \} \geq
  1-\beta_{n} $ for each $n$. Then using the result from
  Lemma~\ref{lemma:closeof2P} we have $\mathbf{P}^n\{
  d_W(\prob,\tilde{\prob}^n) \leq d_W(\tilde{\prob}^n,\hat{\prob}^n) +
  d_W(\prob,\hat{\prob}^n) \leq \epsilon(\beta_{n}) +\omega \} \geq
  1-\beta_{n}$, i.e., $\mathbf{P}^n\{ d_W(\prob,\tilde{\prob}^n) \leq
  \tilde{\epsilon}(\beta_{n}) \} \geq 1-\beta_{n}$ for each $n$. The
  rest of the proof follows directly from that in
  Lemma~\ref{lemma:certgener} and Theorem~\ref{thm:cvxredu}.
\end{IEEEproof}

\bibliographystyle{IEEEtran}
\bibliography{../../../../bib/alias,../../../../bib/SMD-add,../../../../bib/JC,../../../../bib/SM}

\begin{IEEEbiography}[{\includegraphics[width=1in,height=1.25in,clip,keepaspectratio]{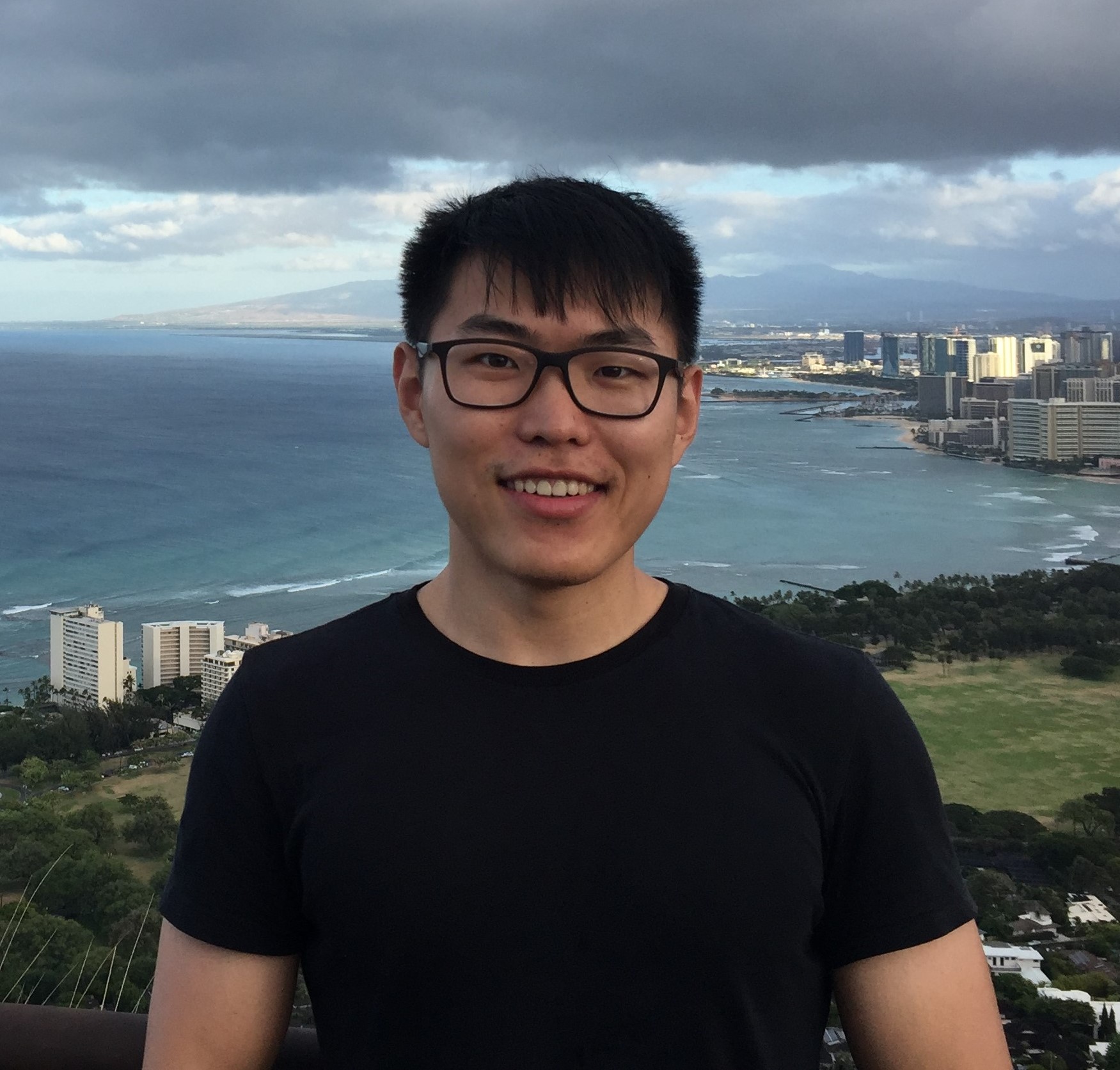}}]{Dan
    Li} received the B.E. degree in automation from the Zhejiang
  University, Hangzhou, China, in 2013, the M.Sc. degree in chemical
  engineering from Queen's University, Kingston, Canada, in 2016.  He
  is currently a Ph.D. student at University of California, San Diego,
  CA, USA. His current research interests include data-driven systems
  and optimization, dynamical systems and control, optimization
  algorithms, applied computational methods, and stochastic systems. He
  received Outstanding Student Award from Zhejiang University in 2012,
  Graduate Student Award from Queen's University in 2014, and
  Fellowship Award from University of California, San Diego, in 2016.
\end{IEEEbiography}
\begin{IEEEbiography}[{\includegraphics[width=1in,height=1.25in,clip,keepaspectratio]{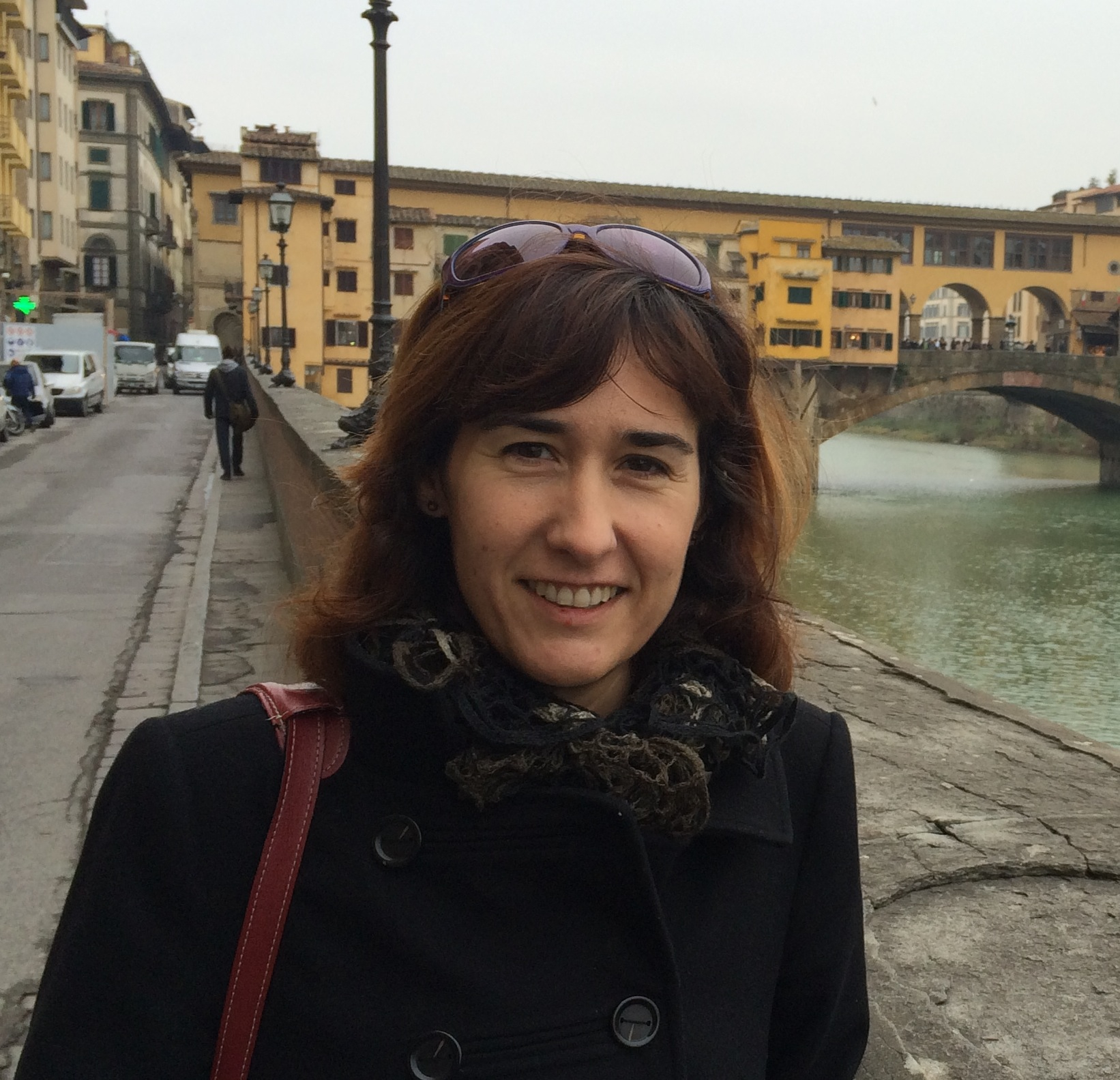}}]{Sonia
   Martínez} is a Professor at the Department of Mechanical and
 Aerospace Engineering at the University of California, San
 Diego. She received her Ph.D. degree in Engineering Mathematics from
 the Universidad Carlos III de Madrid, Spain, in May 2002. Following
 a year as a Visiting Assistant Professor of Applied Mathematics at
 the Technical University of Catalonia, Spain, she obtained a
 Postdoctoral Fulbright Fellowship and held appointments at the
 Coordinated Science Laboratory of the University of Illinois,
 Urbana-Champaign during 2004, and at the Center for Control,
 Dynamical systems and Computation (CCDC) of the University of
 California, Santa Barbara during 2005.

 Her research interests include networked control systems,
 multi-agent systems, and nonlinear control theory with applications
 to robotics and cyber-physical systems. For her work on the control
 of underactuated mechanical systems she received the Best Student
 Paper award at the 2002 IEEE Conference on Decision and Control. She
 co-authored with Jorge Cortés and Francesco Bullo "Motion
 coordination with Distributed Information" for which they received
 the 2008 Control Systems Magazine Outstanding Paper Award. She is a
 Senior Editor of the IEEE Transactions on Control of Networked
 Systems and an IEEE Fellow.
\end{IEEEbiography}

\end{document}